\documentclass[10pt]{article}
\usepackage{amsmath,amsthm,amssymb,dsfont,graphicx,xspace,epsfig,xcolor}
\usepackage[plain]{fullpage}
\usepackage{thm-restate}
\usepackage{hyperref}

\usepackage{algorithm}\usepackage{algorithmic}
\usepackage{tikz, tkz-graph, tkz-berge}
\usetikzlibrary{decorations.pathreplacing}
\usetikzlibrary{patterns}
\usetikzlibrary{patterns.meta}
\usepackage{color}
\usepackage{comment}
\usepackage{authblk}
\usepackage[shortlabels]{enumitem}
\usepackage{float}

\definecolor{g-brown}{RGB}{128,0, 0}

\newtheorem{theorem}{Theorem}[section]
\newtheorem{corollary}[theorem]{Corollary}
\newtheorem{proposition}[theorem]{Proposition}
\newtheorem{lemma}[theorem]{Lemma}
\newtheorem{claim}{Claim}[theorem]

\theoremstyle{definition}

\newcommand{\ind}[1]{[#1]}
\newcommand{\defproblem}[3]{
\renewcommand{\arraystretch}{1}
 \vspace{3mm}
\noindent\fbox{
 \begin{minipage}{0.96\textwidth}
 \begin{tabular*}{\textwidth}{@{\extracolsep{\fill}}lr} #1  \\ \end{tabular*}
 {\bf{Input:}} #2 \\
 {\bf{Question:}} #3
 \end{minipage}
 }
 \vspace{3mm}
}

\newcommand{\true}{\mathrm{\tt true}}
\newcommand{\false}{\mathrm{\tt false}}

\DeclareMathOperator{\UG}{UG}

\newenvironment{proofclaim}[1][]
	{\par\noindent {\it Proof of claim}. }{ \hfill$\lozenge$\par\vspace{10pt}}

\title{Complexity results on the decomposition of a digraph into directed linear forests and out-stars}

\author[,1]{Florian H\"orsch\thanks{Research supported by the European Research Council (ERC) consolidator grant No.~725978 SYSTEMATICGRAPH.}}
\author[,2]{Lucas Picasarri-Arrieta\thanks{Research supported by Japan Science and Technology Agency (JST) as part of Adopting Sustainable Partnerships for Innovative Research Ecosystem (ASPIRE), Grant Number JPMJAP2302.}}

\affil[1]{CISPA Saarbrücken, Germany}
\affil[2]{National Institute of Informatics, Tokyo, Japan}

\date{}

\begin{document}
\maketitle

\begin{abstract}
We consider two decomposition problems in directed graphs. We say that a digraph is $k$-bounded for some $k \in \mathbb{Z}_{\geq 1}$ if each of its connected components contains at most $k$ arcs. 

For the first problem, a directed linear forest is a collection of vertex-disjoint directed paths and we consider the problem of decomposing a given digraph into a $k$-bounded and an $\ell$-bounded directed linear forest for some fixed $k,\ell \in \mathbb{Z}_{\geq 1}\cup \{\infty\}$. We give a full dichotomy for this problem by showing that it can be solved in polynomial time if $k+\ell \leq 3$ and is NP-complete otherwise. This answers a question of Campbell, Hörsch, and Moore.

For the second problem, we say that an out-galaxy is a vertex-disjoint collection of out-stars. Again, we give a full dichotomy of when a given digraph can be arc-decomposed into a $k$-bounded and an $\ell$-bounded out-galaxy for fixed $k,\ell \in \mathbb{Z}_{\geq 1}\cup \{\infty\}$. More precisely, we show that the problem can be solved in polynomial time if $\min\{k,\ell\}\in \{1,\infty\}$ and is NP-complete otherwise.

\medskip

\noindent{}{\bf Keywords:} Digraphs, arc-partitions, directed arboricity, branchings, directed linear forests, directed stars, complexity.
\end{abstract}

\bibliographystyle{abbrv}

\section{Introduction}

Given a graph $G$ (resp. a digraph), we say that
a collection $(H_1, \dots , H_t)$ of spanning subgraphs (resp. subdigraphs) of $G$ is a {\bf decomposition} of $G$ if
$(A(H_1),\dots, A(H_t))$ is a partition of $E(G)$ (resp. $A(G)$).
Given two properties $\mathcal{P}_1$ and $\mathcal{P}_2$ on digraphs, a general problem consists of asking whether a digraph decomposes into two subdigraphs $H_1,H_2$ such that $H_i$ satisfies $\mathcal{P}_i$. We refer the reader to~\cite{bang2022complexity} for an extensive study of such problems. 

Given an undirected graph $G$, a classical decomposition problem asks for the minimum integer $t$ such that $G$ admits a decomposition $(H_1,\dots,H_t)$ in which $H_i$ is a forest for $i \in [t]$. This parameter is called the {\bf arboricity} of $G$, and the following celebrated theorem by Nash-Williams characterizes graphs with arboricity at most $t$.
\begin{theorem}[Nash-Williams~\cite{NashWilliams1964}]
  \label{thm:NW}
  For some positive integer $t$, a graph $G$ decomposes into $t$ forests if and only if for every $X\subseteq V(G)$, the subgraph of $G$ induced by $X$ contains at most $t\cdot (|X|-1)$ edges.
\end{theorem}

Although the condition in Theorem~\ref{thm:NW} involves an exponential family of vertex sets, it is well-known that one can actually compute the arboricity of a graph in polynomial time.

An analogous problem has been defined for digraphs.
A digraph is a {\bf branching} if it is an orientation of a forest in which every vertex has in-degree at most one. 
The {\bf directed arboricity} of a digraph $D$ is then the minimum integer $t$ for which $D$ decomposes into $t$ branchings.
Using the celebrated Edmonds' Branching Theorem~\cite{edmonds1973}, Frank proved the following characterization of digraphs with directed arboricity at most $t$.

\begin{theorem}[Frank~\cite{frank1979covering}]
  \label{thm:frank}
  A digraph $D$ decomposes into $t$ branchings if and only if
  \begin{itemize}
      \item $d_D^-(v) \leq t$ holds for every vertex $v\in V(D)$, and
      \item the underlying graph $G$ of $D$ decomposes into $t$ forests.
  \end{itemize}
\end{theorem}

Theorem~\ref{thm:frank} together with the undirected result implies that the directed arboricity of a digraph can be computed in polynomial time.
A natural question is then to look at the notion of arboricity when restricted to some classes of branchings. This is the topic of this work, and we will consider both (bounded) directed linear forests and (bounded) out-galaxies. 

Recall that a {\bf directed  linear forest} is a vertex-disjoint union of directed paths. A directed linear forest is {\bf $k$-bounded} for some $k \in \mathbb{Z}_{\geq 1} \cup \{\infty\}$ if each of its connected components is a directed path of length at most $k$, where the length of a directed path refers to its number of arcs.
Given a digraph $D$, the minimum integer $t$ such that $D$ decomposes into $t$ directed linear forests is known as the {\bf directed linear-arboricity} of $D$, introduced by Nakayama and Peroche~\cite{nakayama1987linear}, see also~\cite{he2017linear}. For every integer $k\geq 1$, when restricted to $k$-bounded directed linear forests, following the notion recently introduced by Zhou et al.~\cite{zhou2022linear}, we define the {\bf directed linear-$k$-arboricity} of a given digraph $D$ to be the minimum integer $t$ such that $D$ decomposes into $t$ $k$-bounded directed linear forests.
Clearly, the problem of decomposing into one $k$-bounded linear forest is trivial. Further, the problem of decomposing into three or more directed linear forests generalizes the 3-edge-colouring problem, so there is little hope to obtain positive algorithmic results for this class of problems. We hence focus on the possibility of decomposing a digraph into two bounded directed linear forests.
We thus consider the following class of problems, for which not much is known.

\defproblem
{$(k,\ell)$-bounded directed linear forest decomposition ($(k,\ell)$-BDLFD)}
{A digraph $D$}
{Does $D$ decompose into a $k$-bounded directed linear forest and an $\ell$-bounded directed linear forest?}

In the following result, we settle the complexity of the problem above for all fixed values of $k$ and $\ell$, which answers a question raised by Campbell, Moore, and the first author~\cite[Problem~3]{campbell2023decompositions}.

\begin{theorem}\label{maindblfd}
    \label{thm:main_1}
    For $k,\ell \in \mathbb{Z}_{\geq 1}\cup \{\infty\}$, the $(k,\ell)$-BDLFD problem is solvable in polynomial time when $k + \ell \leq 3$ and it is NP-complete when $k+\ell\geq 4$.
\end{theorem}

Observe that it is straightforward that the (1,1)-BDLFD problem is solvable in polynomial time: it asks whether the underlying graph of the input digraph is $2$-edge-colourable. The $(\infty,\infty)$-BDLFD problem has been shown to be NP-complete by Nakayama and Peroche in~\cite[Theorem~3.2]{nakayama1987linear}.
We prove the remaining cases in Section~\ref{sec:decomposition_dlf}.
Interestingly, the dichotomy proved in Theorem~\ref{maindblfd} matches exactly the one proved for the corresponding problem in undirected graph which was partially proved in~\cite{campbell2023decompositions} and partially proved by Banerjee et al. in \cite{banerjee2023complexity}. Notice also that the following result directly follows from Theorem~\ref{thm:main_1} when $k=\ell$.
\begin{corollary}
    \label{cor:linear_k_arboricity}
    For every integer $k\geq 2$, deciding whether the directed linear-$k$-arboricity of a digraph is $2$ is NP-complete.
\end{corollary}

\medskip

An {\bf out-star} is a branching in which at most one vertex has positive out-degree. An {\bf out-galaxy} is a collection of vertex-disjoint out-stars. Again, it is {\bf $k$-bounded} for some $k\in \mathbb{Z}_{\geq 1}\cup \{\infty\}$ if each of its connected component contains at most $k+1$ vertices.
Given a digraph $D$, the minimum integer $t$ such that $D$ decomposes into $t$ out-galaxies is known as the {\bf directed star-arboricity} of $D$, and has been introduced by Algor and Alon in ~\cite{algor1989star}, see also~\cite{amini2010wdm,gonccalves2012spanning,du2017bounds}.
Analogously to the linear-arboricity, for every integer $k\geq 1$, when restricted to $k$-bounded out-galaxies, we obtain the notion of {\bf directed star-$k$-arboricity}. To the best of the authors' knowledge, it has not been defined before.
We then consider the following class of problems, which includes a directed analogue of~\cite[Problem 2]{campbell2023decompositions}. The undirected version of the problem was solved in \cite[Theorem 1.4]{banerjee2023complexity}.

\defproblem
{$(k,\ell)$-bounded out-galaxy decomposition ($(k,\ell)$-BOGD)}
{A digraph $D$}
{Does $D$ decompose into a $k$-bounded out-galaxy and an $\ell$-bounded out-galaxy?}

In the following result, we settle the complexity of the problem above for all fixed values of $k$ and $\ell$.
\begin{theorem}\label{bogdmain}
    $(k,\ell)$-BOGD is solvable in polynomial time when $\min\{k,\ell\}=1$  or $k = \ell = \infty$. In every other case, it is NP-complete.
\end{theorem}

Again, we can observe that (1,1)-BOGD is easily solvable in polynomial time. We prove the other cases in Section~\ref{sec:galaxy}. 
Observe the following analogue of Corollary~\ref{cor:linear_k_arboricity} for star-arboricity.
\begin{corollary}
    For every fixed $k\geq 2$, deciding if the directed star-$k$-arboricity of a digraph is $2$ is NP-complete.
\end{corollary}

\paragraph{Outline of the paper.}
The paper is organised as follows. 
We first give our notation on digraphs and recall some preliminary results in Section~\ref{sec:preliminaries}.
Section~\ref{sec:decomposition_dlf} is devoted to the proof of Theorem~\ref{maindblfd}. We start by proving, in Section~\ref{sec:1_2_case}, that $(k,\ell)$-BDLFD reduces to $2$-SAT when $k+\ell \leq 3$, using some ideas which are similar to the ones of~\cite[Theorem~2]{campbell2023decompositions}. When both $k$ and $\ell$ are integers, we prove the NP-completeness of $(k,\ell)$-BDLFD when $\ell = 1$ and $k\geq 3$ in Section~\ref{sec:1_3_case} and its NP-completeness when $\min \{k,\ell\} \geq 2$ in Section~\ref{sec:2_2_case}. In both cases, our proof consists of building specific gadgets, called forcers, which force the existence, in every decomposition, of a directed path of specific length on a vertex. We then use these forcers to build larger gadgets that modelize variables and clauses, hence allowing us to reduce satisfiability problems to $(k,\ell)$-BDLFD. Finally, we prove the NP-completeness of $(k,\ell)$-BDLFD when one of $\{k,\ell\}$ is equal to $\infty$ by reducing from hamiltonicity in $2$-diregular digraphs.
Section~\ref{sec:galaxy} is devoted to the proof of Theorem~\ref{bogdmain}. We first show in Section~\ref{out1} that $(\infty,\infty)$-BOGD reduces to a $2$-colourability problem, which is pretty simple. In Section~\ref{sec:galaxy_1_k}, we give a more involved proof that $(1,k)$-BOGD is solvable in polynomial time when $k\in\mathbb{Z}_{\geq 2} \cup \{\infty\}$. The proof is based on a reduction to a matching problem in undirected graphsand is also inspired by the proof of ~\cite[Theorem~2]{campbell2023decompositions}. We finally prove the NP-completeness of $(k,\ell)$-BOGD for the remaining cases in Section~\ref{sec:galaxy_2_2} by reducing from specific satisfiability problems.

\section{Preliminaries}
\label{sec:preliminaries}

\subsection{Notation on digraphs}

Our notation follows~\cite{bang2009}. 
For some integer $k$, we denote by $\mathbb{Z}_{\geq k}$ the set of integers which are not smaller than $k$. We use $[k]$ for $\{1,\ldots,k\}$.
Let $D$ be a digraph. The {\bf underlying graph} $\UG(D)$ of $D$ is the graph obtained from $D$ by removing the orientations. We also say that $D$ is an {\bf orientation} of $\UG(D)$. Note that $\UG(D)$ may contain parallel edges. A digraph is {\bf connected} if its underlying graph is a connected graph.  A {\bf connected component} of $D$ is a subdigraph $H$ of $D$ such that $\UG(H)$ is a connected component of $\UG(D)$. 
A {\bf digon} in $D$ is a pair of opposite arcs between two vertices. If $D$ is a digraph without any digon.
 The {\bf out-degree}, $d_D^+(v)$ (resp. {\bf in-degree}, $d_D^-(v)$) of a vertex $v\in V(D)$ is the number of arcs in $A(D)$ of the form $vw$ (resp, $uv$) and the {\bf degree} of $v$ is $d_D(v)=d^+_D(v)+d_D^-(v)$. A digraph $D$ is {\bf $k$-diregular} if every vertex has in-degree and out-degree equal to $k$.
Given a set $X\subseteq V(D)$ of vertices, we denote by $\delta_D^-(X)$ the set of arcs with tail in $V(D)- X$ and head in $X$. We let $\delta_D^+(X)$ denote $\delta_D^-(V(D)- X)$. When $X=\{x\}$ for a single vertex $x$, with a slight abuse of notation we denote $\delta_D^+(X)$ and $\delta_D^-(X)$ respectively by $\delta_D^+(x)$ and $\delta_D^-(x)$.
 A {\bf matching} is a graph in which every vertex has degree at most one. A {\bf directed matching} is an orientation of a matching. 
A {\bf directed cycle} is a connected digraph in which every vertex has in-degree and out-degree exactly $1$.
A {\bf directed path} is obtained from a directed cycle by the removal of exactly one arc. Let $P$ be an orientation of a path. The {\bf endvertices} of $P$ are the vertices of degree one in $\UG(P)$.  When $P$ has length at least two, the {\bf endarcs} of $P$ are the arcs incident to its endvertices.
A connected branching is called an {\bf arborescence}, and its unique vertex with in-degree $0$ is called its {\bf root}.

We now recursively define a graph called {\bf binary tree} with a special vertex called its {\bf tip} that has an integer $k \geq 0$ as parameter called its {\bf depth}. A binary tree of depth 0 only consists of the tip. For $k \geq 1$, a binary tree of depth $k$ is constructed from two disjoint binary trees $T_1,T_2$ of depth $k-1$ with tips $x_1,x_2$, respectively, by adding a vertex $x$ and the edges $x_1x$ and $x_2x$. Further, $x$ is the tip of the binary tree.

\subsection{Known complexity results}
We recall a collection of well-known complexity results that we will use all along the paper.

\begin{proposition}(see~\cite{garey1979})
    \label{prop:3SAT}
    $k$-SAT is NP-complete for any $k \geq 3$.
\end{proposition}

\begin{proposition}(see~\cite{garey1979})
    \label{prop:2SAT}
    $2$-SAT can be solved in polynomial time. Moreover, a satisfying assignment for a positive instance can be found in polynomial time.
\end{proposition}

\begin{proposition}(see~\cite{garey1979})
    \label{checkbip}
    We can check in polynomial time whether a given graph is bipartite.
\end{proposition}

\begin{proposition}(see~\cite{schrijver2003combinatorial})
    \label{checkmatch}
    Let $G$ be a graph and $Z \subseteq V(G)$. We can decide in polynomial time whether $G$ contains a matching covering $Z$. Moreover, if such a matching exists, it can be computed in polynomial time.
\end{proposition}

A {\bf hamiltonian cycle} in a digraph $D$ is a directed cycle of length $|V(D)|$.

\begin{theorem}[{\cite[Theorem~6.1.2]{bang2009}}]
    \label{thm:hardness_hamiltonian}
    It is NP-complete to decide whether a 2-diregular digraph contains a hamiltonian cycle.
\end{theorem}

For some integer $k \geq 1$, an instance of the {\sc Monotone Equitable $k$-SAT} problem (ME-$k$-SAT for conciseness) consists of a set of variables $X$ and a set of clauses $C$ each of which contains exactly $2k + 1$ non-negated variables and the question is whether there is a truth assignment $\phi : X \xrightarrow{} \{\true, \false\}$ such that every
clause in $C$ contains at least $k$ $\true$ and $k$ $\false$ variables with respect to $\phi$. Note that {\sc ME-$1$-SAT} is often referred to as {\sc Monotone-Not-All-Equal-3-SAT}.

\begin{proposition}[{\cite[Proposition~4.1]{DHHRarxiv23}}]
    \label{prop:MEkSAT_NP_hard}
    For every integer $k\geq 1$, 
    {\sc ME-$k$-SAT} is NP-complete.
\end{proposition}

The {\sc $(3,B2)$-SAT} problem is the restriction of $3$-SAT in which every litteral appears exactly twice.

\begin{theorem}[{\cite[Theorem~1]{berman2004approximation}}]
    \label{thm:3_B2_SAT_NPc}
    {\sc $(3,B2)$-SAT} is NP-complete.
\end{theorem}

\section{Decompositions into bounded directed linear forests}
\label{sec:decomposition_dlf}

This section is dedicated to proving our results concerning the computational complexity of decomposing a given digraph into two (bounded) directed linear forests. More precisely, we prove Theorem~\ref{maindblfd}. 

First, in Section~\ref{sec:1_2_case}, we prove the positive algorithmic result contained in Theorem~\ref{maindblfd}, namely the case that $k+\ell \leq 3$. The negative results are split into several parts. In Section~\ref{sec:1_3_case}, we prove the complexity result for the case when both $k$ and $\ell$ are finite integers such that $\ell=1$ and $k\geq 3$. In Section~\ref{sec:2_2_case}, we consider the case that both $k$ and $\ell$ are finite integers and $\min\{k,\ell\}\geq 2$. 
Finally, the case that exactly one of the directed linear forests is unbounded is considered in Section~\ref{sec:unbounded_case}. Together with the result in~\cite{nakayama1987linear} and the fact that $(1,1)$-BDLFD reduces to 2-colourability, we obtain Theorem~\ref{maindblfd}. 

Throughout this section, for $k,\ell \in \mathbb{Z}_{\geq 1} \cup \{\infty\}$, a {\bf $(k,\ell)$-decomposition} of a digraph $D$ is a decomposition $(F_k,F_\ell)$ of $D$ such that $F_k$ is a $k$-bounded directed linear forest, and $F_\ell$ is an $\ell$-bounded directed linear forest.

\subsection{Decomposing into a matching and a 2-bounded directed linear forest}
\label{sec:1_2_case}

This section is dedicated to proving the main positive algorithmic result on directed linear forest decompositions. More precisely, we prove the following theorem. 
\begin{theorem}\label{21dblfd}
    (2,1)-BDLFD is solvable in polynomial time.
\end{theorem}

The proof of Theorem~\ref{21dblfd} contains some ideas which are similar to the proof of~\cite[Theorem~2]{campbell2023decompositions}. We first need a collection of easy preliminary results, that we prove for completeness, which deal with decompositions of orientations of paths and cycles.

\begin{proposition}\label{chembound}
    Let $P$ be an orientation of a path of length at least 2, let $a_1,a_2$ be the endarcs of $P$, and let $A_0\subseteq \{a_1,a_2\}$. Then we can decide in polynomial time whether there exists a $(2,1)$-decomposition $(F_2,F_1)$ of $P$ with $\{a_1,a_2\}\cap A(F_1)=A_0$. Further, such a decomposition can be constructed in polynomial time if it exists.
\end{proposition}
\begin{proof}
    If the length of $\UG(P)$ is exactly 2, we can solve the problem by a brute force approach. We may hence suppose that $\UG(P)$ is a path $v_1\ldots v_q$ with $q \geq 4$.  Let $P'=P-v_1$ and let $a_1'$ be the arc in $A(P)$ whose endvertices are $v_2$ and $v_3$. We distinguish three cases.

    \begin{description}
        \item[Case 1:] $A_0 \neq  \emptyset$.

        Assume without loss of generality that $a_1 \in A_0$. Then there exists a $(2,1)$-decomposition $(F_2,F_1)$ of $P$ with $\{a_1,a_2\}\cap A(F_1)=A_0$ if and only if there exists a $(2,1)$-decomposition $(F'_2,F'_1)$ of $P'$ with $\{a_1',a_2\}\cap A(F'_1)=A_0-a_1$. We can hence recursively solve this smaller problem. 
        
        \item[Case 2:] \textit{$A_0 =  \emptyset$ and $P\ind{\{v_1,v_2,v_3\}}$ is not a directed path.}

        In this case, there exists a $(2,1)$-decomposition $(F_2,F_1)$ of $D$ with $\{a_1,a_2\}\cap A(F_1)=A_0$ if and only if there exists a $(2,1)$-decomposition $(F'_2,F'_1)$ of $P'$ with $\{a_1,a_2\}\cap A(F'_1)=A_0 \cup a'_1$. We can hence recursively solve this smaller problem.

        \item[Case 3:] \textit{$A_0 =  \emptyset$ and $P\ind{\{v_1,v_2,v_3\}}$ is a directed path.}

        In this case, $P$ has the desired decomposition. If $q$ is even, we define $(F_2,F_1)$ so that $A(F_2)$ contains the orientation of $v_iv_{i+1}$ for all odd $i \in [q-1]$ and $A(F_1)=A(P)-A(F_2)$. If $q$ is odd, we define $(F_2,F_1)$ so that $A(F_2)$ contains the orientation of the edge $v_1v_2$ and the orientation of the edge $v_iv_{i+1}$ for all even $i \in [q-1]$ and $A(F_1)=A(P)-A(F_2)$. In either case, we have that $(F_2,F_1)$ is a $(2,1)$-decomposition of $P$.
        \qedhere
    \end{description}
\end{proof}

\begin{proposition}\label{chembound2}
    Let $P$ be an orientation of a path of length at least 2, let $a_1,a_2$ be the endarcs of $P$, and let $A_0\subseteq \{a_1,a_2\}$. Then we can decide in polynomial time whether there exists a $(2,1)$-decomposition $(F_2,F_1)$ of $P$ with $A(F_1)\cap \{a_1,a_2\}=A_0$ and every $a \in \{a_1,a_2\}-A_0$ is the only arc in a connected component of $F_2$. Further, such a decomposition can be constructed in polynomial time if it exists.
\end{proposition}
\begin{proof}
If the length of $\UG(P)$ is at most $3$, we can solve the problem by a brute force approach. We may hence suppose that $\UG(P)$ is a path $v_1\ldots v_q$ where $q\geq 5$.
Let $P'$ be the unique connected component of $P-(\{a_1,a_2\}-A_0)$ with $v_2 \in V(P')$ and let $a_1',a_2'$ be the endarcs of $P'$ such that $a'_1$ is incident to $v_2$ and $a'_2$ is incident to $v_{q-1}$. Then there exists a $(2,1)$-decomposition $(F_2,F_1)$ of $P$ with $A(F_1)\cap \{a_1,a_2\}=A_0$ and every $a \in \{a_1,a_2\}-A_0$ is the only arc in a connected component of $F_2$ if and only if there exists a $(2,1)$-decomposition $(F'_2,F'_1)$ of $P'$ with $ \{a'_1,a'_2\}\subseteq A(F'_1)$. By Proposition~\ref{chembound}, in polynomial time, we can decide the existence of such a decomposition of $P'$ and find such a decomposition if it exists in polynomial time. Clearly, this yields a constructive polynomial time algorithm for finding the desired decomposition of $P$. 
\end{proof}
\begin{proposition}\label{cycchem2}
    Let $P$ be an orientation of a path and let $a$ be an endarc  of $P$. Then there exists a $(2,1)$-decomposition $(F_2,F_1)$ of $P$ with $a\in A(F_1)$ and a $(2,1)$-decomposition $(F'_2,F'_1)$ of $P$ with $a\in A(F'_2)$ and $a$ is the only arc in a connected component of $F'_2$. Further, these decompositions can be constructed in polynomial time.
\end{proposition}
\begin{proof}
    Let $\UG(P)$ be a path $v_1 \ldots v_q$ such that $a$ is the orientation of $v_1v_2$. Let $(F_2,F_1)$ be the decomposition of $P$ defined so that $A(F_1)$ contains the orientation of $v_iv_{i+1}$ for all odd $i \in [q-1]$ and $A(F_2)=A(P)-A(F_1)$. Further, let $(F_2',F_1')$ be defined by $A(F_2')=A(F_1)$ and $A(F_1')=A(F_2)$. It is easy to see that $(F_2,F_1)$ and $(F_2',F_1')$ have the desired properties. The proof is clearly algorithmic.
    \end{proof}
    \begin{proposition}\label{cycchem}
        Let $C$ be an orientation of a cycle $v_1\ldots v_q v_1$. Then $C$ admits a $(2,1)$-decomposition. Further, such a decomposition can be found in polynomial time.
    \end{proposition}
    \begin{proof}
    If $q$ is even, let the decomposition $(F_2,F_1)$ of $C$ be defined so that the orientation of $v_iv_{i+1}$ is contained in $A(F_2)$ for all odd $i \in [q-1]$  and $A(F_1)=A(C)-A(F_2)$. It is easy to see that $(F_2,F_1)$ is a $(2,1)$-decomposition of $C$. Now suppose that $q$ is odd. As $\sum_{v \in V(C)}d_C^+(v)=|A(C)|$ is odd, there exists some $v \in V(C)$ with $d_C^+(v)=1$. By symmetry, we may suppose that $C\ind{\{v_1,v_2,v_3\}}$ is a directed path $v_1v_2v_3$. We now consider the decomposition $(F_2,F_1)$ of $C$ which is defined so that $A(F_2)$ contains $v_1v_2$ and $v_iv_{i+1}$ for all even $i \in [q-1]$ and $A(F_1)=A(C)-A(F_2)$. It is easy to see that $(F_2,F_1)$ is a $(2,1)$-decomposition of $C$.
    Observe that this proof is clearly algorithmic.   
\end{proof}

We are now ready to give the main proof of Theorem~\ref{21dblfd}. It is based on excluding a collection of configurations which clearly exclude the desired decomposition and then reducing the problem to 2-SAT.

\begin{proof}[Proof of Theorem~\ref{21dblfd}]
  Let $D$ be an instance of $(2,1)$-BDLFD. Clearly, if there exists some $v \in V(D)$ with $d_{D}(v)\geq 4$, then $D$ is a negative instance of $(2,1)$-BDLFD. We may hence suppose that $d_{D}(v)\leq 3$  for all $v \in V(D)$. Let $V_3$ denote the set of vertices $v \in V(D)$ with $d_{D}(v)=3$ and let $V_1$ denote the set of vertices $v \in V$ with $d_{D}(v)=1$. Observe that, if there is some $v \in V_3$ with $\max\{d_D^+(v),d_D^-(v)\}=3$, then $D$ is a negative instance of $(2,1)$-BDLFD. We may hence suppose that $\max\{d_D^+(v),d_D^-(v)\}=2$  for every $v \in V_3$. By Propositions~\ref{cycchem2} and~\ref{cycchem}, we may suppose that every connected component of $D$ contains at least one vertex of $V_3$. Let $\mathcal{P}_0$ be the collection of subdigraphs of $D$ whose underlying graph is either a path connecting two vertices in $V_3$ and none of whose interior vertices is contained in $V_3$ or a cycle containing exactly one vertex in $V_3$. For some $P \in \mathcal{P}_0$, we define the {\bf associated path $P'$} to be $P$ if $P$ is a path and to be the path obtained from $P$ by detaching the unique vertex in $V(P)\cap V_3$ into two vertices of degree 1 in $\UG(P')$ if $P$ is cycle. Further, let $\mathcal{P}_1$ be the collection of subdigraphs of $D$ whose underlying graph is a path connecting a vertex in $V_3$ and a vertex in $V_1$ and none of whose interior vertices is contained in $V_3$. Observe that $\{A(P) \mid P \in \mathcal{P}_0 \cup \mathcal{P}_1\}$ is a partition of $A(D)$. We now create an instance $(X,\mathcal{C})$ of 2-SAT. We let $X$ consist of a variable $x_a$ for every arc $a \in A$ which is incident to at least one vertex in $V_3$ in $D$. We let $\mathcal{C}$ consist of a collection of clauses $\mathcal{C}_v$ for every $ v\in V_3$, a collection of clauses $\mathcal{C}_P$ for every $P \in \mathcal{P}_0$  and a collection of clauses $\mathcal{C}_Q$ for every subdigraph $Q \in \mathcal{Q}$ where $\mathcal{Q}$ denotes the collection of subdigraphs of $D$ which are isomorphic to a digon. First consider some $v \in V_3$. Let $a_1,a_2,a_3$ be the arcs incident to $v$ such that $|\delta_D^+(v)\cap \{a_2,a_3\}| \in \{0,2\}$. We then let $\mathcal{C}_v$ consist of the clauses $\{\bar{x}_{a_1}\},\{\bar{x}_{a_2}, \bar{x}_{a_3}\}$, and $\{x_{a_2}, x_{a_3}\}$. Now consider some $P \in \mathcal{P}_0$ and let $P'$  be the associated path of $P$. If $A(P)$ contains a single arc $a$, we set $\mathcal{C}_P=\{x_a\}$. Now suppose that the length of $\UG(P')$ is at least 2 and let $a_1$ and $a_2$ be the endarcs of $P'$.
  For every $A_0 \subseteq \{a_1,a_2\}$, we test whether there is a $(2,1)$-decomposition $(F_2^P,F_1^P)$ of $P'$ such that $A(F_1^P)\cap \{a_1,a_2\}=A_0$ and every $a \in \{a_1,a_2\}-A_0$ is the only arc in a connected component of $F_2^P$.
  Observe that this can be tested in polynomial time by Proposition~\ref{chembound2}. If no such decomposition exists, we add the clause $\{y_{a_1},y_{a_2}\}$ to $\mathcal{C}_P$ where $y_{a_i}=\bar{x}_{a_i}$ if $a_i \in A_0$ and $y_{a_i}=x_{a_i}$ otherwise, for $i \in [2]$. We do this for every $A_0\subseteq \{a_1,a_2\}$, thus creating $\mathcal{C}_P$. Further, for every subdigraph $Q$ of $D$ which is a digon containing two arcs $a_1,a_2$, we set $\mathcal{C}_q=\{\{x_{a_1},x_{a_2}\},\{\bar{x}_{a_1},\bar{x}_{a_2}\}\}$. Finally, we set $\mathcal{C}=\bigcup_{v \in V_3 \cup \mathcal{P}_0 \cup \mathcal{Q}}\mathcal{C}_v$.
   This finishes the description of $(X,\mathcal{C})$.
  \begin{claim}\label{formdec}
      $(X,\mathcal{C})$ is a positive instance of 2-SAT if and only if $D$ is a positive instance of $(2,1)$-BDLFD. Further, a $(2,1)$-decomposition of $D$ can be obtained from a satisfying assignment for $(X,\mathcal{C})$ in polynomial time.
  \end{claim}
  \begin{proofclaim}
      First suppose that $(X,\mathcal{C})$ is a positive instance of 2-SAT, so there exists a satisfying assignment $\phi \colon X \rightarrow \{\true, \false\}$ for $(X,\mathcal{C})$. Consider some $P \in \mathcal{P}_0$ and let $P'$ be the associated path of $P$. If $A(P)$ contains only one arc $a$, we define a decomposition $(F_1^P,F_2^P)$ of $P$ by $A(F_1^P)=A(P)$ and $A(F_2^P)=\emptyset$ if $\phi(x_a)=\true$ and by $A(F_2^P)=A(P)$ and $A(F_1^P)=\emptyset$ if $\phi(x_a)=\false$. Now suppose that the length of $\UG(P')$ is at least 2 and let $a_1$ and $a_2$ be the endarcs of $P'$. As $\phi$ is satisfying for $(X,\mathcal{C})$, there exists a $(2,1)$-decomposition $(F_2^P,F_1^P)$ of $P'$ such that for $i \in [2]$, we have that $a_i$ is contained in $A(F_1)$ if $\phi(x_{a_i})=\true$ and $a_i$ is the only arc in a connected component of $F_2$ if $\phi(x_{a_i})=\false$. 

      Now consider some $P \in \mathcal{P}_1$ and let $a$ be the unique arc in $A(P)$ which is incident to a vertex in $V_3$. If $\phi(x_a)=\true$, we choose a $(2,1)$-decomposition $(F_2^P,F_1^P)$ of $P$ with $a \in A(F_1^P)$. Otherwise, we choose a $(2,1)$-decomposition $(F_2^P,F_1^P)$ of $P$ with $a \in A(F_2^P)$ and $a$ is the only arc in a connected component of $F_2$. 
    By Proposition~\ref{cycchem2}, such a decomposition exists and can be computed in polynomial time. We now define a decomposition $(F_2,F_1)$ of $D$ by $A(F_2)=\bigcup_{P \in \mathcal{P}_0 \cup \mathcal{P}_1}A(F_2^P)$ and $A(F_1)=\bigcup_{P \in \mathcal{P}_0 \cup \mathcal{P}_1}A(F_1^P)$. 

      We now show that $(F_2,F_1)$ is a $(2,1)$-decomposition of $D$. Let $K$ be a connected component of $F_2$ or $F_1$. If $V(K)\subseteq V(P)-V_3$ for some $P \in \mathcal{P}_0 \cup \mathcal{P}_1$, then $K$ is a 2-bounded directed path if $K$ is a connected component of $F_2$ and $1$-bounded directed path if $F$ is a connected component of $F_1$ by construction. We may hence suppose that $V(K)$ contains a vertex $v \in V_3$. Let $a_1,a_2,$ and $a_3$ be the arcs incident to $v$ such that $|\delta_D^+(v)\cap \{a_2,a_3\}| \in \{0,2\}$. As the clauses in $\mathcal{C}_v$ are satisfied by  $\phi$, we obtain that $\phi(x_{a_1})=\false$ and exactly one of $x_{a_2}$ and $x_{a_3}$ is $\true$ under $\phi$, say $\phi(x_{a_2})=\true$ and $\phi(x_{a_3})=\false.$ As $\phi$ satisfies $(X,\mathcal{C})$, we obtain that that none of $a_1$ and $a_3$ are contained in a path of $\mathcal{P}_0$ of length 1. Hence if $K$ is a connected component of $F_2$, we have $A(K)=\{a_1,a_3\}$. Further, as $\phi$ satisfies $(X,\mathcal{C})$, we obtain that $a_1$ and $a_3$ do not form a digon and so $K$ is a 2-bounded directed path. Next, observe that, as the choice of $v$ was arbitrary, the second endvertex of $a_2$ is not incident to any arc in $A(F_1)$ except $a_2$. It follows that $a_2$ is the only arc of $A(K)$ if $K$ is a connected component of $F_1$. Hence $(F_2,F_1)$ is a $(2,1)$-decomposition of $D$. Observe that $(F_2,F_1)$ can be obtained in polynomial time from $\phi$.

      Now suppose that $D$ is a positive instance of $(2,1)$-BDLFD, so there exists a $(2,1)$-decomposition $(F_2,F_1)$ of $D$. We now define a truth assignment $\phi \colon X \rightarrow \{\true,\false\}$ in the following way: for every $a \in A$ which is incident to at least one vertex in $V_3$ in $D$, we set $\phi(x_a)=\true$ if $a \in A(F_1)$ and $\phi(x_a)=\false$ if $a \in A(F_2)$. We show in the following that $\phi$ is a satisfying assignment for $(X,\mathcal{C})$. First consider some $v \in V_3$ and let $a_1,a_2,a_3$ be the arcs incident to $v$ such that $|\delta_D^+(v)\cap \{a_2,a_3\}| \in \{0,2\}$. As $F_2$ is a directed linear forest, we obtain that one of $a_2$ and $a_3$, say $a_2$, is contained in $A(F_1)$. As $F_1$ is a directed matching, we obtain that $\{a_1,a_3\}\subseteq A(F_2)$. By construction, this yields $\phi(x_{a_1})=\phi(x_{a_3})=\false$ and $\phi(x_{a_2})=\true.$ In particular, all clauses in $\mathcal{C}_v$ are satisfied by $\phi$. Now consider some $P \in \mathcal{P}_0$. If $P$ is a path containing a single arc $a$, then observe that $a \in A(F_1)$, as $a$ is incident to two vertices in $V_3$ and $(F_2,F_1)$ is a $(2,1)$-decomposition. Hence $\phi(x_a)=\true$ and so the unique clause in $\mathcal{C}_P$ is satisfied by $\phi$. Now suppose that $P$ contains at least two arcs, let $P'$ be the associated path of $P$ and let $a_1$ and $a_2$ be the endarcs of $P'$. Observe that by construction, we have that $\mathcal{C}_P$ does not contain the clause $\{y_{a_1},y_{a_2}\}$ where $y_{a_i}=x_{a_i}$ if $a_i \in A(F_2)$ and $y_{a_i}=\bar{x}_{a_i}$ if $a_i \in A(F_1)$ for $i \in [2]$. It follows that all the clauses in $\mathcal{C}_P$ are satisfied by $\phi$. 
      Next consider some $q \in Q$. As $F_2$ and $F_1$ are directed linear forests, we obtain that exactly one of the two arcs in the digon is contained in $A(F_1)$. Hence the clauses in $\mathcal{C}_Q$ are satisfied by $\phi.$
      As $\mathcal{C}=\bigcup_{v \in V_3 \cup \mathcal{P}_0 \cup Q}\mathcal{C}_v$, we obtain that $\phi$ is a satisfying assignment for $(X,\mathcal{C})$. Hence $\phi$ is a satisfying assignment for $(X,\mathcal{C})$.
  \end{proofclaim}
  By Claim~\ref{formdec}, it suffices to decide whether $(X, \mathcal{C})$ is a positive instance of 2-SAT. By Proposition~\ref{prop:2SAT}, this can be done in polynomial time. Further, as the proof is algorithmic, a $(2,1)$-decomposition of $D$ can be constructed in polynomial time.
\end{proof}

\subsection{Decomposing into a matching and a (\texorpdfstring{$\geq 3$}{>= 3})-bounded directed linear forest}
\label{sec:1_3_case}

This section is dedicated to proving our hardness results for the case that $k \geq 3$ is finite and $\ell=1$. We first need the following preliminary constructions.
For some integer $k\geq 2$, a {\bf short $k$-in-forcer} is a digraph $D$ together with an arc $a$ whose head is a vertex $z$ with $d_D^+(z)=0$ and $d_D^-(z)=1$ such that there is a $(k,1)$-decomposition of $D$ and for every $(k,1)$-decomposition $(F_k,F_1)$ of $D$, we have that $a \in A(F_1)$. We call $z$ the {\bf tip} of the short $k$-in-forcer.

\begin{proposition}
    \label{prop:short_k_forcer}
    For every $k \geq 2$, there exists a short $k$-in-forcer.
\end{proposition}
\begin{proof}
Let $D$ be the digraph defined by $V(D)=\{v_1,\ldots,v_4,z\}$ and $A(D)=\{v_1v_2,v_2v_3,v_2v_4,v_1z\}$ and let $a=v_1z$.
For an illustration, see Figure~\ref{short1}.

 \begin{figure}[hbt!]
    \begin{center}	
          \begin{tikzpicture}[thick,scale=1, every node/.style={transform shape}]
            \tikzset{vertex/.style = {circle,fill=black,minimum size=5pt, inner sep=0pt}}
            \tikzset{edge/.style = {->,> = latex'}}
            \begin{scope}[xshift=5cm]
                \node[vertex, label=left:$z$] (zb) at (-2,0) {};
                \node[vertex, label=below:$v_1$] (v1b) at (-1,0) {};
                \node[vertex, label=below:$v_2$] (v2b) at (0,0) {};
                \node[vertex, label=right:$v_3$] (v3b) at (45:1) {};
                \node[vertex, label=right:$v_4$] (v4b) at (-45:1) {};
                \draw[edge, green] (v1b) to (v2b);
                \draw[edge, red, dashed] (v1b) to (zb);
                \draw[edge, red, dashed] (v2b) to (v3b);
                \draw[edge, green] (v2b) to (v4b);
            \end{scope}
          \end{tikzpicture}
      \caption{
      A $(k,1)$-decomposition $(F_k,F_1)$ of the short $k$-in-forcer $D$. The dashed red arcs are in $A(F_1)$ and the solid green arcs are in $A(F_k)$.
      }
      \label{short1}
    \end{center}
\end{figure}
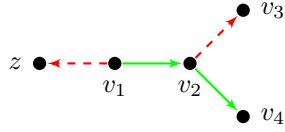

  Let $(F_k,F_1)$ be a $(k,1)$-decomposition of $D$. As $F_k$ is a directed linear forest, we have $v_2v_i \in A(F_1)$ for some $i \in \{3,4\}$. As $F_1$ is a directed matching, we obtain $v_1v_2 \in A(F_k)$. As $F_k$ is a directed linear forest, we obtain $a \in A(F_1)$. Further, it is easy to see that $(F_k,F_1)$ is indeed a $(k,1)$-decomposition. 
\end{proof}

For some integers $k\geq \alpha \geq 1$, a {\bf long $(k,\alpha)$-in-forcer} is a digraph $D$ together with a vertex $z$ satisfying $d_D^+(z)=0$ and $d_D^-(z)=1$ such that there is a $(k,1)$-decomposition $(F_k,F_1)$ of $D$ in which $z$ is not the last vertex of a path of length $\alpha+1$ in $F_k$ and for every $(k,1)$-decomposition $(F_k,F_1)$ of $D$, we have that $z$ is the last vertex of a path of length $\alpha$ in $F_k$. We call $z$ the {\bf tip} of the long $(k,\alpha)$-in-forcer.

\begin{proposition}
    For all integers $k, \alpha$ with $k > \alpha \geq 1$, there exists a long $(k,\alpha)$-in-forcer.
\end{proposition}
\begin{proof}
Let $D$ be obtained from a directed path $v_1\ldots v_{\alpha+1}$ of length $\alpha$ by identifying each vertex of $\{v_1,\ldots,v_{\alpha}\}$ with the tip of a short $k$-in-forcer (the existence of which is guaranteed by Proposition~\ref{prop:short_k_forcer}). Further, we let $z=v_{\alpha+1}$. For an illustration, see Figure~\ref{long}.

 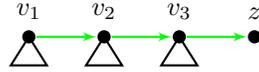
\begin{figure}[hbt!]
    \begin{center}	
          \begin{tikzpicture}[thick,scale=1, every node/.style={transform shape}]
            \tikzset{vertex/.style = {circle,fill=black,minimum size=5pt, inner sep=0pt}}
            \tikzset{edge/.style = {->,> = latex'}}

            \begin{scope}[xshift=6cm]
                \node[vertex, label=above:$v_1$] (v1b) at (-1,0) {};
                \node[vertex, label=above:$v_2$] (v2b) at (0,0) {};
                \node[vertex, label=above:$v_3$] (v3b) at (1,0) {};
                \node[vertex, label=above:$z$] (v4b) at (2,0) {};
                \draw[edge, green] (v1b) to (v2b);
                \draw[edge, green] (v2b) to (v3b);
                \draw[edge, green] (v3b) to (v4b);

                \draw[] (0,0) -- (0.24, -0.38) -- (- 0.24, - 0.38) -- (0,0) {};
                \begin{scope}[xshift=1cm]
                \draw[] (0,0) -- (0.24, -0.38) -- (- 0.24, - 0.38) -- (0,0) {};
                \end{scope}
                \begin{scope}[xshift=-1cm]
                \draw[] (0,0) -- (0.24, -0.38) -- (- 0.24, - 0.38) -- (0,0) {};
                \end{scope}
            \end{scope}
          \end{tikzpicture}
      \caption{
      A $(k,1)$-decomposition $(F_k,F_1)$ of a long $(k,3)$-in-forcer for some $k \geq 3$ with tip $z$. The triangles indicate short $k$-in-forcers. The solid green arcs are in $A(F_k)$ and the short $k$-in-forcers are decomposed as in Figure~\ref{short1}.
      }
      \label{long}
    \end{center}
\end{figure}

   Let $(F_k,F_1)$ be a $(k,1)$-decomposition of $D$. By the definition of short $k$-in-forcers and as $F_1$ is a directed matching, we obtain that $v_iv_{i+1}\in A(F_k)$ for $i\in [\alpha]$. Hence $z$ is the last vertex of a path of length $\alpha$ in $F_k$. Further, it is not the last vertex of a path of length $\alpha+1$ since $v_1$ has no in-neighbour in $F_k$. 
\end{proof}

In the following, we describe some gadgets we need for our reduction. A {\bf $k$-variable gadget} for some $k \geq 3$ is a digraph $D$ together with four arcs $a_1,\ldots,a_4 \in A(D)$ satisfying the following properties:
\begin{enumerate}[$(a)$]
    \item $d_D^+(z_i)=0$ and $d_D^-(z_i)=1$ where $z_i$ is the head of $a_i$ for $i\in [4]$,
    \item for every $(k,1)$-decomposition $(F_k,F_1)$ of $D$, we have $\{a_1,\ldots,a_4\}\cap A(F_1) \subseteq \{a_{1+i},a_{3+i}\}$ for some $i\in \{0,1\}$, and
    \item for every $i\in \{0,1\}$, there is a $(k,1)$-decomposition $(F_k^{i},F_1^{i})$ of $D$ such that $\{a_1,\ldots,a_4\}\cap A(F_1^{i}) = \{a_{1+i},a_{3+i}\}$. 
\end{enumerate}
\begin{lemma}\label{vark1}
    For every $k \geq 3$, there exists a $k$-variable gadget.
\end{lemma}\begin{proof}
We describe a $k$-variable gadget $D$ for some integer $k \geq 3$.
We first let $V(D)$ contain sets $\{v_1,\ldots,v_{12}\}$, $\{z_1,\ldots,z_4\}$, and $\{y_1,\ldots,y_4\}$ and we let $A(D)$ contain the arcs $v_iv_{i+1}$ for $i\in [11]$, the arc $v_{12}v_1$, and the arcs $y_iz_i$ and $y_iv_{3i-2}$ for $i\in[4]$.
Further, for $i \in [4]$,  we identify $y_i$ with the tip of a long $(k,k-2)$-in-forcer. Finally, for $i\in[4]$, we set $a_i=y_iz_i$. This finishes the description of $D$, see Figure~\ref{var3} for an illustration.

\begin{figure}[hbt!]\begin{center}
    \begin{tikzpicture}[thick,scale=1, every node/.style={transform shape}]
        \tikzset{vertex/.style = {circle,fill=black,minimum size=5pt, inner sep=0pt}}
        \tikzset{edge/.style = {->,> = latex'}}
        \foreach \i in {0,...,3}{
            \pgfmathtruncatemacro{\j}{\i+1}
            \node[vertex, label=above:$v_\j$] (v\i) at (\i, 0) {};
        }
        \foreach \i in {1,2,3}{
            \pgfmathtruncatemacro{\j}{\i+4}
            \pgfmathtruncatemacro{\k}{\i+3}
            \node[vertex, label=right:$v_\j$] (v\k) at (3, -\i) {};
        }
        \foreach \i in {1,2,3}{
            \pgfmathtruncatemacro{\j}{\i+7}
            \pgfmathtruncatemacro{\k}{\i+6}
            \node[vertex, label=below:$v_{\j}$] (v\k) at (3-\i, -3) {};
        }
        \foreach \i in {1,2}{
            \pgfmathtruncatemacro{\j}{\i+10}
            \pgfmathtruncatemacro{\k}{\i+9}
            \node[vertex, label=left:$v_{\j}$] (v\k) at (0, -3+\i) {};
        }
        \foreach \i in {1,5,7,11}{
            \pgfmathtruncatemacro{\j}{Mod(\i+1,12)}
            \draw[edge, red, dashed] (v\i) to (v\j) {};
        }
        \foreach \i in {0,2,3,4,6,8,9,10}{
            \pgfmathtruncatemacro{\j}{Mod(\i+1,12)}
            \draw[edge,green] (v\i) to (v\j) {};
        }
        \node[vertex, label=right:$y_1$] (y1) at (-1/1.414,1/1.414) {};
        \node[vertex, label=right:$z_1$] (z1) at (-2/1.414,2/1.414) {};
        \draw[edge,green] (y1) to (v0) {};
        \draw[edge,red,dashed] (y1) to (z1) {};
        \draw (y1) -- (-1/1.414, 1/1.414 - 0.3535) -- (-1/1.414 - 0.3535, 1/1.414 - 0.3535) -- (-1/1.414 - 0.3535, 1/1.414) -- (y1);

        \node[vertex, label=right:$y_2$] (y2) at (3+1/1.414,1/1.414) {};
        \node[vertex, label=right:$z_2$] (z2) at (3+2/1.414,2/1.414) {};
        \draw[edge,red,dashed] (y2) to (v3) {};
        \draw[edge,green] (y2) to (z2) {};
        \draw (y2) -- (3+1/1.414, 1/1.414 + 0.3535) -- (3+1/1.414 - 0.3535, 1/1.414 + 0.3535) -- (3+1/1.414 - 0.3535, 1/1.414) -- (y2);

        \node[vertex, label=left:$y_3$] (y3) at (3+1/1.414,-3-1/1.414) {};
        \node[vertex, label=left:$z_3$] (z3) at (3+2/1.414,-3-2/1.414) {};
        \draw[edge,green] (y3) to (v6) {};
        \draw[edge,red,dashed] (y3) to (z3) {};
        \draw (y3) -- (3+1/1.414, -3-1/1.414 + 0.3535) -- (3+1/1.414+ 0.3535, -3-1/1.414 + 0.3535) -- (3+1/1.414+ 0.3535, -3-1/1.414) -- (y3);

        \node[vertex, label=left:$y_4$] (y4) at (-1/1.414,-3-1/1.414) {};
        \node[vertex, label=left:$z_4$] (z4) at (-2/1.414,-3-2/1.414) {};
        \draw[edge,red,dashed] (y4) to (v9) {};
        \draw[edge,green] (y4) to (z4) {};
        \draw (y4) -- (-1/1.414, -3-1/1.414 - 0.3535) -- (-1/1.414 + 0.3535, -3-1/1.414 - 0.3535) -- (-1/1.414 + 0.3535, -3-1/1.414) -- (y4);
    \end{tikzpicture}
    \caption{An illustration of $D$ and one of its $(k,1)$-decompositions. The squares indicate long $(k,k-2)$-in-forcers. The red dashed arcs are in $A(F_1)$, the solid green arcs are in $A(F_k)$ and the long $(k,k-2)$-in-forcers are decomposed as in Figure~\ref{long}.}\label{var3}
\end{center}
\end{figure}
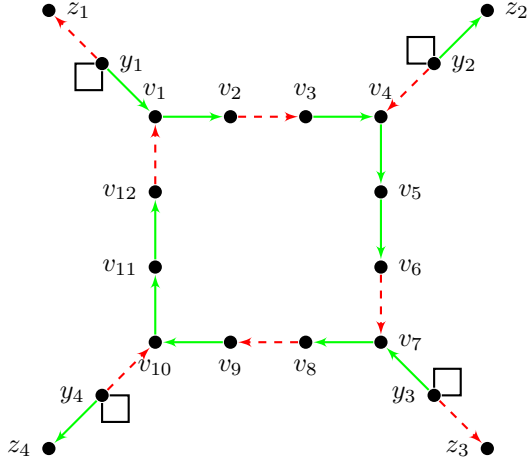

We now show that $D$ is a $k$-variable gadget. It follows by construction that $(a)$ is satisfied.
For $(b)$, by symmetry, it suffices to show that there is no $(k,1)$-decomposition $(F_k,F_1)$ of $D$ such that $\{a_1,a_2\}\subseteq A(F_1)$. Suppose otherwise. As $F_1$ is a directed matching, we obtain that $\{y_1v_1,y_2v_4\}\subseteq A(F_k)$. As $F_k$ is a directed linear forest, we obtain that $\{v_{12}v_1,v_3v_4\}\subseteq A(F_1)$.
As $F_1$ is a directed matching, we obtain that $\{v_1v_2,v_2v_3\}\subseteq A(F_k)$. Hence the directed path obtained from concatenating the directed path of length $k-2$ fully contained in the long $(k,k-2)$-in-forcer incident to $y_1$ whose last vertex is $y_1$ with the directed path $y_1v_1v_2v_3$ is a directed path of length $k+1$ in $F_k$, a contradiction to $F_k$ being a $k$-bounded directed linear forest. This proves $(b)$.

For $(c)$, by symmetry, it suffices to prove that there exists a $(k,1)$-decomposition $(F_k,F_1)$ of $D$ with $\{a_1,a_3\}\subseteq A(F_1)$. Consider the following decomposition $(F_k,F_1)$ of $D$. For $i \in [4]$, we choose a $(k,1)$-decomposition $(F_k^{i},F_1^{i})$ of the $(k,k-2)$-in-forcer incident to $y_i$ such that $y_i$ is not contained in a directed path of length $k-1$ of $F_k^{i}$. We formally define $(F_k,F_1)$, illustrated in Figure~\ref{var3}, in the following way:
\begin{align*}
    A(F_1) &= \bigcup_{i \in [4]}A(F_1^{i})\cup \{v_2v_3,v_6v_7,v_8v_9,v_{12}v_1,y_1z_1,y_2v_4,y_3z_3,y_4v_{10}\}, \\
    \text{and~~} A(F_k) &= A(D) - A(F_1).
\end{align*}
It is straightforward to check that $(F_k,F_1)$ has the desired properties. This proves $(c)$.
\end{proof}

A {\bf $k$-clause gadget} for some integer $k \geq 3$ is a digraph $D$ together with three distinct arcs $b_1,b_2,b_3 \in A(D)$ satisfying the following properties:
\begin{enumerate}[$(a)$]
   
    \item $d_D^+(y_i)=0$ and $d_D^-(y_i)=1$ where $y_i$ is the head of $b_i$ for $i\in[3]$,
    \item for every $(k,1)$-decomposition $(F_k,F_1)$ of $D$, we have $\{b_1,b_2,b_3\}\cap A(F_k) \neq \emptyset$, and
    \item for every nonempty set $S \subseteq \{b_1,b_2,b_3\}$, there is a $(k,1)$-decomposition $(F_k,F_1)$ of $D$ such that $\{b_1,b_2,b_3\}\cap A(F_k) = S$. 
\end{enumerate}

\begin{figure}[hbt!]
    \begin{center}	
          \begin{tikzpicture}[thick,scale=1, every node/.style={transform shape}]
        
            \tikzset{vertex/.style = {circle,fill=black,minimum size=5pt, inner sep=0pt}}
            \tikzset{edge/.style = {->,> = latex'}}

            \pgfmathtruncatemacro{\scaledown}{-3}
            \pgfmathtruncatemacro{\scaleright}{2.2}
            \foreach \i in {0,...,6}{
                \pgfmathtruncatemacro{\k}{Mod(\i,4)}
                \pgfmathtruncatemacro{\l}{(\i - \k)*\scaleright}
                \foreach \j in {1,...,7}{
                    \node[vertex, label=above: $v_\j$] (v\j\i) at (\j + \l, \scaledown*\k) {};
                }
                \node[vertex, label=below:$y_1$] (y1\i) at (1 + \l, \scaledown*\k-1) {};
                \node[] (b1\i) at (0.7 + \l, \scaledown*\k-0.5) {$b_1$};
                \node[vertex, label=below:$y_2$] (y2\i) at (5 + \l, \scaledown*\k-1) {};
                \node[] (b2\i) at (5.3 + \l, \scaledown*\k-0.5) {$b_2$};
                \node[vertex, label=below:$y_3$] (y3\i) at (7 + \l, \scaledown*\k-1) {};
                \node[] (b3\i) at (7.3 + \l, \scaledown*\k-0.5) {$b_3$};
                \draw[] (v3\i) -- (3+\l - 0.25, \scaledown*\k-0.25) -- (3+\l, \scaledown*\k-0.5) -- (3+\l + 0.25, \scaledown*\k-0.25) -- (v3\i) {};

                \ifthenelse{\i=0 \OR \i=3 \OR \i=4 \OR \i=6}{
                    \draw[edge,dashed,red] (v1\i) to (v2\i) {};
                    \draw[edge,green] (v1\i) to (y1\i) {};
                    \draw[edge,green] (v3\i) to (v2\i) {};
                    \draw[edge,dashed,red] (v3\i) to (v4\i) {};
                    \draw[edge,green] (v4\i) to (v5\i) {};
                }{
                    \draw[edge,green] (v1\i) to (v2\i) {};
                    \draw[edge,dashed,red] (v1\i) to (y1\i) {};
                    \draw[edge,dashed,red] (v3\i) to (v2\i) {};
                    \draw[edge,green] (v3\i) to (v4\i) {};
                    \draw[edge,green] (v4\i) to (v5\i) {};
                }
                \ifthenelse{\i=0 \OR \i=2 \OR \i=4}{
                    \draw[edge,dashed,red] (v5\i) to (y2\i) {};
                    \draw[edge,green] (v5\i) to (v6\i) {};
                }{
                    \draw[edge,green] (v5\i) to (y2\i) {};
                    \draw[edge,dashed,red] (v5\i) to (v6\i) {};
                }
                \ifthenelse{\i=2 \OR \i=4}{
                    \draw[edge,dashed,red] (v6\i) to (v7\i) {};
                }{
                    \draw[edge,green] (v6\i) to (v7\i) {};
                }
                \ifthenelse{\i=0 \OR \i=1 \OR \i=3}{
                    \draw[edge,dashed,red] (v7\i) to (y3\i) {};
                }{
                    \draw[edge,green] (v7\i) to (y3\i) {};
                }
            }
            \node[] (S0) at (4, \scaledown*0 - 1.8) {$S=\{b_1\}$};
            \node[] (S1) at (4, \scaledown*1 - 1.8) {$S=\{b_2\}$};
            \node[] (S2) at (4, \scaledown*2 - 1.8) {$S=\{b_3\}$};
            \node[] (S3) at (4, \scaledown*3 - 1.8) {$S=\{b_1,b_2\}$};
            \node[] (S4) at (4+ 4*\scaleright, \scaledown*0 - 1.8) {$S=\{b_1,b_3\}$};
            \node[] (S5) at (4+ 4*\scaleright, \scaledown*1 - 1.8) {$S=\{b_2,b_3\}$};
            \node[] (S6) at (4+ 4*\scaleright, \scaledown*2 - 1.8) {$S=\{b_1,b_2,b_3\}$};
            \end{tikzpicture}
      \caption{The different decompositions of $D$. The green solid arcs belong to $A(F_k)$ and the dashed red arcs belong to $A(F_1)$. Each square marks a long $(k,k-3)$-in-forcer, whose decomposition is the one of Figure~\ref{long}.}
      \label{gadcl42}
    \end{center}
\end{figure}
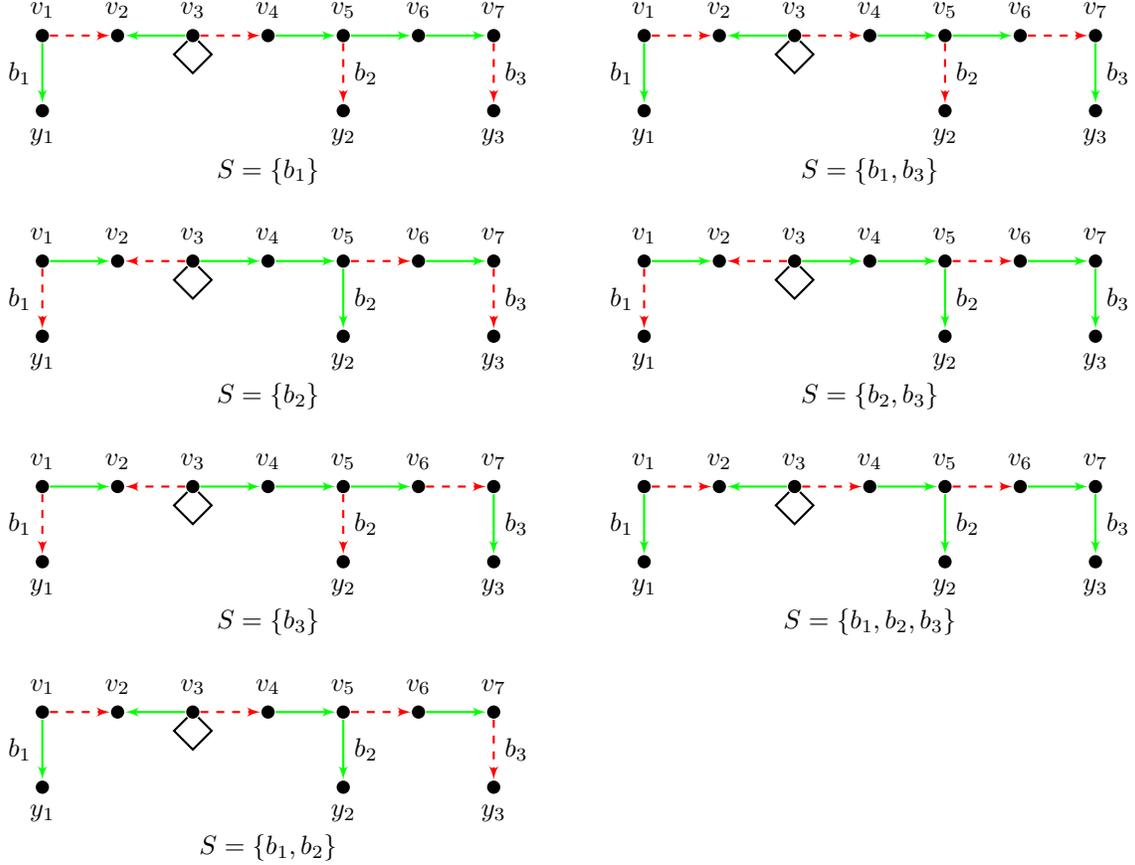

\begin{lemma}\label{clk1}
    For every $k \geq 3$, there exists a $k$-clause gadget.
\end{lemma}
\begin{proof}

We create a $k$-clause gadget $D$ for some fixed $k \geq 3$. We first let $V(D)$ contain a set $\{v_1,\ldots,v_7,y_1,y_2,y_3\}$ and we let $A(D)$ contain the arc set $\{v_1v_2,v_3v_2,v_3v_4,v_4v_5,v_5v_6,v_6v_7,v_1y_1,v_5y_2,v_7y_3\}$. If $k\geq 4$, we further identify $v_3$ with the tip of a long $(k,k-3)$-in-forcer. Finally, let $b_1=v_1y_1,b_2=v_5y_2$, and $b_3=v_7y_3$. This finishes the description of $D$. An illustration can be found in Figure~\ref{gadcl4}.

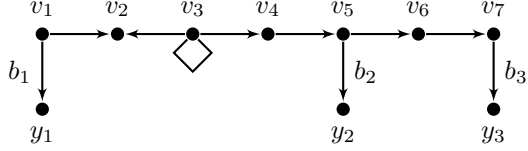
\begin{figure}[hbt!]\begin{center}
  \begin{tikzpicture}[thick,scale=1, every node/.style={transform shape}]
        \tikzset{vertex/.style = {circle,fill=black,minimum size=5pt, inner sep=0pt}}
        \tikzset{edge/.style = {->,> = latex'}}
        \foreach \j in {1,...,7}{
            \node[vertex, label=above: $v_\j$] (v\j) at (\j, 0) {};
        }
        \node[vertex, label=below:$y_1$] (y1) at (1,-1) {};
        \node[] (b1) at (0.7,-0.5) {$b_1$};
        \node[vertex, label=below:$y_2$] (y2) at (5,-1) {};
        \node[] (b2) at (5.3,-0.5) {$b_2$};
        \node[vertex, label=below:$y_3$] (y3) at (7,-1) {};
        \node[] (b3) at (7.3,-0.5) {$b_3$};
        \draw[] (v3) -- (2.75, -0.25) -- (3,-0.5) -- (3.25,-0.25) -- (v3) {};

        \draw[edge] (v1) to (v2) {};
        \draw[edge] (v1) to (y1) {};
        \draw[edge] (v3) to (v2) {};
        \draw[edge] (v3) to (v4) {};
        \draw[edge] (v4) to (v5) {};
        \draw[edge] (v5) to (y2) {};
        \draw[edge] (v5) to (v6) {};
        \draw[edge] (v6) to (v7) {};
        \draw[edge] (v7) to (y3) {};
    
    \end{tikzpicture}
    \caption{An illustration of $D$ where the square marks a long $(k,k-3)$-in-forcer. When $k=3$, the square is deleted.}\label{gadcl4}
\end{center}
\end{figure}

We now show that $(D,b_1,b_2,b_3)$ is a $k$-clause gadget.
By construction, we have that $(a)$ is satisfied.
In order to prove $(b)$, suppose for the sake of a contradiction that there is a $(k,1)$-decomposition $(F_k,F_1)$ of $D$ with $\{b_1,b_2,b_3\}\subseteq A(F_1)$. As $F_1$ is a directed matching, we obtain that $\{v_1v_2,v_4v_5,v_5v_6,v_6v_7\}\subseteq A(F_k)$. We obtain that $v_3v_2 \in A(F_1)$ and so $v_3v_4 \in A(F_k)$. Further, by the definition long $(k,k-3)$-in-forcers, we obtain that $v_3$ is the last vertex of a path of length $k-3$ in $F_k$. Concatenating this path with $v_3v_4v_5v_6v_7$, we obtain that $F_k$ contains a directed path of length $k+1$, a contradiction to $F_k$ being a $k$-bounded linear forest. This proves $(b)$.

We now consider a decomposition of the long $(k,k-3)$-in-forcer in which $v_3$ is not the last vertex of a directed path of length $k-2$ in $F_k$. A small case analysis shows that for every nonempty $S \subseteq \{b_1,b_2,b_3\}$, this can be extended to a $(k,1)$-decomposition $(F_k,F_1)$ of $D$ such that $\{b_1,b_2,b_3\}\cap F_k = S$. All these cases are illustrated in Figure~\ref{gadcl42}. This proves $(c)$.
\end{proof}

We are now ready to prove the following theorem which is the main result of this section.
\begin{theorem}\label{k1}
    $(k,1)$-BDLFD is NP-complete for every integer $k \geq 3$.
\end{theorem}

\begin{proof} 
    We fix some integer $k \geq 3$. Clearly, $(k,1)$-BDLFD is in NP. We prove the hardness by a reduction from $(3,B2)$-SAT, which is NP-complete by Theorem~\ref{thm:3_B2_SAT_NPc}.
    Let $(X,\mathcal{C})$ be an instance of $(3,B2)$-SAT. We now create an instance $D$ of $(k,1)$-BDLFD.
    
    For every $x \in X$, let $C_1,\ldots,C_4$ be an ordering of the clauses containing $x$ or $\bar{x}$ such that $x \in C_1 \cap C_3$ and $\bar{x}\in C_2 \cap C_4$. We add a $k$-variable gadget $(D^x,a_{C_1}^x,\ldots,a_{C_4}^x)$ and for $i=1,\ldots,4$, we let $z_{C_i}^x$ be the head of $a_{C_i}^x$. Observe that this gadget exists by Lemma~\ref{vark1}.

 For every $C \in \mathcal{C}$, let $x_1,x_2,x_3$ be an arbitrary ordering of the variables $x$ for which $x \in C$ or $\bar{x}\in C$ hold. We add a $k$-clause gadget $(D^C,b_{C}^{x_1},b_{C}^{x_2},b_{C}^{x_3})$ and for $i\in[4]$, we let $y_{C_i}^x$ be the head of $b_{C_i}^x$. Observe that this gadget exists by Lemma~\ref{clk1}.
 We finally obtain $D$ by identifying $z_C^x$ and $y_C^x$ for all $x \in X$ and $C \in \mathcal{C}$ for which $x \in C$ or $\bar{x}\in C$ holds.

 We show in the following that $D$ is a positive instance of $(k,1)$-BDLFD if and only if $(X,\mathcal{C})$ is a positive instance of $(3,B2)$-SAT. First suppose that $D$ is a positive instance of $(k,1)$-BDLFD, so there is a $(k,1)$-decomposition $(F_k,F_1)$ of $D$. For every $x \in X$, as $D^x$ is a variable gadget, we obtain that $\{a_{C_1}^x,\ldots,a_{C_4}^x\}\cap A(F_1) \subseteq \{a_{C_{1+i}}^x,a_{C_{3+i}}^x\}$ for some $i\in \{0,1\}$ where $C_1,\ldots,C_4$ is the ordering of the clauses containing $x$ or $\bar{x}$ chosen in the construction of $D$. We now define an assignment $\phi\colon X \rightarrow \{\true,\false\}$ in the following way: We set $\phi(x)=\true$ if $\{a_{C_1}^x,\ldots,a_{C_4}^x\}\cap A(F_1) \subseteq \{a_{C_{1}}^x,a_{C_{3}}^x\}$, and $\phi(x)=\false$, otherwise. In order to see that $\phi$ is a satisfying assignment for $(X,\mathcal{C})$, let $C \in \mathcal{C}$. Suppose that $C$ is not satisfied by $\phi$, hence $\phi(x)=\false$ for all $x \in X$ with $x \in C$ and $\phi(x)=\true$ for all $x \in X$ with $\bar{x} \in C$. It follows by construction that $\{a_{C}^{x_1},a_{C}^{x_2},a_{C}^{x_3}\}\subseteq A(F_k)$ where $\{x_1,x_2,x_3\}$ is the set of variables $x \in X$ such that $x \in C$ or $\bar{x}\in C$ holds. By construction and as $(F_k,F_1)$ is a $(k,1)$-decomposition of $D$, we obtain that $\{b_{C}^{x_1},b_{C}^{x_2},b_{C}^{x_3}\}\subseteq A(F_1)$. By the definition of $k$-clause gadgets, we obtain a contradiction to $(F_k,F_1)$ being a $(k,1)$-decomposition of $D$. Hence $\phi$ is a satisfying assignment for $(X,\mathcal{C})$ and so $(X,\mathcal{C})$ is a positive instance of $(3,B2)$-SAT.
 
 Now suppose that $(X,\mathcal{C})$ is a positive instance of $(3,B2)$-SAT, so there is a satisfying assignment $\phi \colon X \rightarrow \{\true,\false\}$ for $(X,\mathcal{C})$.
For every $x \in X$, let $C_1,\ldots,C_4$ be the ordering of the clauses in $\mathcal{C}$ containing $x$ or $\bar{x}$. If $\phi(x)=\true$, let $(F_k^x,F_1^x)$ be a $(k,1)$-decomposition of $D^x$ with $\{a^x_{C_1},\ldots,a^x_{C_4}\} \cap A(F_1^x)=\{a_{C_{1}}^x,a_{C_{3}}^x\}$. If $\phi(x)=\false$, let $(F_k^x,F_1^x)$ be a $(k,1)$-decomposition of $D^x$ with $\{a^x_{C_1},\ldots,a^x_{C_4}\} \cap A(F_1^x)=\{a_{C_{2}}^x,a_{C_{4}}^x\}$. Observe that these decompositions exist as $D^x$ is a $k$-variable gadget.

For every $C \in \mathcal{C}$, let $S_C$ contain the arc $b_C^x$ for all $x \in X$ with $x \in C$ and $\phi(x)=\true$ and for all $x \in X$ with $\bar{x} \in C$ and $\phi(x)=\false$. Let $(F_k^C,F_1^C)$ be a $(k,1)$-decomposition of $D^C$ with $\{b^{x_1}_{C},b^{x_2}_{C},b^{x_3}_{C}\} \cap A(F_k^x)=S_C$. Observe that such a decomposition exists as $D^C$ is a clause gadget and $\phi$ satisfies $C$.

Now let $(F_k,F_1)$ be defined by $A(F_k)=\bigcup_{x \in X}A(F_k^x)\cup \bigcup_{C \in \mathcal{C}}A(F_k^C)$ and $A(F_1)=\bigcup_{x \in X}A(F_1^x)\cup \bigcup_{C \in \mathcal{C}}A(F_1^C)$. Observe that by construction, every vertex that is incident to arcs from two different gadgets, is incident to exactly two arcs, one in $A(F_k)$ and one in $A(F_1)$. As $(F_k^x,F_1^x)$ is a $(k,1)$-decomposition of $D^x$ for all $x \in X$ and $(F_k^C,F_1^C)$ is a $(k,1)$-decomposition of $D^C$ for all $C \in \mathcal{C}$, we obtain that $(F_k,F_1)$ is a $(k,1)$-decomposition of $D$. Hence $D$ is a positive instance of $(k,1)$-BDLFD.
\end{proof}

\subsection{Decomposing into two linear forests with length bounds at least 2}
\label{sec:2_2_case}
In this section, we deal with the case that both $k$ and $\ell$ are finite integers and $\min\{k,\ell\}\geq 2$. The reduction is from \textsc{ME-1-SAT} and we use a variable and a clause gadget. Once these gadgets are constructed, our reduction will work for any choice of $k$ and $\ell$ in the considered domain. However, during the construction of the gadgets, the case that $\min\{k,\ell\}=2$ often needs to be treated separately. Again, we first give some preliminary constructions in Section~\ref{prelcon}. After, we describe the variable and clause gadgets in Sections~\ref{variable} and~\ref{clause}, respectively. Finally, we give the reduction in Section~\ref{hardness}.

\subsubsection{Preliminary constructions}\label{prelcon}
We here give some preliminary constructions we need for the gadgets which are described in Sections~\ref{variable} and~\ref{clause}. The first one will play a crucial role in the case $\min\{k,\ell\}\geq 3$.
Let $k,\ell \in \mathbb{Z}_{\geq 3}$ be two integers. A {\bf $(k,\ell,-2)$-in-forcer} is a digraph $D$ with a special arc $a=xz$ such that:
\begin{enumerate}[$(a)$]
    \item $a$ is the only arc incident to $z$,
    \item in every $(k,\ell)$-decomposition $(F_k,F_\ell)$ of $D$, $a$ is either the last arc of a directed path of length $k-2$ in $F_k$ or the last arc of a directed path of length $\ell-2$ in $F_\ell$,
    
    \item there exists a $(k,\ell)$-decomposition $(F_k,F_\ell)$ of $D$ such that $a\in A(F_k)$ and $a$ is not the last arc of a directed path of length $k-1$ in $F_k$, and
    
    \item there exists a $(k,\ell)$-decomposition $(F_k,F_\ell)$ of $D$ such that $a\in A(F_\ell)$ and $a$ is not the last arc of a directed path of length $\ell-1$ in $F_\ell$.
\end{enumerate}

In the following, we show that $(k,\ell,-2)$-in-forcers exist for all integers $k,\ell \in \mathbb{Z}_{\geq 3}$. We first consider the case $k=\ell.$

\begin{lemma}
    \label{lemma:kkm2_forcer}
    For every integer $k \geq 3$, there exists a $(k,k,-2)$-in-forcer.
\end{lemma}
\begin{proof}
We describe a $(k,k,-2)$-in-forcer $(D,a)$ for some fixed $k \geq 3$. We first let $D$ contain the unique orientation of a binary tree of depth $k-3$ that contains a path from $v$ to $x$ for all vertices $v$ of this digraph where $x$ is the tip of the binary tree. We now obtain $D$ by adding another vertex $z$ and the arc $a=xz$.
    See Figure~\ref{fig:k_l_m2_forcer} for an illustration.  Note that $(a)$ clearly holds.

          \begin{figure}[hbt!]
        \begin{center}	
              \begin{tikzpicture}[thick,scale=1, every node/.style={transform shape}]
                \tikzset{vertex/.style = {circle,fill=black,minimum size=5pt, inner sep=0pt}}
                \tikzset{littlevertex/.style = {circle,fill=black,minimum size=4pt, inner sep=0pt}}
                \tikzset{edge/.style = {->,> = latex'}}

                \foreach \i in {1,...,4}{
                    \pgfmathtruncatemacro{\j}{2*\i - 1}
                    \pgfmathtruncatemacro{\k}{2*\i}
                    \node[vertex] (a\j) at (\j,0) {};
                    \node[vertex] (a\k) at (\k,0) {};
                    \node[vertex] (b\i) at (\k-0.5,1) {};
                    \draw[edge, dashed, red] (a\j) to (b\i) {};
                    \draw[edge, green] (a\k) to (b\i) {};
                }
                \node[vertex] (c1) at (2.5,2) {};
                \node[vertex] (c2) at (6.5,2) {};
                \draw[edge, dashed, red] (b1) to (c1) {};
                \draw[edge, green] (b2) to (c1) {};
                \draw[edge, dashed, red] (b3) to (c2) {};
                \draw[edge, green] (b4) to (c2) {};
                \node[vertex, label=below:$x$] (x) at (4.5,3) {};
                \node[vertex, label=above:$z$] (z) at (4.5,4) {};
                \draw[edge, dashed, red] (x) to (z) {};
                \draw[edge, dashed, red] (c1) to (x) {};
                \draw[edge, green] (c2) to (x) {};
              \end{tikzpicture}
          \caption{An illustration of a $(6,6)$-decomposition of a $(6,6,-2)$-in-forcer.}
          \label{fig:k_l_m2_forcer}
        \end{center}
    \end{figure}
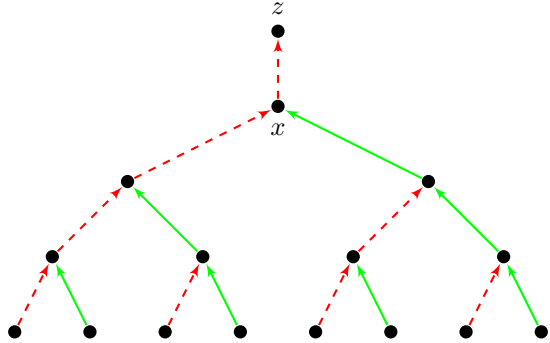
    
    Let $(F_k,F_k')$ be a $(k,k)$-decomposition of $D$. By symmetry, we may suppose that $xz \in A(F_k)$. We will show that $xz$ is the last arc of a directed path of length $k-2$ in $F_k$. Let $P$ be a  directed path of $F_k$ whose length is maximum among all the directed paths in $F_k$ whose last arc is $xz$. If the length of $P$ is smaller than $k-2$, since $D$ is the described orientation of a binary tree of depth $k-3$, the initial vertex $y$ of $P$ has two in-neighbours $y_1,y_2$ in $D$. Hence, as $F_k'$ is a directed linear forest, either $y_1y$ or $y_2y$ belongs to $F_k$, contradicting the maximality of $P$. This proves $(b)$.

    We now describe a $(k,k)$-decomposition $(F_k,F_k')$ of $D$ whose existence proves $(c)$.  For every $v \in V(D)$ with $d_D^-(v)=2$, arbitrarily assign one of its entering arcs to $A(F_k)$ and the other one to $A(F_k')$. Finally assign $xz$ to  $A(F_k)$. This results into a $(k,k)$-decomposition of $D$ with the desired properties, illustrated in Figure~\ref{fig:k_l_m2_forcer}. By symmetry $(d)$ also holds.
   \end{proof}

   For the case $k \neq \ell$, we first need some preliminary constructions. For integers $k > \ell \geq 3$, a {\bf long $(k,\ell)$-out-forcer} is a digraph $D$ together with a special arc $a=xy$ such that $d_D^-(x)=0$, $d_D^+(x)=1$, $D$ admits a $(k,\ell)$-decomposition, and for every $(k,\ell)$-decomposition $(F_k,F_\ell)$ of $D$, $a$ belongs to $A(F_k)$. 
   We say that $x$ is the {\bf origin} of the long $(k,\ell)$-out-forcer.

    \begin{proposition}
    \label{lemma:kllong_forcer}
    For all integers $k,\ell$ with $k>\ell \geq 3$, there exists a long $(k,\ell)$-out-forcer.
\end{proposition}
\begin{proof}
    We describe a long $(k,\ell)$-out-forcer $D$ for some fixed integers $k,\ell$ with $k>\ell \geq 3$. We first let $D$ contain the unique orientation of a binary tree of depth $\ell$ with a root $y$ in which there exists a path from $y$ to $v$ for every vertex $v$. We then create $D$ by adding $x$ and the arc $a=xy$, see Figure~\ref{klout} for an illustration. Clearly, $a$ is the only arc incident to $x$.

    \begin{figure}[hbt!]\begin{center}
        \begin{tikzpicture}[thick,scale=1, every node/.style={transform shape}]
                \tikzset{vertex/.style = {circle,fill=black,minimum size=5pt, inner sep=0pt}}
                \tikzset{littlevertex/.style = {circle,fill=black,minimum size=4pt, inner sep=0pt}}
                \tikzset{edge/.style = {->,> = latex'}}

                \foreach \i in {1,...,4}{
                    \pgfmathtruncatemacro{\j}{2*\i - 1}
                    \pgfmathtruncatemacro{\k}{2*\i}
                    \node[vertex] (a\j) at (\j,0) {};
                    \node[vertex] (a\k) at (\k,0) {};
                    \node[vertex] (b\i) at (\k-0.5,1) {};
                    \draw[edge, green] (b\i) to (a\j){};
                    \draw[edge, dashed, red] (b\i) to (a\k) {};
                }
                \node[vertex] (c1) at (2.5,2) {};
                \node[vertex] (c2) at (6.5,2) {};
                \node[vertex, label=below:$y$] (x) at (4.5,3) {};
                \node[vertex, label=above:$x$] (z) at (4.5,4) {};
                
                \draw[edge, green] (c1) to (b1) {};
                \draw[edge, dashed, red]  (c1) to (b2) {};
                \draw[edge, green] (c2) to (b3){};
                \draw[edge,  dashed, red] (c2) to (b4){};
                \draw[edge, green] (z) to (x) {};
                \draw[edge, green] (x) to (c1) {};
                \draw[edge,  dashed, red] (x) to (c2){};
              \end{tikzpicture}
        \caption{An illustration of a $(k,3)$-decomposition $(F_k,F_3)$ of a long $(k,3)$-out-forcer for some $k\geq 4$. The dashed red arcs are in $A(F_3)$ and the solid green arcs are in $A(F_k)$.}
  \label{klout}
  \end{center}
  \end{figure}
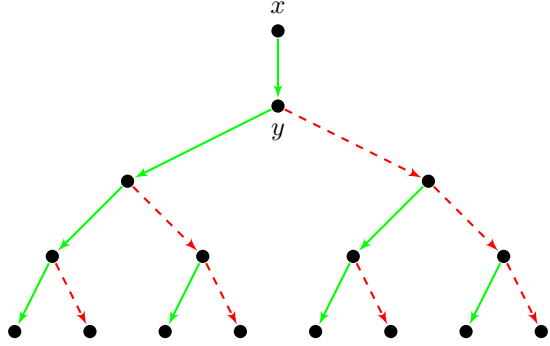 

  Let $(F_k,F_\ell)$ be a $(k,\ell)$-decomposition of $D$ and suppose for the sake of a contradiction that $a \in A(F_\ell)$. Let $P$ be the unique longest path in $F_\ell$ that contains $a$ and let $v$ be the last vertex of $P$. As $F_k$ is a linear directed forest, we obtain that $d_D^+(v)\leq 1$. By construction, we obtain that the length of $P$ is at least $\ell+1$. This contradicts $F_\ell$ being an $\ell$-bounded linear forest.
    
    We now justify the existence of a $(k,\ell)$-decomposition. Let $(F_k,F_\ell)$ be a decomposition of $D$ with $a \in A(F_k)$ and such that for every $v \in V(D)$ with $d_D^+(v)=2$, we have that one of the arcs in $\delta_D^+(v)$ is in $A(F_k)$ and the other one is in $A(F_\ell)$. Then $(F_k,F_\ell)$ is a $(k,\ell)$-decomposition of $D$ with the desired properties. 
\end{proof}

   For some integers $k>\ell \geq 3$, we now define a {\bf short $(k,\ell)$-out-forcer} as a digraph $D$ together with a special arc $a=xy$ such that $d_D^-(x)=0$, $d_D^+(x)=1$, $D$ admits a $(k,\ell)$-decomposition, and for every $(k,\ell)$-decomposition $(F_k,F_\ell)$ of $D$, $a$ belongs to $A(F_\ell)$. 
   Again, we say that $x$ is the {\bf origin} of the short $(k,\ell)$-out-forcer.
   
\begin{proposition}
    \label{lemma:klshort_forcer}
    For all integers $k,\ell$ with $k>\ell \geq 3$, there exists a short $(k,\ell)$-out-forcer.
\end{proposition}
\begin{proof}
    We describe a short $(k,\ell)$-out-forcer $D$ for some fixed integers $k$ and $\ell$ with $k >\ell \geq 3$. We first let $D$ contain the digraph obtained from a long $(k,\ell)$-out-forcer with origin $y$ by reversing all arcs. We then obtain $D$ by adding $x$ and the arc $xy$. Clearly, we have $d_D^-(x)=0$ and $d_D^+(x)=1$. Next observe that there is a direct correspondence between the $(k,\ell)$-decompositions of $D-x$ and the $(k,\ell)$-decompositions of the corresponding long $(k,\ell)$-out-forcer. 

    Let $(F_k,F_\ell)$ be a $(k,\ell)$-decomposition of $D$. Clearly, and $(F_k-x,F_\ell-x)$ is a $(k,\ell)$-decomposition of $D-x$ and so the unique arc incident to $y$ in $D-x$ is contained in $A(F_k)$. As $F_k$ is a linear directed forest, we obtain that $xz \in A(F_\ell)$.
    
    We now justify the existence of a $(k,\ell)$-decomposition of $D$. There exists a $(k,\ell)$-decomposition $(F'_k,F'_\ell)$ of $D-x$ in which the unique arc incident to $y$ is contained in $A(F'_k)$. The decomposition $(F_k,F_\ell)$ of $D$ defined by $A(F_k)=A(F'_k)$ and $A(F_\ell)=A(F'_\ell)\cup xy$ is a $(k,\ell)$-decomposition of $D$ with the desired properties.
\end{proof}

We are now ready to show that $(k,\ell,-2)$-in-forcers exist for all integers $k,\ell \geq 3$. Recall that the case $k=\ell$ has been proved in Lemma~\ref{lemma:kkm2_forcer}.

 \begin{lemma}
    \label{lemma:kldistout_forcer}
    For all integers $k,\ell$ with $k>\ell \geq 3$, there exists a $(k,\ell,-2)$-in-forcer.
\end{lemma}
\begin{proof}
We describe a $(k,\ell,-2)$-in-forcer  $(D,a)$ for some fixed integers $k>\ell \geq 3$. We first let $V(D)$ contain a set  $\{u_1,\ldots,u_{k-3}\}$, a set $\{v_1,\ldots,v_{\ell-3}\}$ if $\ell\geq 4$, and a set $\{x,z\}$. Further, we let $A(D)$ contain the arcs $u_iu_{i+1}$ for $i \in [k-4]$ if $k\geq 5$, the arcs $v_iv_{i+1}$ for $i \in [\ell-4]$ if $\ell \geq 5$ , the arc $u_{k-3}x$, the arc $v_{\ell -3}x$ if $\ell \geq 4$, and the arc $a=xz$. Finally, for $i \in [k-3]$, we identify $u_i$ with the origin of a short $(k,\ell)$-out-forcer and, if $\ell \geq 4$, for $i \in [\ell-3]$, we identify $v_i$ with the origin of a long $(k,\ell)$-out-forcer. This finishes the description of $D$, see Figure~\ref{kl-2} for an illustration. 

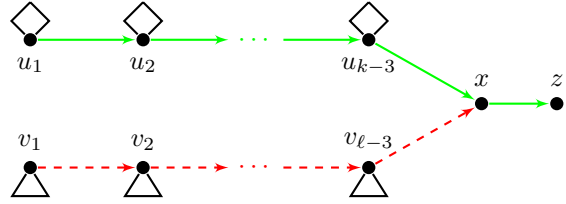
\begin{figure}[hbt!]\begin{center}
   \begin{tikzpicture}[thick,scale=1, every node/.style={transform shape}]
    \tikzset{vertex/.style = {circle,fill=black,minimum size=5pt, inner sep=0pt}}
    \tikzset{littlevertex/.style = {circle,fill=black,minimum size=4pt, inner sep=0pt}}
    \tikzset{edge/.style = {->,> = latex'}}
    
    \foreach \i in {1,2}{
        \node[vertex, label=below:$u_\i$] (u\i) at (\i*1.5,0) {};
        \draw[] (\i*1.5,0.5) -- (\i*1.5-0.25,0.25) -- (u\i) -- (\i*1.5+0.25,0.25) -- (\i*1.5,0.5){};
        
        \node[vertex, label=above:$v_\i$] (v\i) at (\i*1.5,-1.7) {};
        \draw[] (v\i) -- (\i*1.5 + 0.24, -1.7- 0.38) -- (\i*1.5 - 0.24, -1.7 - 0.38) -- (v\i);
    }
    \foreach \i in {4}{
        \node[vertex, label=below:$u_{k-3}$] (u\i) at (\i*1.5,0) {};
        \draw[] (\i*1.5,0.5) -- (\i*1.5-0.25,0.25) -- (u\i) -- (\i*1.5+0.25,0.25) -- (\i*1.5,0.5){};
        
        \node[vertex, label=above:$v_{\ell-3}$] (v\i) at (\i*1.5,-1.7) {};
        \draw[] (v\i) -- (\i*1.5 + 0.24, -1.7 - 0.38) -- (\i*1.5 - 0.24, -1.7 - 0.38) -- (v\i);
    }
    \node[green] (u3) at (4.5,0) {$\cdots$};
    \node[red] (v3) at (4.5,-1.7) {$\cdots$};
    
    \foreach \i in {1,2,3}{
        \pgfmathtruncatemacro{\j}{\i +1}
        \draw[edge,green] (u\i) to (u\j) {};
        \draw[edge,dashed,red] (v\i) to (v\j) {};
    }
    \node[vertex, label=above:$x$] (x) at (7.5,-0.85) {};
    \node[vertex, label=above:$z$] (z) at (8.5,-0.85) {};
    \draw[edge, dashed, red] (v4) to (x);
    \draw[edge, green] (u4) to (x);
    \draw[edge, green] (x) to (z);
  \end{tikzpicture}
  \caption{An illustration of a $(k,\ell)$-decomposition $(F_k,F_\ell)$ of a $(k,\ell,-2)$-in-forcer with $k > \ell\geq 4$ and $xz \in A(F_k)$. The squares indicate short $(k,\ell)$-out-forcers and the triangles indicate long $(k,\ell)$-out-forcers. The dashed red arcs are in $A(F_\ell)$ and the solid green arcs are in $A(F_k)$.}\label{kl-2}
  \end{center}
  \end{figure} 

We now prove that $(a)-(d)$ holds.
Clearly, $(a)$ holds. Now let $(F_k,F_\ell)$ be a $(k,\ell)$-decomposition of $D$. For $i \in [k-3]$, by the definition of short $(k,\ell)$-out-forcers, we have that there exists an arc in $\delta_D^+(u_i)\cap A(F_\ell)$ that is contained in the short $(k,\ell)$-out-forcer attached to $u_i$. It follows that $u_{k-3}z \in A(F_k)$ and, if $k \geq 5$, then $u_iu_{i+1}\in A(F_k)$ for $i \in [k-4]$. It follows that $z$ is the last vertex of a directed path of length $k-3$ in $F_k$. A similar argument shows that $z$ is the last vertex of a directed path of length $\ell-3$ in $F_\ell$. Hence if $a \in A(F_k)$, then $a$ is the last arc of a directed path of length $k-2$ in $F_k$ and if $a \in A(F_\ell)$, then $a$ is the last arc of a path of length $\ell-2$ in $F_\ell$, so $(b)$ holds. Further, by choosing appropriate decompositions of the short and long $(k,\ell)$-out-forcers, if $a \in A(F_k)$, then we have a decomposition as required in $(c)$, and if $a \in A(F_\ell)$, then we have a decomposition as required in $(d)$. 
\end{proof}

  We further need one preliminary construction which will be useful in the case that $\min\{k,\ell\}=2$.
  Let $k,\alpha$ be integers with $k\geq 3$ and $1 \leq \alpha \leq k$. A {\bf $(k,2,\alpha)$-in-forcer} is a digraph $D$ with a special arc $a$ such that:
\begin{enumerate}[$(a)$]
    \item the head of $a$ is not incident to any other arc in $D$, 
    \item in every $(k,2)$-decomposition $(F_k,F_2)$ of $D$, $a$ is the last arc of a path of length at least $\alpha$ in $F_k$, and
    \item $D$ has a $(k,2)$-decomposition $(F_k,F_2)$ in which $a$ is not contained in a path of length $\alpha+1$ in $F_k$.
\end{enumerate}

\begin{lemma}
    \label{lemma:2_k_l_forcer}
    For all integers $k,\alpha$  with $k\geq 3$ and $1 \leq \alpha \leq k$, there exists a $(k,2,\alpha)$-in-forcer.
\end{lemma}
\begin{proof}
   We describe a $(k,2,\alpha)$-in-forcer $(D,a)$ for some fixed integers $k$ and $\alpha$ with $k \geq 3$ and $1 \leq \alpha \leq k$.  
    We let $V(D)$ contain a set of $\alpha+1$ vertices $v_1,\dots,v_{\alpha+1}$. Further, for every $i \in [\alpha]$, we let $V(D)$ contain 6 vertices $u_i,w_i,x_i^1,x_i^2,y_i^1$, and $y_i^2$. For $i\in [\alpha]$, we then let $A(D)$ contain the arcs $v_{i+1}v_i, u_ix_i^1,u_ix_i^2,w_iu_i,v_{i+1}u_i,y_i^1w_i$, and $y_i^2w_i$. We finally set $a=v_2 v_{1}$, which finishes the description of $(D,a)$. For an illustration, see Figure~\ref{vark2}.

    \begin{figure}[hbt!]\begin{center}
        \begin{tikzpicture}[thick,scale=1, every node/.style={transform shape}]
        \tikzset{vertex/.style = {circle,fill=black,minimum size=5pt, inner sep=0pt}}
        \tikzset{littlevertex/.style = {circle,fill=black,minimum size=4pt, inner sep=0pt}}
        \tikzset{edge/.style = {->,> = latex'}}
        
        \foreach \i in {1,2,3,4}{
            \node[vertex, label=above:$v_\i$] (v\i) at (\i*3.5,-0.5) {};
        }
        \foreach \i in {1,2,3}{
            \pgfmathtruncatemacro{\j}{\i +1}
            \draw[edge,green] (v\j) to (v\i) {};

            \node[vertex, label=below:$u_\i$] (u\i) at (\j*3.5, -1.5) {};
            \node[vertex, label=left:$x_\i^1$] (x\i1) at (\j*3.5-0.5, -2) {};
            \node[vertex, label=left:$x_\i^2$] (x\i2) at (\j*3.5-0.5, -1) {};
            
            \node[vertex, label=below:$w_\i$] (w\i) at (\j*3.5+1, -1.5) {};
            \node[vertex, label=right:$y_\i^1$] (y\i1) at (\j*3.5+1.5, -2) {};
            \node[vertex, label=right:$y_\i^2$] (y\i2) at (\j*3.5+1.5, -1) {};

            \draw[edge,dashed, red] (v\j) to (u\i) {};
            \draw[edge,dashed, red] (u\i) to (x\i2) {};
            \draw[edge,green] (u\i) to (x\i1) {};
            \draw[edge,green] (w\i) to (u\i) {};
            \draw[edge,dashed, red] (y\i2) to (w\i) {};
            \draw[edge,green] (y\i1) to (w\i) {};
        }
        \end{tikzpicture}
        \caption{An illustration of a $(k,2)$-decomposition of a $(k,2,3)$-in-forcer. The dashed red arcs are in $A(F_2)$ and the solid green arcs are in $A(F_k)$.
        }\label{vark2}
    \end{center}
    \end{figure}
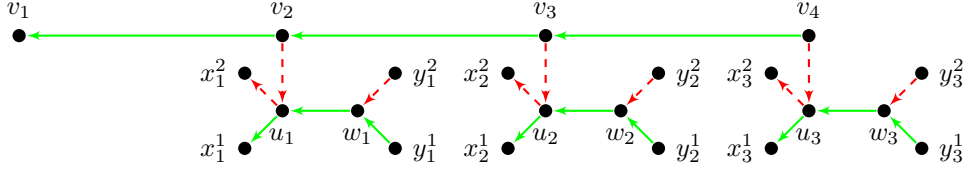

    Note that $(a)$ clearly holds. Now, let $(F_k,F_2)$ be a $(k,2)$-decomposition of $D$. For every $i \in [\alpha]$, as $F_k$ is a directed linear forest, we obtain that one of $u_ix_i^1$ and $u_ix_i^2$ is contained in $A(F_2)$, and one of $y_i^1w_i$, and $y_i^2w_i$ is contained in $A(F_2)$. As $F_2$ is a 2-bounded directed linear forest, we obtain that $w_iu_i\in A(F_k)$. As $F_k$ is a directed linear forest, we obtain that $v_{i+1}u_i\in A(F_2)$. As $F_2$ is a 2-bounded directed linear forest, we obtain that $v_{i+1}v_i\in A(F_k)$. Hence $v_{\ell+1} \ldots v_{1}$ is a path of $F_k$ of length $\alpha$ whose last arc is $a$. This proves $(b)$.

    On the other hand, let $(F_k,F_2)$ be the decomposition of $D$
  with $A(F_k)=\bigcup_{i=1}^\alpha\{v_{i+1}v_{i},u_ix_i^1,w_iu_i,y_i^1w_i\}$ and $A(F_2)=\bigcup_{i=1}^\alpha\{u_ix_i^2,v_{i+1}u_i,y_i^2w_i\}$. Then $(F_k,F_2)$ is a $(k,2)$-decomposition of $D$ with the desired properties. For an illustration, see Figure~\ref{vark2}. This proves $(c)$.
    \end{proof}

\subsubsection{Variable gadgets}\label{variable}
In this section, we describe the variable gadgets. 
For integers $k$,$\ell$, and $t$ with $k\geq \ell\geq 2$ and $t\geq 1$, a {\bf $(k,\ell,t)$-variable gadget} is a digraph $D$ together with a collection of arcs $\{a_1,\ldots,a_t\}$ with the following properties:
  \begin{itemize}
      \item[$(a)$] for each $a \in \{a_1,\ldots,a_t\}$, we have that $a$ is the only arc incident to its head,
      \item[$(b)$] for every $(k,\ell)$-decomposition $(F_k,F_\ell)$ of $D$, we have $A(F_k)\cap\{a_1,\ldots,a_t\}\in \{\{a_1,\ldots,a_t\},\emptyset\}$ and each arc in $A(F_k)\cap\{a_1,\ldots,a_t\}$ is the last arc of a path either of length $k$ in $F_k$ or of length $\ell$ in $F_\ell$,
      \item[$(c)$] there is a $(k,\ell)$-decomposition $(F_k,F_\ell)$ of $D$ such that $\{a_1,\ldots,a_t\} \subseteq A(F_k)$, and
       \item[$(d)$] there is a $(k,\ell)$-decomposition $(F'_k,F'_\ell)$ of $D$ such that $\{a_1,\ldots,a_t\} \subseteq A(F_\ell')$
  \end{itemize}
  We now show that variable gadgets exist for all integers in the considered domain. We start with the case that $\min\{k,\ell\}\geq 3$.
\begin{lemma}
    \label{lemma:k_l_t_gadget}
    For all integers $k,\ell$ and $t\geq 1$ with $k\geq \ell \geq 3$, there exists a $(k,\ell,t)$-variable gadget whose size is polynomial in $t$.
\end{lemma}
\begin{proof}
    We construct a $(k,\ell,t)$-variable gadget $(D,\{a_1,\dots,a_t\})$ for some fixed  integers $k,\ell$, and $t$ with $k\geq \ell \geq 3$ and $t \geq 1$. 
    
    We first let $D$ contain a directed path on $2t+1$ vertices $u_0,\dots,u_{2t}$. Then for every $i\in[2t-1]$, we add a vertex $v_i$ and the arc $u_iv_i$. We further add two more vertices $w_1$ and $w_2$ and the arcs $u_{2t}w_1$ and $u_{2t}w_2$. We now obtain $D$ by identifying $u_i$ with the tip of a $(k,\ell,-2)$-in-forcer for all $i \in [2t-1]$ and identifying $u_0$ with the tips of two $(k,\ell,-2)$-in-forcers, where all these $(k,\ell,-2)$-in-forcers are vertex-disjoint before the identification. Observe that these $(k,\ell,-2)$-in-forcers exist by Lemmas~\ref{lemma:kkm2_forcer} and~\ref{lemma:kldistout_forcer}. 
    For every $i\in[t]$ we further set $a_i = u_{2i-1}v_{2i-1}$. This completes the description of $(D,\{a_1,\dots,a_t\})$, see Figure~\ref{fig:vertex_gadget_k_l_decomposition} for an illustration.

    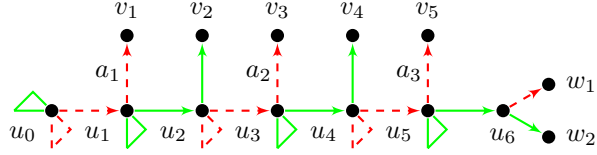
\begin{figure}[hbt!]
        \begin{center}	
              \begin{tikzpicture}[thick,scale=1, every node/.style={transform shape}]
                
                \tikzset{vertex/.style = {circle,fill=black,minimum size=5pt, inner sep=0pt}}
                \tikzset{littlevertex/.style = {circle,fill=black,minimum size=4pt, inner sep=0pt}}
                \tikzset{edge/.style = {->,> = latex'}}
                
                \foreach \i in {1,3,5}{
                    \node[vertex, label={[label distance=0.1]-135:$u_\i$}] (u\i) at (\i,0) {};
                    \draw[green] (\i,-0.5) -- (u\i) -- (\i+0.25,-0.25) -- (\i,-0.5){};
                }
                \foreach \i in {0,2,4}{
                    \node[vertex, label={[label distance=0.1]-135:$u_\i$}] (u\i) at (\i,0) {};
                    \draw[dashed,red] (\i,-0.5) -- (u\i) -- (\i+0.25,-0.25) -- (\i,-0.5){};
                }
                \node[vertex, label=below:$u_6$] (u6) at (6,0) {};
                \draw[green] (-0.5,0) -- (u0) -- (-0.25,0.25) -- (-0.5,0){};
                
                \foreach \i in {0,2,4}{
                    \pgfmathtruncatemacro{\j}{\i +1}
                    \draw[edge,dashed,red] (u\i) to (u\j);
                }
                \foreach \i in {1,3,5}{
                    \pgfmathtruncatemacro{\j}{\i +1}
                    \draw[edge,green] (u\i) to (u\j);
                }
                \foreach \j in {1,2,3}{
                    \pgfmathtruncatemacro{\i}{2*\j-1}
                    \node[vertex, label=above:$v_\i$] (v\i) at (\i,1) {};
                    \draw[edge,dashed,red] (u\i) to (v\i);
                    
                    \node[] (a\i) at (\i-0.25,0.5) {$a_{\j}$};
                }
                \foreach \i in {2,4}{
                    \node[vertex, label=above:$v_\i$] (v\i) at (\i,1) {};
                    \draw[edge,green] (u\i) to (v\i);
                }
                \node[vertex, xshift=6cm, label=right:$w_1$] (w1) at (30:0.7){};
                \node[vertex, xshift=6cm, label=right:$w_2$] (w2) at (-30:0.7){};
                \draw[edge,dashed,red] (u6) to (w1);
                \draw[edge,green] (u6) to (w2);
              \end{tikzpicture}
          \caption{An illustration of a $(k,\ell)$-decomposition of the $(k,\ell,3)$-variable gadget with $k\geq \ell \geq 3$. Triangles indicate $(k,\ell,-2)$-in-forcers. The dashed red arcs are in $A(F_\ell)$ and the solid green arcs are in $A(F_k)$. The colour of a forcer indicates the part of the decomposition the arc incident to its tip is contained in.}
          \label{fig:vertex_gadget_k_l_decomposition}
        \end{center}
    \end{figure}
    
    First observe that $(a)$ clearly holds.
    To show that $(b)$ holds, let us fix a $(k,\ell)$-decomposition $(F_k,F_\ell)$ of $D$.
    We show by induction on $i\in [2t-1]$ that $u_iv_i$ is either the last arc of a directed path of length $\ell$ in $F_\ell$ or the last arc of a directed path of length $k$ in $F_k$. Moreover, we show that this directed path also contains $u_{i-1}u_i$. Observe that this implies $(b)$.
    
    When $i=1$, by definition of a $(k,\ell,-2)$-in-forcer, we know that one entering arc of $u_0$ is the last arc of a directed path of length $k-2$ in $F_k$, and that the other one is the last arc of a directed path of length $\ell-2$ in $F_\ell$. Hence $u_0u_1$ is either the last arc of a directed path of length $k-1$ in $F_k$ or the last arc of a directed path of length $\ell-1$ in $F_\ell$. In both cases, since $u_2$ has out-degree $2$, $u_1v_1$ and $u_0u_1$ must belong to the same part.

   Now assume that $i\in \{2,\ldots,2t-1\}$. Since by induction $u_{i-2}u_{i-1}$ and $u_{i-1}v_{i-1}$ belong to the same part of $(F_k,F_\ell)$, we deduce that $u_{i-1}u_i$ and the arc incident to $u_{i-1}$ contained in the attached $(k,\ell,-2)$-in-forcer are contained in the same part of $(F_k,F_\ell)$. By the definition of $(k,\ell,-2)$-in-forcers, we obtain that $u_{i-1}u_i$ is the last arc of a directed path either of length $k-1$ in $F_k$ or of length $\ell-1$ in $F_\ell$. Since $u_{i+1}$ has out-degree $2$, we deduce that $u_{i-1}u_i$ and $u_{i}v_i$ belong to the same part, and that $u_iv_i$ is either the last arc of a directed path of length $k$ in $F_k$ or the last arc of a directed path of length $\ell$ in $F_\ell$. This shows $(b)$.

    Further observe that the described $(k,\ell)$-decomposition indeed yields the decompositions required in $(c)$ and $(d)$ when choosing appropriate decompositions of the $(k,\ell,-2)$-in-forcers. An illustration can be found in Figure~\ref{fig:vertex_gadget_k_l_decomposition}. Note that the roles of $F_k$ and $F_\ell$ are symmetric.
\end{proof}

We now give a similar result for the case that $\min\{k,\ell\}=2$.

\begin{lemma}
    \label{lemma:k_2_t_gadget}
    For all positive integers $k,t$ with $k \geq 2$, there exists a $(k,2,t)$-variable gadget.
\end{lemma}
\begin{proof}
    We describe a $(k,2,t)$-variable gadget for fixed positive integers $k$ and $t$ with $k\geq 2$.
    First let $V(D)$ contain a vertex $v_i^{j}$ for all $i\in [t]$ and $j\in[7]$ with $(i,j)\notin \{(t,6),(t,7)\}$. We further let $V(D)$ contain one more vertex $w$. Next, for every $i \in [t]$, we let $A(D)$ contain the arcs $v_i^2v_i^1,v_i^2v_i^3, v_i^4v_i^3$, and $v_i^4v_i^5$ and for every $i \in [t-1]$, we let $A(D)$ contain the arcs $v_i^1v_i^6,v_i^7v_i^6$, and $v_i^6v_{i+1}^1$. Next, we add the arc $wv_1^1$. We now complete the construction of $D$ dependent on $k$. If $k=2$, we further add two in-arcs to $w$ and to $v_i^4$ for every $i \in [t]$.
   If $k\geq 3$, we further attach a $(k,2,k-2)$-in-forcer to $v_i^2$ for all $i \in [t]$ and to $v_i^7$ for all $i \in [t-1]$ and we attach a $(k,2,k-1)$-in-forcer and an in-arc to $w$ and to $v_i^4$ for every $i \in [t]$. Observe that these forcers exist by Lemma \ref{lemma:2_k_l_forcer}.
    Finally, for $i\in[t]$, we set $a_i=v_{i}^4v_{i}^5$. This finishes the description of $D$. For an illustration, see Figure~\ref{vark23_neu}.

\begin{figure}[hbt!]\begin{center}
    \begin{tikzpicture}[thick,scale=1, every node/.style={transform shape}]         
    \tikzset{vertex/.style = {circle,fill=black,minimum size=5pt, inner sep=0pt}}
    \tikzset{littlevertex/.style = {circle,fill=black,minimum size=4pt, inner sep=0pt}}
    \tikzset{edge/.style = {->,> = latex'}}
        \foreach \i in {1,...,3}{
            \node[vertex, label=below:$v_\i^1$] (v\i1) at (\i*4,0) {};
            \node[vertex, label=right:$v_\i^2$] (v\i2) at (\i*4,1) {};
            \node[vertex, label=right:$v_\i^3$] (v\i3) at (\i*4,2) {};
            \node[vertex, label=right:$v_\i^4$] (v\i4) at (\i*4,3) {};
            \node[vertex, label=right:$v_\i^5$] (v\i5) at (\i*4,4.5) {};
            \node[] (a\i) at (\i*4-0.3,3.75) {$a_\i$};
        
            \draw[edge] (v\i2) to (v\i1) {};
            \draw[edge] (v\i2) to (v\i3) {};
            \draw[edge] (v\i4) to (v\i3) {};            
            \draw[edge] (v\i4) to (v\i5) {};

            \draw[] (v\i2)  -- (\i*4-0.25,1+0.25) -- (\i*4-0.5,1) -- (\i*4-0.25,1-0.25) -- (v\i2);
            \draw[edge] (\i*4-0.4,3-0.4) to (v\i4);
            \draw[] (v\i4) -- (\i*4-0.5,3) -- (\i*4-0.25,3+0.25) -- (v\i4);
        }
        \foreach \i in {1,...,2}{
            \pgfmathtruncatemacro{\j}{\i +1}
            \node[vertex, label=below:$v_\i^6$] (v\i6) at (\i*4+2,0) {};
            \node[vertex, label=right:$v_\i^7$] (v\i7) at (\i*4+2,1) {};
            \draw[edge] (v\i1) to (v\i6) {};
            \draw[edge] (v\i7) to (v\i6) {};
            \draw[edge] (v\i6) to (v\j1) {};

            \draw[] (v\i7)  -- (\i*4+2-0.25,1+0.25) -- (\i*4+2-0.5,1) -- (\i*4+2-0.25,1-0.25) -- (v\i7);
        }
        \node[vertex, label=below:$w$] (w) at (2,0) {};
        \draw[edge] (w) to (v11) {};
        \draw[edge] (2-0.4,-0.4) to (w);
        \draw[] (w) -- (2-0.5,0) -- (2-0.25,0.25) -- (w);
    \end{tikzpicture}
    \caption{An illustration of a $(k,2,3)$-variable gadget for some $k \geq 2$. When $k\geq 3$, the triangles indicate $(k,2,k-1)$-in-forcers and the squares indicate $(k,2,k-2)$-in-forcers. When $k=2$, the triangles indicate in-arcs and the squares are deleted. }
    \label{vark23_neu}
\end{center}
\end{figure}

    Clearly, $D$ satisfies $(a)$. For $(b)$, let $(F_k,F_2)$ be a $(k,2)$-decomposition of $D$. First suppose for the sake of a contradiction that there is an index $i \in [t-1]$ such that $a_i \in  A(F_2)$ and $a_{i+1}\in A(F_k)$. As $F_k$ and $F_2$ are directed linear forests, we obtain that $v_i^2v_i^1 \in A(F_k)$ and $v_{i+1}^2v_{i+1}^1\in A(F_2)$. As $F_2$ is a directed linear forest, we obtain that $v_i^6v_{i+1}^1\in A(F_k)$. If $v_i^1v_i^6\in A(F_k)$ and $k=2$, we obtain that $v_i^2v_i^1v_i^6v_{i+1}^1$ is a path of length 3 in $F_k$, a contradiction. If $v_i^1v_i^6\in A(F_k)$ and $k\geq 3$, we obtain that  the path which is obtained from concatenating a path in $F_k$ of length $k-2$ fully contained in the $(k,2,k-2)$-in-forcer incident to $v_i^2$ whose last vertex is $v_i^2$ with $v_i^2v_i^1v_i^6v_{i+1}^1$ is a path of length $k+1$ in $F_k$, a contradiction. Hence $v_i^1v_i^6\in A(F_2)$. As $F_k$ is a directed linear forest, we obtain that $xv_i^1 \in A(F_2)$ where $x=v_{i-1}^6$ if $i\geq 2$ and $x=w$ otherwise. Further, as $d_D^-(x)=2$ by construction and $F_k$ is a linear forest, there exists a vertex $y \in N^-_D(x)$ with $yx \in A(F_2)$. It follows that $yxv_i^{1}v_i^6$ is a path of length 3 in $F_2$, a contradiction.

    Now suppose for the sake of a contradiction that there is an index $i \in [t-1]$ such that $a_i \in  A(F_k)$ and $a_{i+1}\in A(F_2)$. By the above, we have $a_j \in  A(F_k)$ for all $j \in [i]$. As $F_k$ and $F_2$ are directed linear forests, we obtain that $v_j^2v_j^1 \in A(F_2)$ for all $j \in [i]$. As $F_2$ is a directed linear forest, we obtain $wv_1^1 \in A(F_k)$ and $v_j^6v_{j+1}^1\in A(F_k)$ for all $j \in [i-1]$. As $F_k$ is a linear forest, we obtain $v_i^6v_{i+1}^1\in A(F_k)$.
    We next show that $v_j^1v_j^6 \in A(F_k)$ for all $j \in [i]$. Suppose otherwise and let $j_0$ be the largest integer with $j_0 \leq i$ and $v_{j_0}^1v_{j_0}^6 \in A(F_2)$. 
    If $j_0=i$, we obtain that $v_i^2v_i^1v_i^6v_{i+1}^1$ is a directed path of length 3 in $F_2$, a contradiction, so assume $j_0<i$.
    As $F_2$ is a directed linear forest, we obtain that $v_{j_0}^7v_{j_0}^6\in A(F_k)$. Since $j_0<i$, we have $v_{j_0}^6v_{j_0+1}^1 \in A(F_k)$, and by maximality of $j_0$ we have $v_{j_0+1}^1v_{j_0+1}^6 \in A(F_k)$. 
    Further, by construction, we have that $v_{j_0}^7$ is the last vertex of a path of length $k-2$ in $F_k$. The concatenation of this path with $v_{j_0}^7v_{j_0}^6v_{j_0+1}^1v_{j_0+1}^6$ is a directed path of length $k+1$ in $F_k$, a contradiction.
    In particular, we obtain that $v_1^1v_1^6\in A(F_k)$. By definition, there exists a path of length $k-1$ in $F_k$ whose last vertex is $w$ and that is fully contained in the $(k,2,k-1)$-in-forcer attached at $w$. Concatenating this path with $wv_1^1v_1^6$ yields a path of length $k+1$ in $F_k$, a contradiction. 

    We have shown that $A(F_k) \cap \{a_1,\dots,a_t\} \in \{\{a_1,\dots,a_t\},\emptyset\}$. Let $i\in[t]$. We will show that $a_i$ is either the last arc of a directed path of length $k$ in $F_k$ or the last arc of a directed path of length $2$ in $F_2$. Assume first that $a_i \in A(F_2)$, since $v_i^4$ has two entering arcs, it has an in-neighbour $x$ such that $xv_i^4 \in A(F_2)$. Thus, $xv_i^4v_i^5$ is a directed path of length 2 in $F_2$. Now assume that $a_i\in A(F_k)$. If $k=2$, $F_2$ and $F_k$ play a symmetric role, so assume $k\geq 3$. The concatenation of the path of $F_k$ of length $k-1$ fully contained in the $(k,2,k-1)$-in-forcer attached in $v_i^4$ and $a_i$ is a path of length $k$ in $F_k$, the last arc of which is $a_i$.
    This yields $(b)$.
    
     For $(c)$, we define a decomposition $(F_k,F_2)$ of $D$ in the following way: If $k \geq 3$, for every $(k,2,k-1)$-in-forcer (respectively $(k,2,k-2)$-in-forcer), we choose a $(k,\ell)$-decomposition such that the tip arc of this forcer is not the last arc of a path of length $k$ (respectively $k-1$) fully contained in the forcer. Next, we let all attached in-arcs be contained in $A(F_2)$. If $k=2$, for every $v \in \{v_i^4 \mid i \in [t]\}\cup w$, we let each of $A(F_k)$ and $A(F_2)$ contain one of the arcs entering $v$.  We then extend this to $(F_k,F_2)$ in the following way: for every $i \in [t]$, we let $A(F_k)$ contain $a_i$ and $v_i^2v_i^3$, for every $i \in [t-1]$, we let $A(F_k)$ contain $v_i^6v_{i+1}^1$ and $v_i^7v_i^6$, and we let $A(F_k)$ contain $wv_1^1$. We then set $A(F_2)=A(D)-A(F_k)$.  This finishes the description of $(F_k,F_2)$. For an illustration, see Figure~\ref{vark24_neu}.  
     
     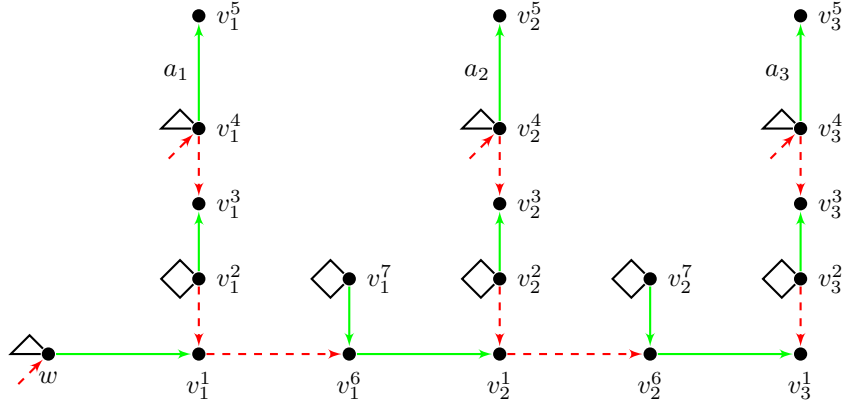
\begin{figure}[hbt!]\begin{center}
        \begin{tikzpicture}[thick,scale=1, every node/.style={transform shape}]    
        \tikzset{vertex/.style = {circle,fill=black,minimum size=5pt, inner sep=0pt}}
        \tikzset{littlevertex/.style = {circle,fill=black,minimum size=4pt, inner sep=0pt}}
        \tikzset{edge/.style = {->,> = latex'}}
        \foreach \i in {1,...,3}{
            \node[vertex, label=below:$v_\i^1$] (v\i1) at (\i*4,0) {};
            \node[vertex, label=right:$v_\i^2$] (v\i2) at (\i*4,1) {};
            \node[vertex, label=right:$v_\i^3$] (v\i3) at (\i*4,2) {};
            \node[vertex, label=right:$v_\i^4$] (v\i4) at (\i*4,3) {};
            \node[vertex, label=right:$v_\i^5$] (v\i5) at (\i*4,4.5) {};
            \node[] (a\i) at (\i*4-0.3,3.75) {$a_\i$};
        
            \draw[edge,dashed,red] (v\i2) to (v\i1) {};
            \draw[edge,green] (v\i2) to (v\i3) {};
            \draw[edge,dashed,red] (v\i4) to (v\i3) {};            
            \draw[edge,green] (v\i4) to (v\i5) {};

            \draw[] (v\i2)  -- (\i*4-0.25,1+0.25) -- (\i*4-0.5,1) -- (\i*4-0.25,1-0.25) -- (v\i2);
            \draw[edge,dashed,red] (\i*4-0.4,3-0.4) to (v\i4);
            \draw[] (v\i4) -- (\i*4-0.5,3) -- (\i*4-0.25,3+0.25) -- (v\i4);
        }
        \foreach \i in {1,...,2}{
            \pgfmathtruncatemacro{\j}{\i +1}
            \node[vertex, label=below:$v_\i^6$] (v\i6) at (\i*4+2,0) {};
            \node[vertex, label=right:$v_\i^7$] (v\i7) at (\i*4+2,1) {};
            \draw[edge,dashed,red] (v\i1) to (v\i6) {};
            \draw[edge,green] (v\i7) to (v\i6) {};
            \draw[edge,green] (v\i6) to (v\j1) {};

            \draw[] (v\i7)  -- (\i*4+2-0.25,1+0.25) -- (\i*4+2-0.5,1) -- (\i*4+2-0.25,1-0.25) -- (v\i7);
        }
        \node[vertex, label=below:$w$] (w) at (2,0) {};
        \draw[edge,green] (w) to (v11) {};
        \draw[edge,dashed,red] (2-0.4,-0.4) to (w);
        \draw[] (w) -- (2-0.5,0) -- (2-0.25,0.25) -- (w);
    \end{tikzpicture}
      \caption{An illustration of the decomposition $(F_k,F_2)$. The dashed red arcs are in $A(F_2)$ and the solid green arcs are in $A(F_k)$. When $k=2$, the triangles are replaced by solid green arcs.}\label{vark24_neu}
      \end{center}
      \end{figure}
  This shows $(c)$. When $k=2$, it also shows $(d)$, so we now assume $k\geq 3$.

     For $(d)$, we define a decomposition $(F'_k,F'_2)$ of $D$ in the following way: Again, for every $(k,2,k-1)$-in-forcer (respectively $(k,2,k-2)$-in-forcer), we choose a $(k,2)$-decomposition such that the tip arc of this forcer is not the last arc of a path of length $k$ (respectively $k-1$) fully contained in the forcer. We then extend this to $(F'_k,F'_2)$ in the following way: for every $i \in [t]$, we let $A(F'_2)$ contain $a_i$ and $v_i^2v_i^3$, for every $i \in [t-1]$, we let $A(F'_2)$ contain $v_i^6v_{i+1}^1$ and $v_i^7v_i^6$, and we let $A(F'_2)$ contain $wv_1^1$. We then set $A(F'_k)=A(D)-A(F'_2)$.  This finishes the description of $(F'_k,F'_2)$. For an illustration, see Figure~\ref{vark25_neu}. This shows $(d)$.
     \begin{figure}[hbt!]\begin{center}
     \begin{tikzpicture}[thick,scale=1, every node/.style={transform shape}]    
        \tikzset{vertex/.style = {circle,fill=black,minimum size=5pt, inner sep=0pt}}
        \tikzset{littlevertex/.style = {circle,fill=black,minimum size=4pt, inner sep=0pt}}
        \tikzset{edge/.style = {->,> = latex'}}
        \foreach \i in {1,...,3}{
            \node[vertex, label=below:$v_\i^1$] (v\i1) at (\i*4,0) {};
            \node[vertex, label=right:$v_\i^2$] (v\i2) at (\i*4,1) {};
            \node[vertex, label=right:$v_\i^3$] (v\i3) at (\i*4,2) {};
            \node[vertex, label=right:$v_\i^4$] (v\i4) at (\i*4,3) {};
            \node[vertex, label=right:$v_\i^5$] (v\i5) at (\i*4,4.5) {};
            \node[] (a\i) at (\i*4-0.3,3.75) {$a_\i$};
        
            \draw[edge,green] (v\i2) to (v\i1) {};
            \draw[edge,dashed,red] (v\i2) to (v\i3) {};
            \draw[edge,green] (v\i4) to (v\i3) {};            
            \draw[edge,dashed,red] (v\i4) to (v\i5) {};

            \draw[] (v\i2)  -- (\i*4-0.25,1+0.25) -- (\i*4-0.5,1) -- (\i*4-0.25,1-0.25) -- (v\i2);
            \draw[edge,dashed,red] (\i*4-0.4,3-0.4) to (v\i4);
            \draw[] (v\i4) -- (\i*4-0.5,3) -- (\i*4-0.25,3+0.25) -- (v\i4);
        }
        \foreach \i in {1,...,2}{
            \pgfmathtruncatemacro{\j}{\i +1}
            \node[vertex, label=below:$v_\i^6$] (v\i6) at (\i*4+2,0) {};
            \node[vertex, label=right:$v_\i^7$] (v\i7) at (\i*4+2,1) {};
            \draw[edge,green] (v\i1) to (v\i6) {};
            \draw[edge,dashed,red] (v\i7) to (v\i6) {};
            \draw[edge,dashed,red] (v\i6) to (v\j1) {};

            \draw[] (v\i7)  -- (\i*4+2-0.25,1+0.25) -- (\i*4+2-0.5,1) -- (\i*4+2-0.25,1-0.25) -- (v\i7);
        }
        \node[vertex, label=below:$w$] (w) at (2,0) {};
        \draw[edge,dashed,red] (w) to (v11) {};
        \draw[edge,dashed,red] (2-0.4,-0.4) to (w);
        \draw[] (w) -- (2-0.5,0) -- (2-0.25,0.25) -- (w);
    \end{tikzpicture}
  \caption{An illustration of the decomposition $(F'_k,F'_2)$. The dashed red arcs are in $A(F'_2)$ and the solid green arcs are in $A(F'_k)$.}\label{vark25_neu}
  \end{center}
  \end{figure}
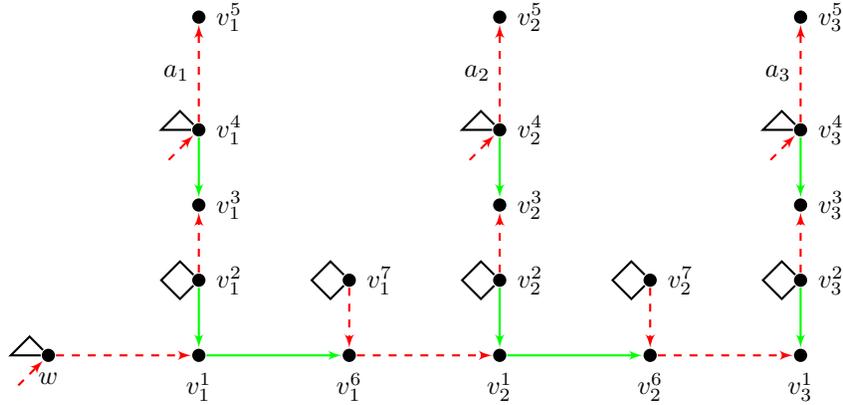 
\end{proof}
\subsubsection{Clause gadgets}\label{clause}

We now describe the clause gadgets we need.
 For two integers $k$ and $\ell$ with $\min\{k,\ell\}\geq 2$, a {\bf $(k,\ell)$-clause gadget} is a digraph $D$ together with three special arcs $a_1,a_2,a_3 \in A(D)$ and three special vertices $t_1,t_2,t_3 \in V(D)$ satisfying the following properties:
\begin{enumerate}[$(a)$]
    \item for $i \in [3]$, we have that $a_i$ is the only arc incident to $t_i$,
    \item for every $(k,\ell)$-decomposition $(F_k,F_\ell)$ of $D$, we have $A(F_k)\cap\{a_1,a_2,a_3\}\neq \emptyset$ and $A(F_\ell)\cap\{a_1,a_2,a_3\}\neq \emptyset$, and
    \item for every nonempty $Z \subsetneq [3]$ there is a $(k,\ell)$-decomposition $(F_k,F_\ell)$ of $D$ with $\{i \in [3] \mid a_i \in A(F_k)\}=Z$.
\end{enumerate}

Again, we first prove the existence of clause gadgets when $\min\{k,\ell\}\geq 3$.
\begin{lemma}\label{clause3}
    For all integers $k,\ell$ with $\min \{k,\ell\}\geq 3$, there exists a $(k,\ell)$-clause gadget.
\end{lemma}
\begin{proof}
    Let $D$ be the digraph on six vertices $y_1,y_2,y_3,t_1,t_2,t_3$ made of the directed triangle $y_1y_2y_3y_1$ and the arcs $a_i=t_iy_i$ for $i\in [3]$, see Figure~\ref{fig:clause_gadget_k_l} for an illustration. 
    \begin{figure}[hbt!]
        \begin{center}	
              \begin{tikzpicture}[thick,scale=1, every node/.style={transform shape}]
                
                \tikzset{vertex/.style = {circle,fill=black,minimum size=5pt, inner sep=0pt}}
                \tikzset{littlevertex/.style = {circle,fill=black,minimum size=4pt, inner sep=0pt}}
                \tikzset{edge/.style = {->,> = latex'}}

                \foreach \i in {1,2}{
                    \node[vertex, label=left:$y_\i$] (y\i) at (\i*120:1) {};
                    \node[vertex, label=left:$t_\i$] (t\i) at (\i*120:2) {};
                    \draw[edge,green] (t\i) to (y\i) {};
                }
                \node[vertex, label=below:$y_3$] (y3) at (0:1) {};
                \node[vertex, label=below:$t_3$] (t3) at (0:2) {};
                \draw[edge, dashed, red] (t3) to (y3) {};
                
                \node[] (a1) at (110:1.5) {$a_1$};
                \node[] (a2) at (-110:1.5) {$a_2$};
                \node[] (a3) at (1.5,0.3) {$a_3$};
                
                \draw[edge, dashed, red] (y1) to (y2) {};
                \draw[edge,green] (y2) to (y3) {};
                \draw[edge, dashed, red] (y3) to (y1) {};
              \end{tikzpicture}
          \caption{An illustration of the $(k,\ell)$-decomposition $(F_k,F_\ell)$ of the $(k,\ell)$-clause gadget when $\min\{k,\ell\} \geq 3$. The dashed red arcs are in $A(F_k)$ and the solid green arcs are in $A(F_\ell)$.}
          \label{fig:clause_gadget_k_l}
        \end{center}
    \end{figure}
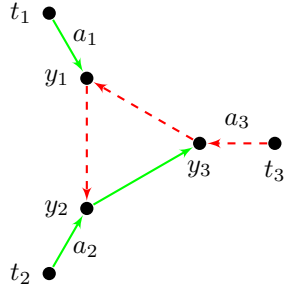
    Clearly, $(a)$ holds. In order to prove $(b)$, let $(F_k,F_\ell)$ be a $(k,\ell)$-decomposition of $D$ and suppose for the sake of a contradiction that $A(F_k)\cap\{a_1,a_2,a_3\}=\emptyset$. As $F_\ell$ is a directed linear forest, we obtain that $\{y_1y_2,y_2y_3,y_3y_1\}\subseteq A(F_k)$. This contradicts $F_k$ being a directed linear forest. A similar argument shows that $A(F_\ell)\cap\{a_1,a_2,a_3\}\neq \emptyset$. This yields $(b)$. For $(c)$, by symmetry, it suffices to prove the statement for $Z\in \{\{1,2\},\{3\}\}$. Let $(F_k,F_\ell)$ be defined by $A(F_k)=\{t_1y_1,t_2y_2,y_2y_3\}$ and $A(F_\ell)=A(D)-A(F_k)$. Then $(F_k,F_\ell)$ has the desired properties for $Z=\{1,2\}$ and $(F_\ell,F_k)$ has the desired properties for $Z=\{3\}$. This proves $(c)$.  
\end{proof}

We now prove a similar result for the case that $\min\{k,\ell\}=2$.

\begin{lemma}\label{clause2}
    For every integer $k \geq 2$, there exists a $(k,2)$-clause gadget.
\end{lemma}
\begin{proof}
    We describe a $(k,2)$-clause gadget $D$ for some $k\geq 2$. We first let $V(D)$ contain some vertices $t_1,t_2,t_3,v_1,\ldots,v_5$ and we let $A(D)$ contain the arcs $v_1t_1,v_2t_2,t_3v_3,v_4v_1,v_4v_2,v_5v_3$, and $v_5v_4$. If $k \geq 3$, we further attach a $(k,2,k-2)$-in-forcer to $v_5$. Observe that this forcer exists by Lemma \ref{lemma:2_k_l_forcer}. We finally set $a_1=v_1t_1,a_2=v_2t_2$, and $a_3=t_3v_3$. This finishes the description of $(D,\{a_1,a_2,a_3\})$. For an illustration, see Figure~\ref{clause_2}.
     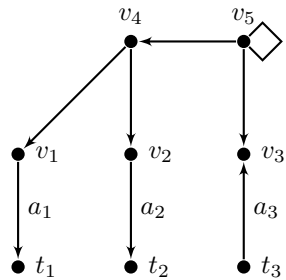
\begin{figure}[hbt!]\begin{center}
     \begin{tikzpicture}[thick,scale=1, every node/.style={transform shape}]
        \tikzset{vertex/.style = {circle,fill=black,minimum size=5pt, inner sep=0pt}}
        \tikzset{littlevertex/.style = {circle,fill=black,minimum size=4pt, inner sep=0pt}}
        \tikzset{edge/.style = {<-,> = latex'}}
        
        \foreach \i in {1,...,3}{
        \node[vertex, label=right:$t_\i$] (x\i) at (\i*1.5,0) {};
        \node[vertex, label=right:$v_\i$] (y\i) at (\i*1.5,1*1.5) {};
        \node[] (a\i) at (\i*1.5+0.3,0.5*1.5) {$a_\i$};
        }
        \node[vertex, label=above:$v_4$] (z12) at (2*1.5,2*1.5) {};
        \node[vertex, label=above:$v_5$] (z3) at (3*1.5,2*1.5) {};
        
        \draw[edge] (x1) to (y1) {};
        \draw[edge] (x2) to (y2) {};
        \draw[edge] (y3) to (x3) {};
        \draw[edge] (y1) to (z12) {};
        \draw[edge] (y2) to (z12) {};
        \draw[edge] (y3) to (z3) {};
        \draw[edge] (z12) to (z3) {};
        \draw[] (z3) -- (3*1.5+0.25,2*1.5-0.25) -- (3*1.5+0.5,2*1.5) -- (3*1.5+0.25,2*1.5+0.25) -- (z3);

     \end{tikzpicture}
  \caption{An illustration of $D$ where the square marks a $(k,2,k-2)$-in-forcer which does not exist if $k=2$.}\label{clause_2}
  \end{center}
  \end{figure}
    Clearly, $D$ satisfies $(a)$. 

    For $(b)$, let $(F_k,F_2)$ be a $(k,2)$-decomposition of $D$ and first suppose for the sake of a contradiction that $\{a_1,a_2,a_3\}\subseteq A(F_k)$. As $F_k$ and $F_2$ are directed linear forests, we obtain that $v_5v_4 \in A(F_k)$ and one of $v_4v_1$ and $v_4v_2$, say $v_4v_1$, is contained in $A(F_k)$. If $k=2$, we obtain that $v_5v_4v_1t_1$ is a directed path of length 3 in $F_k$, a contradiction. If $k \geq 3$, then the concatenation of a path in $F_k$ of length $k-2$ fully contained in the $(k,2,k-2)$-in-forcer attached to $v_5$ whose last vertex is $v_5$ with $v_5v_4v_1t_1$ is a directed path of length $k+1$ in $F_k$, a contradiction. If $k=2$, this yields $(b)$ by symmetry. Now suppose that $k \geq 3$ and, for the sake of a contradiction, that $\{a_1,a_2,a_3\}\subseteq A(F_2)$. As $F_k$ and $F_2$ are directed linear forests, we obtain that $v_5v_4 \in A(F_2)$ and one of $v_4v_1$ and $v_4v_2$, say $v_4v_1$, is contained in $A(F_2)$. Then $v_5v_4v_1t_1$ is a directed path of length $3$ in $F_2$, a contradiction. This yields $(b)$.

    For $(c)$, if $k \geq 3$, we choose a decomposition of the $(k,2,k-2)$-in-forcer attached to $v_5$ in which $v_5$ is not the last vertex of a path of length $k-1$. A small case analysis shows that for any nonempty $S \subsetneq [3]$, this decomposition can be extended into a decomposition $(F_k,F_2)$ of $D$ with $\{i \in [3] \mid a_i \in A(F_k)\}=S$. All these cases are illustrated in Figure~\ref{clause_2_dec}.
\end{proof}
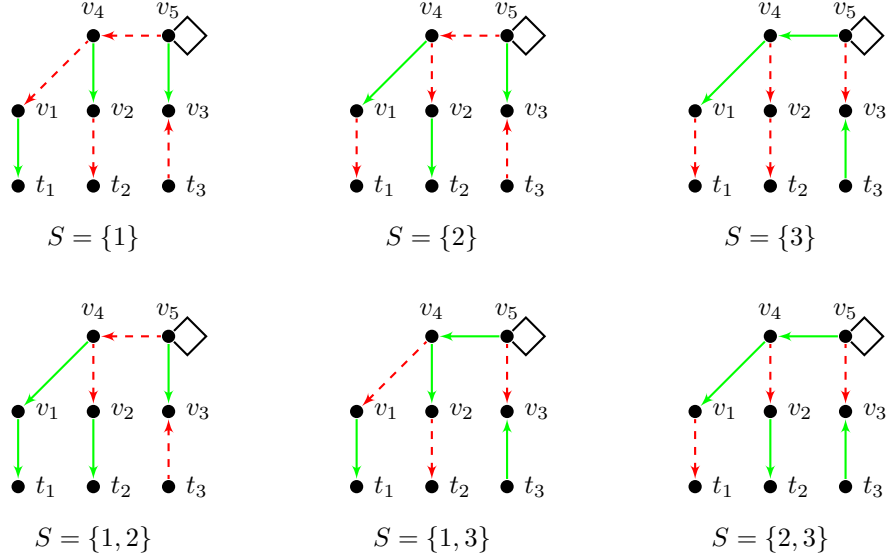
\begin{figure}[hbt!]\begin{center}
  \begin{tikzpicture}[thick,scale=1, every node/.style={transform shape}]
        \tikzset{vertex/.style = {circle,fill=black,minimum size=5pt, inner sep=0pt}}
        \tikzset{littlevertex/.style = {circle,fill=black,minimum size=4pt, inner sep=0pt}}
        \tikzset{edge/.style = {<-,> = latex'}}

        \foreach \i in {1,...,3}{
        \node[vertex, label=right:$t_\i$] (x\i) at (\i,0) {};
        \node[vertex, label=right:$v_\i$] (y\i) at (\i,1) {};
        }
        \node[vertex, label=above:$v_4$] (z1) at (2,2) {};
        \node[vertex, label=above:$v_5$] (z2) at (3,2) {};
        \draw[edge,green] (x1) to (y1) {};
        \draw[edge,dashed, red] (x2) to (y2) {};
        \draw[edge,dashed, red] (y3) to (x3) {};
        \draw[edge,dashed, red] (y1) to (z1) {};
        \draw[edge,green] (y2) to (z1) {};
        \draw[edge,green] (y3) to (z2) {};
        \draw[edge,dashed, red] (z1) to (z2) {};
        \draw[] (z2) -- (3+0.25,2-0.25) -- (3+0.5,2) -- (3+0.25,2+0.25) -- (z2);
        \node[] (S) at (2, -0.7){ $S=\{1\}$};

        \begin{scope}[xshift=4.5cm]
        \foreach \i in {1,...,3}{
        \node[vertex, label=right:$t_\i$] (x\i) at (\i,0) {};
        \node[vertex, label=right:$v_\i$] (y\i) at (\i,1) {};
        }
        \node[vertex, label=above:$v_4$] (z1) at (2,2) {};
        \node[vertex, label=above:$v_5$] (z2) at (3,2) {};
        \draw[edge,dashed, red] (x1) to (y1) {};
        \draw[edge, green] (x2) to (y2) {};
        \draw[edge,dashed, red] (y3) to (x3) {};
        \draw[edge, green] (y1) to (z1) {};
        \draw[edge,dashed, red] (y2) to (z1) {};
        \draw[edge, green] (y3) to (z2) {};
        \draw[edge,dashed, red] (z1) to (z2) {};
        \draw[] (z2) -- (3+0.25,2-0.25) -- (3+0.5,2) -- (3+0.25,2+0.25) -- (z2);
        \node[] (S) at (2, -0.7){ $S=\{2\}$};
        \end{scope}

        \begin{scope}[xshift=9cm]
        \foreach \i in {1,...,3}{
        \node[vertex, label=right:$t_\i$] (x\i) at (\i,0) {};
        \node[vertex, label=right:$v_\i$] (y\i) at (\i,1) {};
        }
        \node[vertex, label=above:$v_4$] (z1) at (2,2) {};
        \node[vertex, label=above:$v_5$] (z2) at (3,2) {};
        \draw[edge,dashed, red] (x1) to (y1) {};
        \draw[edge,dashed, red] (x2) to (y2) {};
        \draw[edge, green] (y3) to (x3) {};
        \draw[edge, green] (y1) to (z1) {};
        \draw[edge,dashed, red] (y2) to (z1) {};
        \draw[edge,dashed, red] (y3) to (z2) {};
        \draw[edge, green] (z1) to (z2) {};
        \draw[] (z2) -- (3+0.25,2-0.25) -- (3+0.5,2) -- (3+0.25,2+0.25) -- (z2);
        \node[] (S) at (2, -0.7){ $S=\{3\}$};
        \end{scope}

        \begin{scope}[yshift=-4cm]
        \foreach \i in {1,...,3}{
        \node[vertex, label=right:$t_\i$] (x\i) at (\i,0) {};
        \node[vertex, label=right:$v_\i$] (y\i) at (\i,1) {};
        }
        \node[vertex, label=above:$v_4$] (z1) at (2,2) {};
        \node[vertex, label=above:$v_5$] (z2) at (3,2) {};
        \draw[edge, green] (x1) to (y1) {};
        \draw[edge, green] (x2) to (y2) {};
        \draw[edge,dashed, red] (y3) to (x3) {};
        \draw[edge, green] (y1) to (z1) {};
        \draw[edge,dashed, red] (y2) to (z1) {};
        \draw[edge, green] (y3) to (z2) {};
        \draw[edge,dashed, red] (z1) to (z2) {};
        \draw[] (z2) -- (3+0.25,2-0.25) -- (3+0.5,2) -- (3+0.25,2+0.25) -- (z2);
        \node[] (S) at (2, -0.7){ $S=\{1,2\}$};
        \end{scope}

        \begin{scope}[xshift=4.5cm,yshift=-4cm]
        \foreach \i in {1,...,3}{
        \node[vertex, label=right:$t_\i$] (x\i) at (\i,0) {};
        \node[vertex, label=right:$v_\i$] (y\i) at (\i,1) {};
        }
        \node[vertex, label=above:$v_4$] (z1) at (2,2) {};
        \node[vertex, label=above:$v_5$] (z2) at (3,2) {};
        \draw[edge, green] (x1) to (y1) {};
        \draw[edge,dashed, red] (x2) to (y2) {};
        \draw[edge, green] (y3) to (x3) {};
        \draw[edge,dashed, red] (y1) to (z1) {};
        \draw[edge, green] (y2) to (z1) {};
        \draw[edge,dashed, red] (y3) to (z2) {};
        \draw[edge, green] (z1) to (z2) {};
        \draw[] (z2) -- (3+0.25,2-0.25) -- (3+0.5,2) -- (3+0.25,2+0.25) -- (z2);
        \node[] (S) at (2, -0.7){ $S=\{1,3\}$};
        \end{scope}

        \begin{scope}[xshift=9cm,yshift=-4cm]
        \foreach \i in {1,...,3}{
        \node[vertex, label=right:$t_\i$] (x\i) at (\i,0) {};
        \node[vertex, label=right:$v_\i$] (y\i) at (\i,1) {};
        }
        \node[vertex, label=above:$v_4$] (z1) at (2,2) {};
        \node[vertex, label=above:$v_5$] (z2) at (3,2) {};
        \draw[edge,dashed, red] (x1) to (y1) {};
        \draw[edge, green] (x2) to (y2) {};
        \draw[edge, green] (y3) to (x3) {};
        \draw[edge, green] (y1) to (z1) {};
        \draw[edge,dashed, red] (y2) to (z1) {};
        \draw[edge,dashed, red] (y3) to (z2) {};
        \draw[edge, green] (z1) to (z2) {};
        \draw[] (z2) -- (3+0.25,2-0.25) -- (3+0.5,2) -- (3+0.25,2+0.25) -- (z2);
        \node[] (S) at (2, -0.7){ $S=\{2,3\}$};
        \end{scope}

     \end{tikzpicture}
  \caption{An illustration of the possible decompositions $(F_k,F_2)$ of $D$ for all different choices of $S$. The dashed red arcs are in $A(F_2)$ and the solid green arcs are in $A(F_k)$.}\label{clause_2_dec}
  \end{center}
\end{figure}
  
\subsubsection{The main reduction}\label{hardness}

We are finally ready to prove the following main result of Section~\ref{sec:2_2_case} making use of the gadgets constructed in Sections~\ref{variable} and~\ref{clause}.

\begin{theorem}
    For all integers $k, \ell$ with $\min\{k,\ell\} \geq 2$, $(k,\ell)$-BDLFD is NP-complete.
\end{theorem}
\begin{proof}
    Let us fix a pair of integers $k$ and $\ell$ with $\min\{k,\ell\} \geq 2$. Note that $(k,\ell)$-BDLFD is clearly in NP. We will show that it is NP-complete through a reduction from ME-$1$-SAT, which is NP-complete by Proposition~\ref{prop:MEkSAT_NP_hard}.

    Let $(X,\mathcal{C})$ be an instance of ME-$1$-SAT.
    We now describe a digraph $D$. For some $x \in X$, let $C_1,\ldots,C_q$ be the clauses containing $x$. We let $D$ contain a $(k,\ell,q)$-variable gadget $(H_x,\{a_1,\ldots,a_q\})$ which exists by Lemmas~\ref{lemma:k_l_t_gadget} and~\ref{lemma:k_2_t_gadget}. Further, for $i \in [q]$, we let $u_{(x,C_i)}$ denote the head of $a_i$. We do this for every $x \in X$. For some $C \in \mathcal{C}$, let $x_1,x_2$, and $x_3$ be the variables contained in $C$.  We let $D$ contain a clause-gadget $(H_C,\{a_1,a_2,a_3\})$ which exists by Lemmas~\ref{clause2} and~\ref{clause3}. Further, for $i \in [3]$, we let $t_{(x_i,C)}$ denote the vertex incident to $a_i$ which is of degree 1 in $H_C$.  We do this for every $C \in \mathcal{C}$.

    Finally, we obtain $D$ by identifying $u_{(x,C)}$ and $t_{(x,C)}$ for all $x \in X$ and $C \in \mathcal{C}$ with $x \in C$.
    We now show that $D$ is a positive instance of $(k,\ell)$-BDLFD if and only if $(X,\mathcal{C})$ is a positive instance of ME-$1$-SAT.

    First suppose that $(X,\mathcal{C})$ is a positive instance of ME-$1$-SAT, so there exists a satisfying assignment $\phi \colon X \rightarrow \{\true,\false\}$ for $(X,\mathcal{C})$. Let $x \in X$. If $\phi(x)=\true$, we choose a $(k,\ell)$-decomposition $(F_k^x,F_\ell^x)$ of $H_x$ with $\{a_1,\ldots,a_q\}\subseteq F_k$ and if $\phi(x)=\false$, we choose a $(k,\ell)$-decomposition $(F_k^x,F_\ell^x)$ of $H_x$ with $\{a_1,\ldots,a_q\}\subseteq F_\ell$. We do this for every $x \in X$. Observe that such a decomposition exists by the definition of $H_x$. Now consider some $C \in \mathcal{C}$ and let $x_1,x_2$, and $x_3$ be the variables contained in $C$ in the order they were used when constructing $H_C$. Let $Z=\{i \in [3] \mid \phi(x_i)=\false\}$. As $\phi$ is a satisfying assignment for $(X,\mathcal{C})$, we have that $Z$ is a nonempty, strict subset of $[3]$. We now choose a $(k,\ell)$-decomposition $(F_k^C,F_\ell^C)$ of $H_C$ with $\{i \in [3] \mid a_i \in A(F_k^C)\}=Z$. Observe that such a decomposition exists by the definition of $H_C$. We now define a decomposition $(F_k,F_\ell)$ of $D$ by $A(F_k)=\bigcup_{x \in X \cup \mathcal{C}}A(F_k^x)$ and $A(F_\ell)=A(D)-A(F_k)$. Observe that every connected component of $F_k$ or $F_\ell$ is fully contained in $H_x$ for some $x \in X \cup \mathcal{C}$. It follows that $(F_k,F_\ell)$ is a $(k,\ell)$-decomposition of $D$, so $D$ is a positive instance of $(k,\ell)$-BDLFD.

    Now suppose that $D$ is a positive instance of $(k,\ell)$-BDLFD, so there exists a $(k,\ell)$-decomposition $(F_k,F_\ell)$ of $D$. We now define a truth assignment $\phi \colon X \rightarrow \{\true,\false\}$. Consider some $x \in X$. By the definition of $H_x$, we have that one of $\{a_1,\ldots,a_q\}\subseteq A(F_k)$ and $\{a_1,\ldots,a_q\}\subseteq A(F_\ell)$ holds. In the former case, we set $\phi(x)=\true$ and in the latter case, we set $\phi(x)=\false$. In order to prove that $\phi$ is a satisfying assignment, we consider some $C \in \mathcal{C}$. We let $x_1,x_2$, and $x_3$ be the variables contained in $C$ and we let $Z=\{i \in [3] \mid \phi(x_i)=\true\}$.  By the definition of $\phi$ and $H_{x_i}$ for $i \in [3]$, we obtain that $\{i \in [3] \mid a_i \in A(F_\ell^C)\}=Z$. It follows by the definition of $H_C$ that $Z$ is a non-empty, strict subset of $[3]$. This yields  that $\phi$ is a satisfying assignment for $(X,\mathcal{C})$, so $(X,\mathcal{C})$ is a positive instance of ME-$1$-SAT.
\end{proof}

\subsection{Decomposing into directed linear forests, one of which is unbounded}
 \label{sec:unbounded_case}
\label{infneg}
Finally, in this section, we show the result for the case that one of the two directed linear forests is unbounded. More precisely, we prove the following result.
\begin{theorem}
    For every integer $k \geq 1$, $(\infty,k)$-BDLFD  is NP-complete.
\end{theorem}
\begin{proof}
    The problem clearly being in NP, we prove the hardness by a reduction from the problem of deciding whether a given 2-diregular digraph admits a hamiltonian cycle. Recall that this problem is NP-complete by Theorem~\ref{thm:hardness_hamiltonian}. 

    Let $D$ be a 2-diregular digraph and let $x$ be an arbitrary vertex in $V(D)$. We now create a digraph $D'$ in the following way: For every $v \in V(D)$, we let $V(D')$ contain two vertices $v_+$ and $v_-$. Further, if $k \geq 2$, for every $v \in V(D)-x$, we let $V(D')$ contain $2k-2$ additional vertices $v_2,\ldots,v_k,v'_2,\ldots,v'_k$. Next, for every arc $uv \in A(D)$, we let $A(D')$ contain the arc $u_+v_-$. If $k=1$, for every $v \in V(D)-x$, we let $A(D')$ contain the arc $v_-v_+$. If $k \geq 2$, for every $v \in V(D)-x$, we let $A(D')$ contain the arcs $v_-v_2,v_-v_2',v_2'v_2,v_kv_+$, and the arcs $v_{i-1}v_i,v_{i-1}v_i'$, and $v_i'v_i$ for every $i \in \{3,\ldots,k\}$. This finishes the description if $D'$. For an illustration, see Figure~\ref{fig:unbounded_reduction}. 

    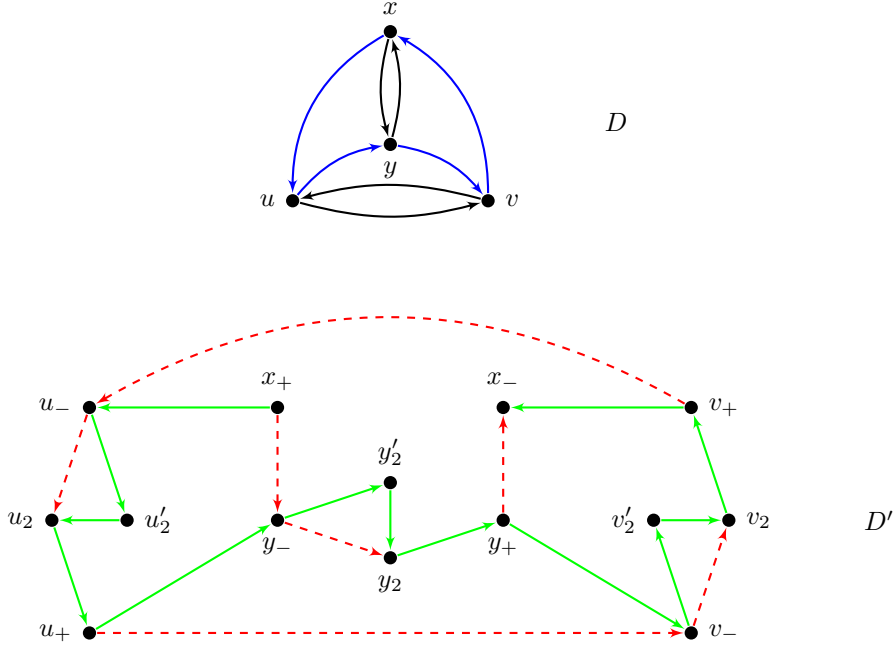
\begin{figure}[hbt!]\begin{center}
  \begin{tikzpicture}[thick,scale=1, every node/.style={transform shape}]
        \tikzset{vertex/.style = {circle,fill=black,minimum size=5pt, inner sep=0pt}}
        \tikzset{edge/.style = {->,> = latex'}}

        \node[vertex,label=above:$x$] (x) at (90:1.5){};
        \node[vertex,label=below:$y$] (y) at (0,0){};
        \node[vertex,label=left:$u$] (u) at (210:1.5){};
        \node[vertex,label=right:$v$] (v) at (-30:1.5){};
        \draw[edge,blue] (x) to[out=-150, in=90] (u){};
        \draw[edge,blue] (u) to[out=50, in=-170] (y){};
        \draw[edge,blue] (y) to[out=-10, in=130] (v){};
        \draw[edge,blue] (v) to[out=90, in=-30] (x){};
        \draw[edge, bend right=15] (x) to (y){};
        \draw[edge, bend right=15] (y) to (x){};
        \draw[edge, bend right=15] (v) to (u){};
        \draw[edge, bend right=15] (u) to (v){};
        \node[] (D) at (3,0.3) {$D$};

        \begin{scope}[yshift=-5cm]
        \node[vertex,label=above:$x_-$] (xm) at (1.5,1.5){};
        \node[vertex,label=above:$x_+$] (xp) at (-1.5,1.5){};
        
        \node[vertex,label=below:$y_-$] (ym) at (-1.5,0){};
        \node[vertex,label=below:$y_+$] (yp) at (1.5,0) {};
        \node[vertex,label=below:$y_2$] (y2) at (0,-0.5){};
        \node[vertex,label=above:$y_2'$] (y22) at (0,0.5){};
        
        \node[vertex,label=left:$u_-$] (um) at (-4,1.5){};
        \node[vertex,label=left:$u_+$] (up) at (-4,-1.5) {};
        \node[vertex,label=left:$u_2$] (u2) at (-4.5,0){};
        \node[vertex,label=right:$u_2'$] (u22) at (-3.5,0){};
            
        \node[vertex,label=right:$v_-$] (vm) at (4,-1.5){};
        \node[vertex,label=right:$v_+$] (vp) at (4,1.5) {};
        \node[vertex,label=right:$v_2$] (v2) at (4.5,0){};
        \node[vertex,label=left:$v_2'$] (v22) at (3.5,0){};

        \draw[edge,green] (xp) to (um){};
        \draw[edge,green] (up) to (ym){};
        \draw[edge,green] (yp) to (vm){};
        \draw[edge,green] (vp) to (xm){};
        \draw[edge,red, dashed] (xp) to (ym){};
        \draw[edge,red, dashed] (yp) to (xm){};
        \draw[edge,red, dashed] (vp) to[out=150, in=30] (um){};
        \draw[edge,red, dashed] (up) to (vm){};
        
        \draw[edge,green] (ym) to (y22){};
        \draw[edge, red, dashed] (ym) to (y2){};
        \draw[edge,green] (y22) to (y2){};
        \draw[edge,green] (y2) to (yp){};
        
        \draw[edge,green] (um) to (u22){};
        \draw[edge, red, dashed] (um) to (u2){};
        \draw[edge,green] (u22) to (u2){};
        \draw[edge,green] (u2) to (up){};
        
        \draw[edge,green] (vm) to (v22){};
        \draw[edge, red, dashed] (vm) to (v2){};
        \draw[edge,green] (v22) to (v2){};
        \draw[edge,green] (v2) to (vp){};
        
        \node[] (Dp) at (6.5,0) {$D'$};
        \node[] (Dp2) at (-6.5,0) {};
        \end{scope}
    
    \end{tikzpicture}
    \caption{An example for the construction of the digraph $D'$ from the digraph $D$ with $k=2$, together with an $(\infty,k)$-decomposition $(F_\infty,F_k)$ of $D'$, built from the hamiltonian cycle of $D$ coloured in blue. The dashed red arcs are in $A(F_k)$ and the solid green arcs are in $A(F_\infty)$.}\label{fig:unbounded_reduction}
\end{center}
\end{figure}
    
    We now prove that $D'$ admits an $(\infty,k)$-decomposition if and only if $D$ contains a hamiltonian cycle.
    First suppose that $D$ contains a hamiltonian cycle $C$. We now define a set $A_\infty \subseteq A(D')$. We first let $A_\infty$ contain the arc $u_+v_-$ for every arc $uv \in A(C)$. Next, if $k=1$,  for every $v \in V(D)-x$, we let $A_\infty$ contain the arc $v_-v_+$ and if $k\geq 2$, we let $A_\infty$ contain the arcs $v_-v_2',v_2'v_2,v_kv_+$, and the arcs $v_{i}'v_i$ and $v_{i-1}v_i'$ for every $i \in \{3,\ldots,k\}$. We now define the decomposition $(F_\infty,F_k)$ of $D'$ by $A(F_\infty)=A_\infty$ and $A(F_k)=A(D')-A_\infty$. In order to see that $F_k$ is a $k$-bounded directed linear forest observe that by construction, every connected component of $F_k$ is a directed path of the form $u_+v_-$ if $k=1$ and a directed path of the form $u_+v_-v_2\ldots v_k$ if $k \geq 2$ for some arc $uv \in A(D)$. In order to see that $F_\infty$ is a directed linear forest, first observe that we have $d_{F_\infty}^+(x_+)=1,d_{F_\infty}^-(x_+)=0,d_{F_\infty}^+(x_-)=0,d_{F_\infty}^-(x_-)=1$ and $d_{F_\infty}^+(v)=d_{F_\infty}^-(v)=1$ for all $v \in V(D')-\{x_+,x_-\}$. It hence suffices to prove that $F_\infty$ does not contain a directed cycle $C'$ with $V(C')\subseteq V(D')-\{x_+,x_-\}$. Suppose for the sake of a contradiction that $C'$ is such a cycle and let $C''$ be the subdigraph of $D$ defined by $A(C'')=\{uv \in A(D) \mid u_+v_- \in A(C')\}$. As $C'$ is a directed cycle, we have $d_{C''}^+(v)\leq 1$ and $d_{C''}^+(v)=d_{C''}^-(v)$ for all $v \in V(D)$. It follows that $C''$ is a disjoint collection of directed cycles and disjoint cycles. By construction, we have that $A(C'')$ is nonempty and $d_{C''}^+(x)=0$. Finally, by construction, we have that $C''$ is a subdigraph of $C$. This contradicts $C$ being a hamiltonian cycle of $D$.

   Now suppose that there exists an $(\infty,k)$-decomposition $(F_\infty,F_k)$ of $D$. Let $C$ be the spanning subdigraph of $D$ defined by $A(C)=\{uv \in A(D) \mid u_+v_- \in A(F_\infty)\}$. As $F_\infty$ and $F_k$ are directed linear forests, we have $d_C^+(v)=d_C^-(v)=1$ for all $v \in V(D)$, so $C$ is a collection of directed cycles. Suppose for the sake of a contradiction that $C$ contains a connected component $C'$ which is a directed cycle $v^1\ldots v^q$ with $x \notin \{v^1,\ldots,v^q\}$. If for every $i \in [q]$, we have that $F_\infty$ contains a directed $v^{i}_-v^{i}_+$-path, then $F_\infty$ contains a connected component containing the arcs of all those paths and $\{u_+v_- \mid uv \in A(C')\}$. This contradicts $F_\infty$ being a directed linear forest. Hence there exists some $w \in \{v^1,\ldots,v^q\}$ such that $F_\infty$ does not contain a directed $w_-w_+$-path. Observe that, as both $F_\infty$ and $F_k$ are directed linear forests, we have that both $F_\infty$ and $F_k$ contain a directed path from $w_-$ to the unique vertex in $N_D^-(w_+)$ and both these paths are of length at least $k-1$. As $F_\infty$ does not contain a directed $w_-w_+$-path, we obtain that the unique edge in $\delta_{D'}^-(w_+)$ is contained in $A(F_k)$. It follows that $w_-$ is the first vertex of a path of length at least $k$ in $F_k$. Hence, as $F_k$ is a $k$-bounded directed linear forest, both arcs entering $w_-$ in $A(D')$ are contained in $A(F_\infty)$. This contradicts $F_\infty$ being a directed linear forest. We obtain a contradiction to the existence of $C'$. It follows that $C$ is a hamiltonian cycle of $D$.
\end{proof}

\section{Decomposing into out-galaxies}
\label{sec:galaxy}

This section is dedicated to proving Theorem~\ref{bogdmain}. In Section~\ref{out1}, we give the simple proof that the problem of decomposing a given digraph into two unbounded out-galaxies can be reduced to finding a bipartition of a certain graph and can hence be solved in polynomial time. In Section~\ref{sec:galaxy_1_k}, we give the somewhat more involved proof that the problem of decomposing a given digraph into a directed matching and a (bounded) out-galaxy can be solved in polynomial time by a reduction to a matching problem. In Section~\ref{sec:galaxy_2_2}, we show that the problem is NP-complete in all remaining cases. The combination of these results yields Theorem~\ref{bogdmain}. Throughout this section, for $k,\ell \in \mathbb{Z}_{\geq 1} \cup \{\infty\}$, a {\bf $(k,\ell)$-factorization} of a digraph $D$ is a decomposition $(F_k,F_\ell)$ of $D$ such that $F_k$ is a $k$-bounded out-galaxy, and $F_\ell$ is an $\ell$-bounded out-galaxy.

\subsection{Decomposing into two unbounded out-galaxies}
\label{out1}

We here prove the following result which covers one of the cases of Theorem~\ref{bogdmain}. The proof consists of a simple reduction to a bipartition problem.

\begin{theorem}
    $(\infty,\infty)$-BOGD is solvable in polynomial time.
\end{theorem}
\begin{proof}
    Let $D$ be an instance of $(\infty,\infty)$-BOGD. We create an undirected graph $G$ in the following way. For every $a \in A(D)$, we let $V(G)$ contain a vertex $v_a$ and we let $E(G)$ contain an edge linking two vertices $v_a$ and $v_{a'}$ if $a$ and $a'$ share a vertex which is not the tail of both of them in $D$.
    \begin{claim}\label{bip}
        $D$ is a positive instance of $(\infty,\infty)$-BOGD if and only if $G$ is bipartite.
    \end{claim}
    \begin{proofclaim}
        First suppose that $D$ is a positive instance of $(\infty,\infty)$-BOGD, so there exists an $(\infty,\infty)$-factorization $(F,F')$ of $D$. Let $X=\{v_a \mid a \in A(F)\}$ and $Y=V(G)-X$. Let $v_a,v_{a'}\in X$. Then, as $F$ is an out-galaxy, we obtain that $a$ and $a'$ are either disjoint or have a common tail. We obtain that $v_a$ and $v_{a'}$ are not linked by an edge in $G$. It follows that $X$ is an independent set. Similarly, $Y$ is an independent set, so $(X,Y)$ is a bipartition of $G$.

        Now suppose that $G$ is bipartite and let $(X,Y)$ be a bipartition of $G$. Let $F$ be the spanning subdigraph of $D$ with $A(F)=\{a \in A \mid v_a \in X\}$ and let $F'$ be the spanning subdigraph of $D$ with $A(F')=A(D)- A(F)$. As $X$ is an independent set in $G$, we obtain that all arcs $a,a'\in A(F)$ are either disjoint or share their tail. It follows that $F$ is an out-galaxy. Similarly, $F'$ is an out-galaxy, so $(F,F')$ is an $(\infty,\infty)$-factorization of $D$.
        
    \end{proofclaim}
    
By Claim~\ref{bip}, it suffices to check whether $G$ is bipartite. By Proposition~\ref{checkbip} and as $G$ can be constructed from $D$ in polynomial time, this can be done in polynomial time.
\end{proof}

\subsection{Decomposing into a matching and an out-galaxy}
\label{sec:galaxy_1_k}

This section is dedicated to proving that a decomposition of a given digraph into a directed matching and a possibly bounded out-galaxy can be found in polynomial time. The proof is based on a reduction to a matching problem in undirected graphs and has certain similarities with the proof of Theorem~\ref{21dblfd} and the proof of~\cite[Theorem~2]{campbell2023decompositions}. We first need a collection of preliminary results which deal with certain decompositions of orientations of cycles and paths.
\begin{proposition}\label{2reg}
    Let $C$ be an orientation of a cycle. Then, for any $k \in \mathbb{Z}_{\geq 2} \cup \{\infty\}$ there exists a $(k,1)$-factorization of $C$ if and only if $k$ is even or $C$ is not oriented as a circuit.
\end{proposition}
\begin{proof}
Let $C$ be an orientation of the cycle $v_1 v_2\ldots v_rv_1$. If $r$ is even, we can define a decomposition  $(F_k,F_1)$ of $C$ in which $A(F_k)$ consists of the orientation of the edge $v_iv_{i+1}$ for all odd $i \in [r]$ and $A(F_1)=A(C)-A(F_k)$. It is easy to see that $(F_k,F_1)$ is a $(k,1)$-factorization of $C$. 

We may hence suppose that $r$ is odd. First suppose $C$ is oriented as a circuit, say $v_1 v_2\ldots v_rv_1$ and that there exists a $(k,1)$-factorization $(F_k,F_1)$ of $C$. At least one arc of $A(C)$ belongs to $F_k$. Hence, by symmetry, we may suppose that $v_rv_1 \in A(F_k)$. As $F_k$ and $F_1$ are out-galaxies, we inductively obtain that $v_iv_{i+1}\in A(F_1)$ for all odd $i \in [r-2]$ and $v_iv_{i+1}\in A(F_k)$ for all even $i \in [r-1]$. This yields $v_{r-1}v_r,v_rv_1 \in A(F_k)$, a contradiction to $F_k$ being an out-galaxy. Hence in this case there exists no $(k,1)$-factorization of $C$.
 Now suppose that $C$ is not a circuit. Then $C$ contains a source, that is some vertex $v \in V(C)$ with $d_C^+(v)=2$, say $v_1$. Then let $A(F_k)$ contain $v_1v_r$ and the orientation of the edge $v_iv_{i+1}$ for all odd $i \in [r-2]$ and let $A(F_1)=A(C)-A(F_k)$. It is easy to see that $(F_k,F_1)$ is a $(k,1)$-factorization of $C$.
 \end{proof}
 
\begin{proposition}\label{chemin}
    Let $P$ be an orientation of a path of length at least 2, $a_1,a_2$ the endarcs of $P$, and let $A_0 \subseteq \{a_1,a_2\}$. Then we can decide in polynomial time whether there exists a $(k,1)$-factorization $(F_k,F_1)$ of $P$ with $A(F_1)\cap \{a_1,a_2\}=A_0$.
\end{proposition}
\begin{proof}
Let the underlying graph of $P$ be the path $v_1\ldots v_q$. If $q = 3$, we can solve the instance by a brute force approach. We may hence suppose that $q \geq 4$. 
We distinguish three cases.
\begin{description}
    \item[Case 1:] \textit{$A_0 \neq \emptyset$.}

    We assume without loss of generality that $a_1 \in A_0$. Then there exists a $(k,1)$-factorization $(F_k,F_1)$ of $P$ with $A(F_1)\cap \{a_1,a_2\}=A_0$ if and only if there exists a $(k,1)$-factorization $(F'_k,F'_1)$ of $P-v_1$ with $A(F'_1)\cap \{a'_1,a_2\}=A_0-a_1$ where $a'_1$ is the orientation of the edge $v_2v_3$ in $P$. We can hence recursively solve this smaller instance.

    \item[Case 2:] \textit{$A_0 = \emptyset$ and $A(P)$ does not contain $\{v_2v_1,v_2v_3\}$.}

    In this case, there exists a $(k,1)$-factorization $(F_k,F_1)$ of $P$ with $A(F_1)\cap \{a_1,a_2\}=\emptyset$ if and only if there exists a $(k,1)$-factorization $(F'_k,F'_1)$ of $P-v_1$ with $A(F'_1)\cap \{a'_1,a_2\}=\{a'_1\}$ where $a'_1$ is the orientation of the edge $v_2v_3$ in $P$. We can hence recursively solve this smaller instance.

    \item[Case 3:] \textit{$A_0 = \emptyset$ and $A(P)$ contains $\{v_2v_1,v_2v_3\}$.}

    We will show that in this case, we can find the desired decomposition. If $q$ is odd, we define $F_k$ by letting $A(F_k)$ consist of $v_2v_1$ and the orientation of $v_iv_{i+1}$ for all even integers $i \in \{2,\ldots,q-1\}$. If $q$ is even, we define $F_k$ by letting $A(F_k)$ consist of the orientation of $v_iv_{i+1}$ for all odd integers $i \in \{1,\ldots,q-1\}$. We define $F_1$ by $A(F_1)=A(P) -  A(F_k)$. In both cases, it is easy to see that $(F_k,F_1)$ is a $(k,1)$-factorization of $P$ and that $A(F_1) \cap \{a_1,a_2\} = \emptyset$. \qedhere
\end{description}
\end{proof}

\begin{figure}[hbt!]\begin{center}
    \begin{tikzpicture}[thick,scale=1, every node/.style={transform shape}]
            \tikzset{vertex/.style = {circle,fill=black,minimum size=5pt, inner sep=0pt}}
            \tikzset{squarevertex/.style = {rectangle,fill=black,minimum size=5pt, inner sep=0pt}}
            \tikzset{edge/.style = {->,> = latex'}}

            \begin{scope}
                \node[vertex, label=left:$v_1$] (v1) at (0, 0) {};
                \node[vertex, label=right:$v_2$] (v2) at (3, 0) {};
                \node[squarevertex, label=left:$z_1$] (z1) at (0, 1.5) {};
                \node[squarevertex, label=right:$z_2$] (z2) at (3, 1.5) {};
                \node[] (e1) at (-0.3,0.75) {$e_1$};
                \node[] (e2) at (3.3,0.75) {$e_2$};
                \node[] (a) at (1.5,-0.3) {$(a)$};
                \draw (v1) -- (z1) -- (z2) -- (v2);
            \end{scope}

            \begin{scope}[xshift=4.5cm]
                \node[vertex, label=left:$v_1$] (v1) at (0, 0) {};
                \node[vertex, label=right:$v_2$] (v2) at (3, 0) {};
                \node[squarevertex, label=above:$z$] (z) at (1.5, 1.5) {};
                \node[] (e1) at (0.4,0.75) {$e_1$};
                \node[] (e2) at (2.6,0.75) {$e_2$};
                \node[] (a) at (1.5,-0.3) {$(b)$};
                \draw (v1) -- (z) -- (v2);
            \end{scope}

            \begin{scope}[xshift=9cm]
                \node[vertex, label=left:$v_1$] (v1) at (0, 0) {};
                \node[vertex, label=right:$v_2$] (v2) at (3, 0) {};
                \node[vertex, label=above:$w$] (z) at (1.5, 1.5) {};
                \node[] (e1) at (0.4,0.75) {$e_1$};
                \node[] (e2) at (2.6,0.75) {$e_2$};
                \node[] (a) at (1.5,-0.3) {$(c)$};
                \draw (v1) -- (z) -- (v2);
            \end{scope}

            \begin{scope}[yshift=-3cm]
                \node[vertex, label=left:$v_1$] (v1) at (0, 0) {};
                \node[vertex, label=right:$v_2$] (v2) at (3, 0) {};
                \node[squarevertex, label=left:$z$] (z1) at (0, 1.5) {};
                \node[vertex, label=right:$w$] (z2) at (3, 1.5) {};
                \node[] (e1) at (-0.3,0.75) {$e_1$};
                \node[] (e2) at (3.3,0.75) {$e_2$};
                \node[] (a) at (1.5,-0.3) {$(d)$};
                \draw (v1) -- (z1) -- (z2) -- (v2);
            \end{scope}

            \begin{scope}[yshift=-3cm, xshift=4.5cm]
                \node[vertex, label=left:$v_1$] (v1) at (0, 0) {};
                \node[vertex, label=right:$v_2$] (v2) at (3, 0) {};
                \node[squarevertex, label=left:$z_1$] (z1) at (0, 1.5) {};
                \node[vertex, label=above:$w$] (w) at (1.5, 1.5) {};
                \node[squarevertex, label=right:$z_2$] (z2) at (3, 1.5) {};
                \node[] (e1) at (-0.3,0.75) {$e_1$};
                \node[] (e2) at (3.3,0.75) {$e_2$};
                \node[] (a) at (1.5,-0.3) {$(e)$};
                \draw (v1) -- (z1) -- (w) -- (z2) -- (v2);
            \end{scope}

            \begin{scope}[yshift=-3cm, xshift=9.75cm]
                \node[vertex, label=left:$v_1$] (v1) at (0, 0) {};
                \node[vertex, label=right:$v_2$] (v2) at (1.5, 0) {};
                \node[vertex, label=left:$w_1$] (w1) at (0, 1.5) {};
                \node[vertex, label=right:$w_2$] (w2) at (1.5, 1.5) {};
                \node[] (e1) at (-0.3,0.75) {$e_1$};
                \node[] (e2) at (1.8,0.75) {$e_2$};
                \node[] (a) at (0.75,-0.3) {$(f)$};
                \draw (v1) -- (w1);
                \draw (v2) -- (w2);
            \end{scope}
        \end{tikzpicture}
  \caption{An illustration of the $X$-gadgets for all sets $X$ which are relevant by Proposition~\ref{chemin2}. The vertices in $Z$ are marked by squares and the remaining vertices by disks.}\label{gadgmatch}
  \end{center}
\end{figure}
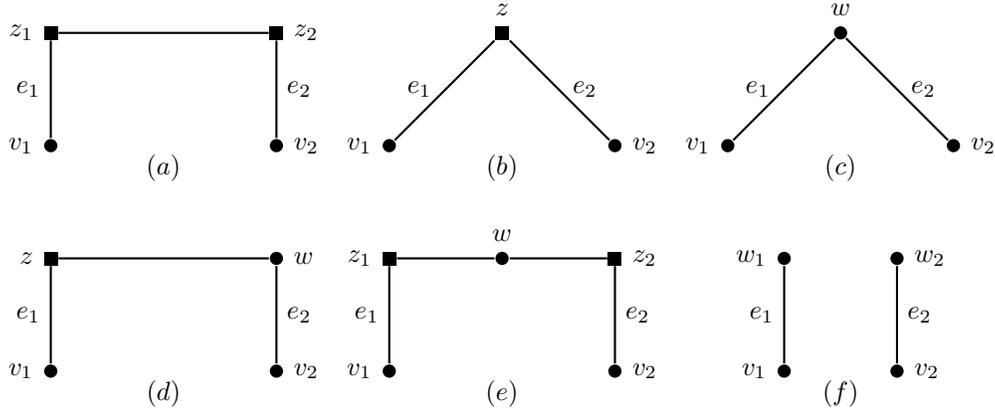

\begin{proposition}\label{chemin2}
    Let $P$ be an orientation of a path of length at least 2, $a_1,a_2$ the endarcs of $P$ and let $X \subseteq 2^{\{1,2\}}$ be the set such that for every $A_0 \subseteq\{a_1,a_2\}$, there exists a $(k,1)$-factorization $(F_k,F_1)$ of $P$ with $A(F_1)\cap \{a_1,a_2\}=A_0$ if and only if $\{i \in \{1,2\} \mid a_i \in A_0\}\in X$. Then at least one of $\{\emptyset,\{1,2\}\}\subseteq X$ and $\{\{1\},\{2\}\}\subseteq X$ holds. Moreover, we can compute $X$ in polynomial time.
\end{proposition}
\begin{proof}
    Let $P$ be an orientation of the path $v_1 \ldots v_q$ for some positive integer $q\geq 3$. We can define a decomposition $(F_k,F_1)$ of $P$ where $A(F_k)$ contains the orientation of $v_iv_{i+1}$ for all odd integers $i$ and $A(F_1)$ contains the orientation of $v_iv_{i+1}$ for all even integers $i$. Further, we can define $(F'_k,F'_1)$ by $A(F'_k)=A(F_1)$ and $A(F'_1)=A(F_k)$. It is easy to see that both $(F_k,F_1)$ and $(F'_k,F'_1)$ are $(k,1)$-factorizations of $P$. This proves that $\{\emptyset,\{1,2\}\}\subseteq X$ if $q$ is even and $\{\{1\},\{2\}\}\subseteq X$ if $q$ is odd. In order to compute $X$, we test whether $I \in X$ for every $I \subseteq \{1,2\}$. As there are only four sets to test, this can be done in polynomial time by Proposition~\ref{chemin}.
\end{proof}
In order to reduce our decomposition problem to a matching problem, we need a collection of gadgets. Each of these gadgets will correspond to a path in the underlying graph of the input digraph and the gadget will reflect the possible decompositions of the corresponding subdigraph of the input digraph. The following result shows that all the desired gadgets exist.
\begin{lemma}\label{couplagegadget}
    Let $X \subseteq 2^{\{1,2\}}$ with $\{\emptyset,\{1,2\}\}\subseteq X$ or $\{\{1\},\{2\}\}\subseteq X$. Then there exists a graph $G$ together with some $v_1,v_2 \in V(G), e_1,e_2 \in E(G)$ and $Z \subseteq V(G)-\{v_1,v_2\}$ with the following properties:
    \begin{itemize}
    \item for $i=1,2$, the only edge incident to $v_i$ is $e_i$,
    \item for every $I \in X$, there exists a matching $M$ in $G$ covering $Z$ with $M \cap \{e_1,e_2\}=\{e_i \mid i \in I\}$,
    \item for every matching $M$ in $G$ covering $Z$, we have $\{i \in \{1,2\} \mid e_i \in M\}\in X$.
    \end{itemize}
\end{lemma}
\begin{proof}
    We explicitly give the gadgets for every possible set $X$, omitting sets which are fully symmetric to previously considered ones. As it is easily visible that these gadgets have the desired properties, we do not prove that in detail.

    \begin{description}
        \item[Case~1:] $X=\{\emptyset,\{1,2\}\}$.

        We set $V(G)=\{v_1,v_2,z_1,z_2\}$, $E(G)=\{e_1=v_1z_1,e_2=v_2z_2,z_1z_2\}$ and $Z=\{z_1,z_2\}$, see Figure~\ref{gadgmatch}~$(a)$.  

        \item[Case~2:] $X=\{\{1\},\{2\}\}$.

        We set $V(G)=\{v_1,v_2,z\}$, $E(G)=\{e_1=v_1z,e_2=v_2z\}$, and $Z=\{z\}$, see Figure~\ref{gadgmatch}~$(b)$. 

        \item[Case~3:] $X=\{\emptyset,\{1\},\{2\}\}$.
        We set $V(G)=\{v_1,v_2,w\}$, $E(G)=\{e_1=v_1w,e_2=v_2w\}$ and $Z=\emptyset$, see Figure~\ref{gadgmatch}~$(c)$.

        \item[Case~4:] $X=\{\emptyset,\{1\},\{1,2\}\}$.

        We set $V(G)=\{v_1,v_2,z,w\}$, $E(G)=\{e_1=v_1z,e_2=v_2w,wz\}$ and $Z=\{z\}$, see Figure~\ref{gadgmatch}~$(d)$.

        \item[Case~5:] $X=\{\emptyset,\{2\},\{1,2\}\}$.

        This case is symmetric to the previous one.
        
        \item[Case~6:] $X=\{\{1\},\{2\},\{1,2\}\}$.

        We set $V(G)=\{v_1,v_2,z_1,z_2,w\}$, $E(G)=\{e_1=v_1z_1,e_2=v_2z_2,z_1w,z_2w\}$ and $Z=\{z_1,z_2\}$, see Figure~\ref{gadgmatch} $(e)$.

        \item[Case~7:] $X=\{\emptyset,\{1\},\{2\},\{1,2\}\}$.
    We set $V(G)=\{v_1,v_2,w_1,w_2\}$, $E(G)=\{e_1=v_1w_1,e_2=v_2w_2\}$ and $Z=\emptyset$, see Figure~\ref{gadgmatch}~$(f)$. \qedhere
    \end{description}
\end{proof}

We say that $(G,v_1,v_2,e_1,e_2,Z)$ as described in Lemma~\ref{couplagegadget} is an $X$-gadget.
We are now ready to prove the main result of this section. We first exclude some local configurations which clearly make the instance negative and then reduce the problem to a matching problem.
\begin{theorem}
    For every $k \in \mathbb{Z}_{\geq 1}\cup \{\infty\}$, $(k,1)$-BOGD is solvable in polynomial time, even if $k$ is part of the input.
\end{theorem}
\begin{proof}Let $D$ be an instance of $(k,1)$-BOGD for some $k \in \mathbb{Z}_{\geq 1}\cup \{\infty\}$. 
    Observe that $D$ is a positive instance of $(\infty,1)$-BOGD if and only if $D$ is a positive instance of $(k,1)$-BOGD for $k=\max_{v \in V}d_{D}(v)$. We may hence suppose that $k$ is finite.

    A vertex $v \in V$ is called {\bf big} if it satisfies $d_{D}(v)\geq 3$ and {\bf tiny} if it satisfies $d_{D}(v)=1$. We denote by $B$ and $T$ the set of big and tiny vertices, respectively. By Proposition~\ref{2reg}, we may suppose that every connected component of $D$ that contains at least one arc contains a vertex in $B \cup T$. Further, for $i=0,1$, we denote by $B_i$ the set of vertices in $B$ with $d_D^-(v)=i$ and $d_{D}(v)\leq k+1$. Clearly, if $D$ contains a vertex in $B-(B_0 \cup B_1)$, then $D$ is a negative instance  of $(k,1)$-BOGD. As this property can be checked in polynomial time, we may suppose that $B=B_0 \cup B_1$. We further let $\mathcal{P}$ be the set of subdigraphs of $D$ whose underlying graphs forms either a path in $\UG(D)$ connecting two vertices of $B \cup T$ and none of whose interior vertices is contained in $B \cup T$ or a cycle containing exactly one vertex in $B \cup T$. We say that $P$ is {\bf incident} to these vertices in $B \cup T$. Observe that $\{A(P) \mid P \in \mathcal{P} \}$ is a partition of $A(D)$.

    We now create an undirected graph $G$. First, we let $V(G)$ contain $B_0 \cup T$. For every $b \in B_0 \cup T$ and every arc $a \in A$ which is incident to $b$, we also refer to $b$ as $u_b^a$. Next, for every $b \in B_1$ and every $a \in \delta_D^+(b)$, we let $V(G)$ contain two vertices $u_b^a$ and $z_b^a$, and we let $E(G)$ contain an edge $u_b^az_b^a$. Further, for every $b \in B_1$ and every $a \in \delta_D^-(b)$, we let $V(G)$ contain a vertex $u_b^a$.

    Now consider some $P \in \mathcal{P}$. Let $b_1,b_2 \in B \cup T$ be the  vertices incident to $P$ with $b_1=b_2$ if $P$ is a cycle. We define the {\bf associated path $P'$} of $P$ by $P'=P$ if $P$ is a path and the path obtained by detaching $b_1$ into two vertices $b_1,b_2$ with $d_{P'}(b_1)=d_{P'}(b_2)=1$ if $P$ is a cycle. Further, let $a_1,a_2$ be the arcs in $A(P')$ such that $a_i$ is incident to $b_i$. We now compute the set $X_P\subseteq 2^{\{1,2\}}$ such that for every $I \subseteq \{1,2\}$, we have that $P'$ admits a $(k,1)$-decomposition $(F_k,F_1)$ with $\{i \in \{1,2\} \mid a_i \in A(F_1)\}=I$ if and only if $I \in X^P$. Observe that $X^P=\{\emptyset,\{1,2\}\}$ if $P'$ is of length 1. By Proposition~\ref{chemin}, we obtain that $X^P$ can be computed in polynomial time.
    
    We now add an $X^P$-gadget $(G^P,v_1^P,v_2^P,e_1^P,e_2^P,Z^P)$ and identify $v_1^P$ with $u_{b_1}^{a_1}$ and $v_2^P$ with $u_{b_2}^{a_2}$. The existence of this gadget is guaranteed by Proposition~\ref{chemin2} and Lemma~\ref{couplagegadget}. We do this for every $P \in \mathcal{P}$. This finishes the description of $G$. We now let $Z$ consist of $\bigcup_{P \in \mathcal{P}}Z^P$, the vertex $u_b^a$ for every $b \in B_1$ and $a \in \delta^-(b)$, the vertex $z_b^a$ for every $b \in B_1$ and $a \in \delta^+(b)$  and the vertex $b$ for every $b \in B_0$ with $d_D^+(b)=k+1$. For an illustration, see Figure~\ref{drzftgz}.
    
    \begin{figure}[hbt!]
     \begin{center}
     \begin{tikzpicture}[thick,scale=1, every node/.style={transform shape}]
            \tikzset{vertex/.style = {circle,fill=black,minimum size=5pt, inner sep=0pt}}
            \tikzset{smallvertex/.style = {circle,fill=black,minimum size=3pt, inner sep=0pt}}
            \tikzset{squarevertex/.style = {rectangle,fill=black,minimum size=5pt, inner sep=0pt}}
            \tikzset{edge/.style = {->,> = latex'}}

            \node[vertex, label=above:$v_1$] (v1) at (0, 0) {};
            \node[vertex, label=above:$v_2$] (v2) at (2, 0) {};
            \node[vertex, label=above:$v_3$] (v3) at (6, 0) {};
            \node[vertex, label=below:$v_4$] (v4) at (2, -2) {};
            \node[vertex, label=below:$v_5$] (v5) at (6, -2) {};

            \node[smallvertex] (P1a) at (0.666, 0) {};
            \node[smallvertex] (P1b) at (1.333, 0) {};
            \draw[edge, orange] (P1a) -- (v1) {};
            \draw[edge, orange] (P1a) -- (P1b);
            \draw[edge, orange] (P1b) -- (v2) {};
            \node[orange] (P1) at (1,0.3) {$P_1$};

            \node[smallvertex] (P2a) at (4, 0.8) {};
            \draw[edge, green] (v2) to[out=30, in=180] (P2a) {};
            \draw[edge, green] (v3) to[out=150, in=0] (P2a);
            \node[green] (P2) at (4,1.1) {$P_2$};

            \foreach \i in {1,...,6}{
                \node[smallvertex] (P3\i) at (2+0.571428*\i, 0) {};
            }
            \draw[edge, teal] (v2) to (P31) {};
            \draw[edge, teal] (P32) to (P31) {};
            \draw[edge, teal] (P32) to (P33) {};
            \draw[edge, teal] (P33) to (P34) {};
            \draw[edge, teal] (P35) to (P34) {};
            \draw[edge, teal] (P35) to (P36) {};
            \draw[edge, teal] (v3) to (P36) {};
            \node[teal] (P3) at (4,0.3) {$P_3$};

            \foreach \i in {1,...,5}{
                \node[smallvertex] (P4\i) at (2+0.66666*\i, -0.333333*\i) {};
            }
            \draw[edge, violet] (v2) to (P41) {};
            \draw[edge, violet] (P41) to (P42) {};
            \draw[edge, violet] (P43) to (P42) {};
            \draw[edge, violet] (P43) to (P44) {};
            \draw[edge, violet] (P45) to (P44) {};
            \draw[edge, violet] (v5) to (P45) {};
            \node[violet] (P4) at (4,-0.6) {$P_4$};

            \foreach \i in {1,...,2}{
                \node[smallvertex] (P5\i) at (2.2+1.3333333*\i, -2.2 + 0.666666*\i) {};
            }
            \draw[edge, pink] (v4) to (P51) {};
            \draw[edge, pink] (P51) to (P52) {};
            \draw[edge, pink] (v3) to (P52) {};
            \node[pink] (P5) at (4,-1.7) {$P_5$};

            \draw[edge, g-brown] (v5) to (v3) {};
            \node[g-brown] (P6) at (6.3,-1) {$P_6$};

            \node[smallvertex] (P7a) at (4, -2.2) {};
            \draw[edge, blue] (v4) to[out=-10, in=180] (P7a) {};
            \draw[edge, blue] (v5) to[out=-170, in=0] (P7a);
            \node[blue] (P7) at (4,-2.5) {$P_7$};

            \node[smallvertex] (P8a) at (2, -1) {};
            \node[smallvertex] (P8b) at (1, -1) {};
            \node[smallvertex] (P8c) at (1, -2) {};
            \draw[edge, red] (v4) to (P8a) {};
            \draw[edge, red] (v4) to (P8c) {};
            \draw[edge, red] (P8b) to (P8a) {};
            \draw[edge, red] (P8b) to (P8c) {};
            \node[red] (P8) at (0.5,-1.5) {$P_8$};

        \node[] (G) at (-4.5, -1) {$D$};
        \node[] (D) at (-4.5, -5) {$G$};
        \begin{scope}[xshift=3cm, yshift=-3.5cm]
            \node[squarevertex, label=right:$z_{v_2}^{a_{22}}$] (z22) at (4.5,0) {};
            \node[squarevertex, label=right:$z_{v_2}^{a_{23}}$] (z23) at (3,0) {};
            \node[squarevertex, label=right:$z_{v_2}^{a_{24}}$] (z24) at (1.5,0) {};
            
            \node[vertex, label=right:$u_{v_2}^{a_{22}}$] (u22) at (4.5,-1) {};
            \node[vertex, label=right:$u_{v_2}^{a_{23}}$] (u23) at (3,-1) {};
            \node[vertex, label=right:$u_{v_2}^{a_{24}}$] (u24) at (1.5,-1) {};
            \node[squarevertex, label=right:$u_{v_2}^{a_{21}}$] (u21) at (5.5,-1) {};
            \node[vertex, label=above:$v_1$] (v1) at (6.5,-1) {};
            
            \node[squarevertex, label=right:$z_{v_3}^{a_{32}}$] (z32) at (0,0) {};
            \node[squarevertex, label=right:$z_{v_3}^{a_{33}}$] (z33) at (-1.5,0) {};
            \node[squarevertex, label=right:$z_{v_3}^{a_{35}}$] (z35) at (-3,0) {};
            
            \node[vertex, label=right:$u_{v_3}^{a_{32}}$] (u32) at (0,-1) {};
            \node[vertex, label=right:$u_{v_3}^{a_{33}}$] (u33) at (-1.5,-1) {};
            \node[vertex, label=right:$u_{v_3}^{a_{35}}$] (u35) at (-3,-1) {};
            \node[squarevertex, label=above:$u_{v_3}^{a_{36}}$] (u36) at (-4,-1) {};
            \node[squarevertex, label=above:$v_4$] (v4) at (-5,-1) {};
            \node[vertex, label=above:$v_5$] (v5) at (-6,-1) {};

            \node[squarevertex, red] (x1) at (-6,-3) {};
            \node[vertex, red] (x2) at (-5.25,-3) {};
            \node[squarevertex, red] (x3) at (-4.5,-3) {};
            \node[squarevertex, blue] (x4) at (-3.5,-3) {};
            \node[squarevertex, g-brown] (x5) at (-2.5,-3) {};
            \node[squarevertex, g-brown] (x6) at (-1.5,-3) {};
            \node[squarevertex, pink] (x7) at (-0.5,-3) {};
            \node[squarevertex, pink] (x8) at (0.5,-3) {};
            \node[vertex, violet] (x9) at (1.5,-3) {};
            \node[vertex, teal] (x10) at (2.5,-3) {};
            \node[vertex, teal] (x11) at (3.5,-3) {};
            \node[squarevertex, green] (x12) at (4.5,-3) {};
            \node[squarevertex, orange] (x13) at (5.5,-3) {};
            \node[vertex, orange] (x14) at (6.5,-3) {};
        \end{scope}

        \draw[] (z35) -- (u35){};
        \draw[] (z33) -- (u33){};
        \draw[] (z32) -- (u32){};
        \draw[] (z22) -- (u22){};
        \draw[] (z23) -- (u23){};
        \draw[] (z24) -- (u24){};

        \draw[red] (v4) -- (x1) -- (x2) -- (x3) -- (v4){};
        \draw[blue] (v5) -- (x4) -- (v4){};
        \draw[g-brown] (v5) -- (x5) -- (x6) -- (u36){};
        \draw[pink] (v4) -- (x7) -- (x8) -- (u35){};
        \draw[violet] (v5) -- (x9) -- (u24){};
        \draw[teal] (x10) -- (u33){};
        \draw[teal] (x11) -- (u23){};
        \draw[green] (u32) -- (x12) -- (u22){};
        \draw[orange] (u21) -- (x13) -- (x14) -- (v1){};
        \end{tikzpicture}
  \caption{An example for the construction of the graph $G$ from the digraph $D$ for $k=3$. In $D$, we have $T=\{v_1\}, B_0=\{v_4,v_5\}$, and $B_1=\{v_2,v_3\}$.
  For integers $i \in [5]$ and $j\in [7]$ such that $v_i$ is an endvertex of $P_j$, the unique arc of $P_j$ incident to $v_i$ is denoted by $a_{ij}$.
  Vertices in $Z$ are marked by squares.
  The different elements of $\mathcal{P}=\{P_1,\ldots,P_8\}$ are marked in different colours. For $i\in [8]$, the gadget in $G$ corresponding to $P_i$ in $D$ is marked in the same colour as $P_i$.
  We have $X^{P_1}=\{\emptyset,\{2\},\{1,2\}\}$, $X^{P_2}=\{\{1\},\{2\}\}$, $X^{P_3}=\{\emptyset,\{1\},\{2\},\{1,2\}\}$, $X^{P_4}=\{\emptyset,\{1\},\{2\}\}$, $X^{P_5}=\{\emptyset,\{1,2\}\}$, $X^{P_6}=\{\emptyset,\{1,2\}\}$, $X^{P_7}=\{\{1\},\{2\}\}$, and $X^{P_8}=\{\{1\},\{2\}\,\{1,2\}\}$, where for $i \in [7]$, we let $(b_1,b_2)$ be the endvertices of $P_i$ in increasing order with respect to $(v_1,\ldots,v_5)$.}
  \label{drzftgz}
\end{center}
\end{figure}
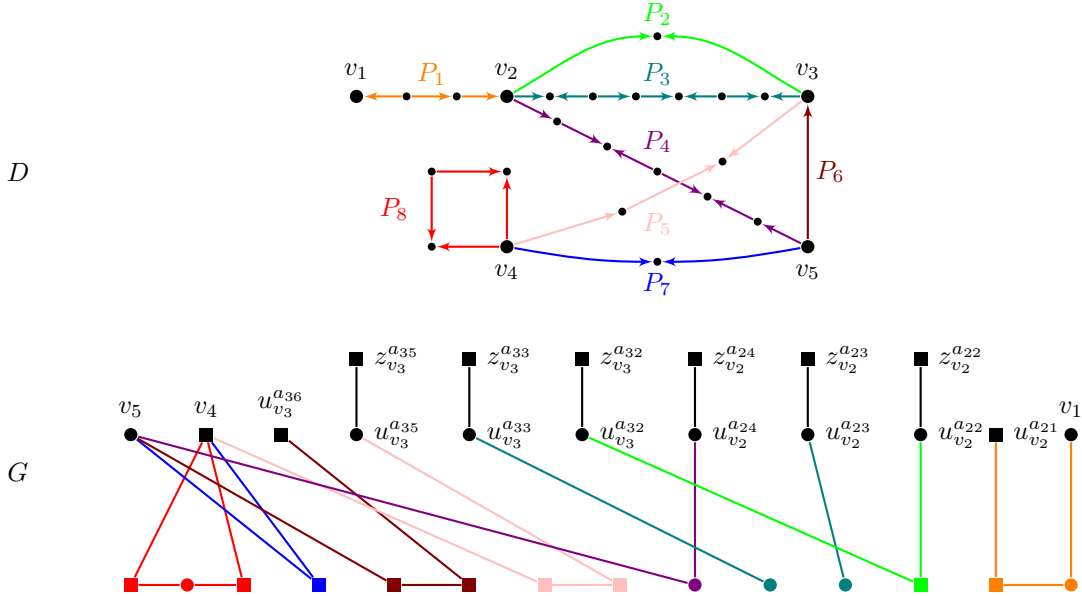
\begin{claim}\label{decmatch}
    $D$ is a positive instance of $(k,1)$-BOGD if and only if $G$ contains a matching covering $Z$.
\end{claim}
    \begin{proofclaim}
        First suppose that $D$ is a positive instance of $(k,1)$-BOGD, so there exists a $(k,1)$-factorization $(F_k,F_1)$ of $D$. We now construct a set $M \subseteq E(G)$. We first let $M$ contain the set of edges $u_b^az_b^a$ for all $b \in B_1$ and $a \in \delta_D^+(b)$. Now let $P \in \mathcal{P}$ and let $a_1,a_2\in A$ be the arcs of $P$ incident to vertices in $B \cup T$ in the order they were used when constructing $G$. We now let $M^P$ be a matching in $G^P$ that covers $Z^P$ and, for $i=1,2$,  contains $e_i^P$ if and only if $a_i \in A(F_1)$. Observe that such a matching exists by the choice of $G^P$ and because the restriction of $(F_k,F_1)$ to $P$ is a $(k,1)$-factorization of $P$. We now add $M^P$ to $M$ and do this for all $P \in \mathcal{P}$. This finishes the description of $M$. We next show that $M$ is a matching that covers $Z$. By construction, every vertex contained in $V(G_P)-\{v_1^P,v_2^P\}$ for some $P \in \mathcal{P}$ is incident to at most one edge in $M$ and $M$ covers $Z^P$ for all $P \in \mathcal{P}$. Next, every $t\in T$ is incident to only one edge in $E(G)$, hence in particular at most one edge in $M$. Next, consider some $b \in B_0$. Every edge of $M$ incident to $b$ corresponds to an arc of $A(F_1)$ incident to $b$. As $F_1$ is a directed matching, there is at most one such arc. Moreover, if $b \in Z$, then $d_D^+(b)=k+1$. As $F_k$ is a $k$-bounded out-galaxy and $(A(F_1),A(F_k))$ is a partition of $A(D)$, we obtain that $b$ is incident to exactly one arc in $A(F_1)$. Hence $b$ is covered by $M$. Finally, consider some $b \in B_1$ and let $a^*$ be the unique arc in $\delta_D^-(b)$. As $(F_k,F_1)$ is a $(k,1)$-out galaxy decomposition of $D$, we obtain that $a^*$ is contained in $A(F_1)$ and all other arcs incident to $b$ in $D$ are contained in $A(F_k)$. Hence the unique edge in $E(G)$ incident to $u_b^{a^*}$ is contained in $M$ and for all $a \in \delta_D^+(b)$, we have that $u_b^az_b^a \in M$ and no other edge incident to $u_b^a$ in $G$ is contained in $M$. It follows that each vertex in $V(G)$ is incident to at most one edge of $M$ and moreover, every vertex in $Z$ is incident to exactly one edge of $M$. Hence $M$ is a matching covering $Z$.    

        Now suppose that $G$ contains a matching $M$ covering $Z$. Consider some $P \in \mathcal{P}$ and let $a_1,a_2$ be the arcs incident to vertices $b_1,b_2 \in B \cup T$, respectively, in the order they were used when constructing $G$. We now let $A_0 \subseteq \{a_1,a_2\}$ be the set that, for $i=1,2$, contains the arc $a_i$ if $e_i^P \in M$. By construction, we obtain that there exists a $(k,1)$-factorization  $(F_k^P,F_1^P)$ of $P$ such that $A(F_1^P)\cap \{a_1,a_2\}=A_0$. We now consider the decomposition $(F_k,F_1)$ of $D$ which is defined by $A(F_k)=\bigcup_{P \in \mathcal{P}}A(F_k^P)$ and $A(F_1)=A(D) -A(F_k)$. Let $K$ be a connected component of $F_k$. If $V(K)\subseteq V(P)-B$ for some $P \in \mathcal{P}$, it follows from the fact that $F_k^P$ is a $k$-bounded out-galaxy that $K$ is a $k$-bounded out-star. We may hence suppose that $V(K)$ contains some $b \in B$. First suppose that $d_K^-(b)\geq 1$, so there exists an arc $a \in \delta_K^-(b)\cap \delta_P^-(b)$ for some $P \in \mathcal{P}$. In particular, we have $b \in B_1$ and by construction, the unique edge in $G$ incident to $u_b^a$ is not contained in $M$. As $u_b^a \in Z$, this contradicts $M$ covering $Z$. We may hence suppose that $d_K^-(b)=0$ and, as $b$ was chosen arbitrarily and $(F_k^p,F_1^P)$ is a $(k,1)$-factorization of $P$ for all $P \in \mathcal{P}$, that $b$ is the only vertex in $V(K)\cap (B \cup T)$ and the underlying graph of $K$ is a star. If $d_D^+(b)\leq k$, we obtain $d_K^+(b)\leq d_D^+(b)\leq k$, so $K$ is a $k$-bounded tree. If $d_D^+(b)= k+1$, then $b \in Z$, so $b$ is incident to an edge of $M$. By construction, this means that $b$ is incident to an arc of $A(F_1)$. It follows that $d_K^+(b)= d_D^+(b)-d_{F_1}^+(b)= (k+1)-1=k$, so $K$ is a $k$-bounded star.

        Now let $K$ be a connected component of $F_1$. It follows from the fact that $F_1^P$ is a directed matching for every $P \in \mathcal{P}$ that $d_{K}(v)\leq 1$ for all $v \in V(G)-(B \cup T)$. For all $b \in B_0 \cup T$, the fact that $b$ is incident to at most one edge in $M$ by construction yields $d_{K}(b)\leq 1$. Finally consider some $b \in B_1$. For all $a \in \delta_D^+(b)$, we obtain from $z_b^a \in Z$ that $u_b^az_b^a\in M$. As $M$ is a matching, we obtain that no other edge incident to $u_b^a$ is contained in $M$. Hence, by construction, we have $d_{K}(b)\leq 1$. This yields that $F_1$ is a directed matching. Hence $(F_k,F_1)$ is a $(k,1)$-factorization of $D$.
    \end{proofclaim}
By Claim~\ref{decmatch}, it suffices to check whether $G$ contains a matching covering $Z$. By Proposition~\ref{checkmatch}, this can be done in polynomial time. Further observe that the algorithm is fully constructive, hence the desired decomposition can be found in polynomial time if it exists.
\end{proof}

\subsection{Decomposing into two (\texorpdfstring{$\geq 2$}{>=2})-bounded out-galaxies}
\label{sec:galaxy_2_2}

In this section, we prove the hardness results contained in Theorem~\ref{bogdmain}. We first describe a gadget which will prove useful.
For every $q\geq 1$, a {\bf $q$-variable gadget} is a digraph $D$ together with a set $S \subseteq V(D)$ of size $q$ with the following properties:
\begin{enumerate}[$(a)$]
    \item $d_D^+(v)=0$ and $d_D^-(v)=1$ hold for all $s \in S$,
    \item for every $(\infty,\infty)$-decomposition $(F,F')$ of $D$, we have $\delta_D^-(S)\subseteq A(F)$ or $\delta_D^-(S)\subseteq A(F')$, and
    \item there exists a $(2,2)$-factorization of $D$.
\end{enumerate}

\begin{figure}[hbt!]
        \begin{center}	
              \begin{tikzpicture}[thick,scale=1, every node/.style={transform shape}]
                \tikzset{vertex/.style = {circle,fill=black,minimum size=5pt, inner sep=0pt}}
                \tikzset{littlevertex/.style = {circle,fill=black,minimum size=4pt, inner sep=0pt}}
                \tikzset{edge/.style = {->,> = latex'}}

                \foreach \i in {1,...,8}{
                    \node[vertex, label=below:$u_\i$] (u\i) at (\i,0) {};
                }
                
                \foreach \i in {1,3,5,7}{
                    \pgfmathtruncatemacro{\j}{\i +1}
                    \draw[dashed,edge,red] (u\i) to (u\j);
                }
                \foreach \i in {2,4,6}{
                    \pgfmathtruncatemacro{\j}{\i +1}
                    \draw[edge,green] (u\i) to (u\j);
                }
                \foreach \i in {1,...,4}{
                    \pgfmathtruncatemacro{\j}{2*\i}
                    \node[vertex, label=above:$v_\i$] (v\i) at (\j,1) {};
                    \draw[edge, green] (u\j) to (v\i);
                }
              \end{tikzpicture}
          \caption{A 4-variable gadget together with a $(2,2)$-factorization. }
          \label{fig:vertex_gadget_star_k_k}
        \end{center}
    \end{figure}
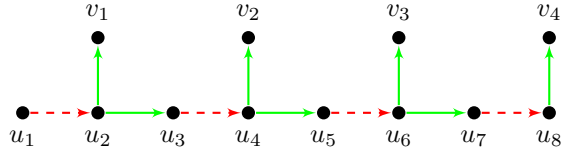

\begin{proposition}\label{gadgex}
    For every integer $q\geq 1$, there exists a $q$-variable gadget whose size is polynomial in $q$.
\end{proposition}
\begin{proof}
  Let $D$ be the digraph on $3q$ vertices $u_1,\dots,u_{2q},v_1,\dots,v_q$ that consists of a directed path on $2q$ vertices $u_1,\dots,u_{2q}$ and the arc $u_{2i}v_i$ for $i\in [q]$, see Figure~\ref{fig:vertex_gadget_star_k_k} for an illustration. Further, let $S=\{v_1,\ldots,v_q\}$. 
  Observe that the size of $D$ is clearly polynomial in $q$ and that $D$ satisfies $(a)$. 
  
  Next, let $(F,F')$ be an $(\infty,\infty)$-factorization of $D$. By symmetry, we may suppose that $u_1u_2 \in A(F)$. As $F$ and $F'$ are out-galaxies, we inductively obtain that $u_{i}u_{i+1}\in A(F)$ for all odd $i \in [2q-1]$ and $u_{i}u_{i+1}\in A(F')$ for all even $i \in [2q-2]$. As $F$ is an out-galaxy, we obtain that $u_{2i}v_i\in A(F')$ for all $i \in [q]$. Hence $(b)$ holds.
  
  We now define $(F,F')$ by $A(F) = \{u_{2i-1}u_{2i} \mid i\in q\}$ and $A(F')=A(D) - A(F)$. It is easy to see that $(F,F')$ is a $(2,2)$-factorization of $D$, so $(c)$ holds. An illustration can be found in Figure \ref{fig:vertex_gadget_star_k_k}.
\end{proof}
We are now ready to prove our NP-completeness results. The proof deals with two cases separately. The first one is the case that the bounds imposed on the out-galaxies are equal and the second one is the case that they are distinct.
\begin{theorem}
    For every integer $k\geq 2$, $(k,k)$-BOGD is NP-complete.
\end{theorem}
\begin{proof}
 First observe that $(k,k)$-BOGD is clearly in NP. We will show that $(k,k)$-BOGD is NP-complete by a reduction from {\sc ME-$(k-1)$-SAT}, which is NP-complete by Proposition~\ref{prop:MEkSAT_NP_hard}.

 Let $(X,\mathcal{C})$ be an instance of {\sc ME-$(k-1)$-SAT}. We now describe a digraph $D$. For every $x \in X$, we let $D$ contain a $q(x)$-variable gadget $(D_x,S_x)$ whose size is polynomial in $q(x)$, where $q(x)$ denotes the number of clauses in $\mathcal{C}$ containing $x$ and $S_x$ contains a vertex $s_{(x,C)}$ for every $C \in \mathcal{C}$ with $x \in C$. Observe that such a gadget exists by Proposition~\ref{gadgex}. Further, for every $C \in \mathcal{C}$, we let $D$ contain an out-star $D_C$ whose root is a vertex $v_C$ and that contains a leaf $t_{(x,C)}$  for every $x \in C$. We finally obtain $D$ by identifying $s_{(x,C)}$ and $t_{(x,C)}$ into a vertex $v_{(x,C)}$ for every $C \in \mathcal{C}$ and every $x \in C$.

 We show in the following that $D$ is a positive instance of $(k,k)$-BOGD if and only if $(X,\mathcal{C})$ is a positive instance of {\sc ME-$(k-1)$-SAT}.

 First suppose that $(X,\mathcal{C})$ is a positive instance of {\sc ME-$(k-1)$-SAT}, so there exists a mapping $\phi \colon X \rightarrow \{\true, \false\}$ such that every $C \in \mathcal{C}$ contains at least $k-1$ variables assigned $\true$ and at least $k-1$ variables assigned $\false$ by $\phi$. We now define a decomposition $(F,F')$ of $D$. For every $x \in X$ with $\phi(x)=\true$ (respectively $\false$), we choose a $(2,2)$-factorization $(F_x,F'_x)$ of $D_x$ such that $\delta_{D_x}^-(S_x)\in A(F_x)$ (respectively $A(F'_x)$). Observe that such a decomposition exists as $D_x$ is a $q(x)$-variable gadget. Next for every $C \in \mathcal{C}$, we define a decomposition $(F_C,F'_C)$ of $D_C$ in the following way: for every $x \in C$, if $\phi(x)=\true$ we let the arc $v_C t_{(x,C)}$ be contained in $A(F'_C)$, otherwise we let it be contained in $A(F_C)$. We now define the decomposition $(F,F')$ of $D$ by $A(F)=\bigcup_{x \in X \cup \mathcal{C}}A(F_x)$ and $A(F')=A(D)- A(F)$. 

 In order to show that $(F,F')$ is a $(k,k)$-factorization of $D$, by symmetry, it suffices to prove that every connected component of $F$ is a $k$-bounded out-star. Let $K$ be a connected component of $F$. It follows by construction that $K$ is a subdigraph of $D_x$ for some $x \in X \cup \mathcal{C}$. If $K$ is a subdigraph of $D_x$ for some $x \in X$, then, it follows by construction that $K$ is a $k$-bounded out-star. If $K$ is a subdigraph of $D_C$ for some $C \in \mathcal{C}$, then it follows by construction that $K$ is an out-star whose root is $v_C$. Further, every leaf of $K$ corresponds to a variable $x \in C$ with $\phi(x)=\false$. As $\phi$ is satisfying assignment for $(X,\mathcal{C})$, there are at most $k$ such variables and hence $K$ is $k$-bounded. It follows that $D$ is a positive instance of $(k,k)$-BOGD.

 Now suppose that $D$ is a positive instance of $(k,k)$-BOGD, so there exists a $(k,k)$-factorization  $(F,F')$ of $D$. For every $x \in X$, as $D_x$ is a $q(x)$-variable gadget, we obtain that $\delta_{D_x}^-(S_x)\subseteq A(F)$ or $\delta_{D_x}^-(S_x)\subseteq A(F')$ holds. We now define $\phi \colon X \rightarrow \{\true, \false\}$ by $\phi(x)=\true$ if $\delta_{D_x}^-(S_x)\subseteq A(F)$ and $\phi(x)=\false$, otherwise. In order to see that $\phi$ is a satisfying assignment for $(X,\mathcal{C})$, consider some $C \in \mathcal{C}$. As $F_k$ is an out-galaxy, for every $x \in C$ with $\phi(x)=\true$, we have that $v_C v_{(x,C)}\in A(F')$. As $F'$ is a $k$-bounded out-galaxy, we obtain that there are at most $k$ variables $x \in C$ with $\phi(x)=\true$. A similar argument shows that there are at most $k$ variables $x \in C$ with $\phi(x)=\false$. Hence $\phi$ is a satisfying assignment for $(X,\mathcal{C})$ and so $(X,\mathcal{C})$ is a positive instance of {\sc ME-$(k-1)$-SAT}.
\end{proof}

    For the case that $k \neq \ell$, we need a slightly more complex clause gadget. More concretely, for some $\ell \in \mathbb{Z}_{\geq 2}$, some $k \in \mathbb{Z}_{\geq \ell +1}\cup \{\infty\}$ and nonnegative integers $\alpha_1,\alpha_2$ with $\alpha_1+\alpha_2=\ell+1$, a {\bf $(k,\ell,\alpha_1,\alpha_2)$-clause gadget} is a digraph $D$ together with two disjoint sets $S_1,S_2\subseteq V(D)$ with $|S_i|=\alpha_i$ for $i \in [2]$ satisfying the following properties:

\begin{enumerate}[$(a)$]   
    \item $d_D^+(v)=0$ and $d_D^-(v)=1$ hold for all $v \in S_1 \cup S_2$,
    \item for every $(k,\ell)$-factorization $(F_k,F_\ell)$ of $D$, we have $\delta_D^-(S_1 \cup S_2)\cap A(F_k)\neq \delta_D^-(S_2)$, and
    \item for every $S_0 \subseteq S_1 \cup S_2$ with $S_0 \neq S_2$, there exists a $(k,\ell)$-factorization of $D$ with $A(F_k)\cap \delta_D^-(S_1 \cup S_2)=\delta_D^-(S_0)$.
\end{enumerate}

\begin{lemma}
Let $\ell \in \mathbb{Z}_{\geq 2}$, $k \in \mathbb{Z}_{\geq \ell +1}\cup \{\infty\}$ and $\alpha_1,\alpha_2$ be nonnegative integers with $\alpha_1+\alpha_2=\ell+1$. There exists a $(k,\ell,\alpha_1,\alpha_2)$-clause gadget.
\end{lemma}
\begin{proof}
    We create a digraph $D$ where $V(D)$ consists of a set $\{r,s_1,\ldots,s_{\alpha_1},s'_1,\ldots,s'_{\alpha_2},u_1,\ldots,u_{\alpha_2}\}$. Further, we let $A(D)$ consist of the arcs $rs_i$ for $i \in [\alpha_1]$ and the arcs $ru_i$ and $u_is'_i$ for $i \in [\alpha_2]$. Finally, we set $S_1=\{s_1,\ldots,s_{\alpha_1}\}$ and $S_2=\{s'_1,\ldots,s'_{\alpha_2}\}$. This finishes the description of $D$. For an illustration, see Figure~\ref{fig:k_l_alpha1_alpha2_clause_gadget}.
\begin{figure}[hbt!]
        \begin{center}	
              \begin{tikzpicture}[thick,scale=1, every node/.style={transform shape}]
                \tikzset{vertex/.style = {circle,fill=black,minimum size=5pt, inner sep=0pt}}
                \tikzset{littlevertex/.style = {circle,fill=black,minimum size=4pt, inner sep=0pt}}
                \tikzset{edge/.style = {->,> = latex'}}

                \node[vertex, label=above:$r$] (r) at (0,0) {};
                \foreach \i in {1,2,3,4}{
                    \node[vertex, label=above:$s_\i$] (s\i) at (-1.5,2.5-\i) {};
                }
                \foreach \i in {1,2,3}{
                    \node[vertex, label=above:$u_\i$] (u\i) at (1.5,2-\i) {};
                    \node[vertex, label=above:$s_\i'$] (sp\i) at (3,2-\i) {};
                }
                \foreach \i in {1,2,3,4}{
                    \draw[edge] (r) to (s\i);
                }
                \foreach \i in {1,2,3}{
                    \draw[edge] (r) to (u\i);
                    \draw[edge] (u\i) to (sp\i);
                }
                \draw[dotted] (-1.5,0.15) ellipse (0.5 and 2);
                \node[] (S1) at (-2.3, 0.15) {$S_1$};
                
                \draw[dotted] (3,0.15) ellipse (0.5 and 1.6);
                \node[] (S1) at (3.8, 0.15) {$S_2$};
              \end{tikzpicture}
          \caption{A $(k,6,4,3)$-clause gadget for some $k \geq 7$.}
          \label{fig:k_l_alpha1_alpha2_clause_gadget}
        \end{center}
    \end{figure}
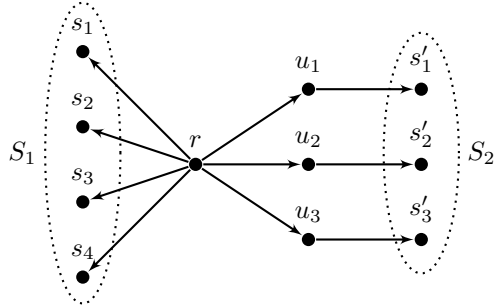
    
    Clearly, $S_1$ and  $S_2$ are disjoint, we have $|S_i|=\alpha_i$ for $i \in [2]$, and $(a)$ is satisfied. For $(b)$, suppose for the sake of a contradiction that there exists a $(k,\ell)$-factorization $(F_k,F_\ell)$ of $D$ with $A(F_k)\cap \delta_D^-(S_1 \cup S_2)=\delta_D^-(S_2)$. As $F_k$ is an out-galaxy, we obtain that $\delta_D^+(r)\subseteq F_\ell$. As $d_D^+(r)=\alpha_1+\alpha_2=\ell+1$, this contradicts $F_\ell$ being an $\ell$-bounded out-galaxy. This proves $(b)$.

    Now consider some $S_0 \subseteq S_1 \cup S_2$ with $S_0\neq S_2$. We define the factorization $(F_k,F_\ell)$ of $D$ by $A(F_k)=\{rs_i \mid s_i \in S_0 \cap S_1\} \cup \{u_is'_i \mid s'_i \in S_0 \cap S_2\}\cup \{ru_i \mid s'_i \in S_2 - S_0\}$ and $A(F_\ell)=A(D)-A(F_k)$. It is easy to see that $(F_k,F_\ell)$ has the desired properties. This proves $(c)$.
\end{proof}

We are now ready to prove the second part of our hardness results.
\begin{theorem}
    For every integer $\ell \geq 2$ and every $k\in  \mathbb{Z}_{\geq \ell+1} \cup \{ \infty \}$,  $(k,\ell)$-BOGD is NP-complete.
\end{theorem}
\begin{proof}
    Let us fix $\ell \in \mathbb{Z}_{\geq 2}$ and $k\in  \mathbb{Z}_{\geq \ell+1} \cup \{ \infty \}$.
    First observe that $(k,\ell)$-BOGD is clearly in NP. We will show that $(k,\ell)$-BOGD is NP-complete by a reduction from $(\ell+1)$-SAT, which is known to be NP-complete since $\ell+1 \geq 3$ and by Proposition \ref{prop:3SAT}.

    Let $(X,\mathcal{C})$ be an instance of $(\ell+1)$-SAT. We now create a digraph $D$. For every $x \in X$, we let $D$ contain a $q(x)$-variable gadget $(D^x,S^x)$ where $q(x)$ is the total number of occurences of $x$ or $\bar{x}$ in $(X,\mathcal{C})$ and $S^x$ contains a vertex $t_{(x,C)}$ for every $C \in \mathcal{C}$ containing $x$ or $\bar{x}$. Now consider some clause $C \in \mathcal{C}$, let $\alpha_1$ be the number of positive literals in $C$ and let $\alpha_2$ be the number of negative literals in $C$. Observe that $\alpha_1+\alpha_2=\ell+1$ by definition. We now let $D$ contain a $(k,\ell,\alpha_1,\alpha_2)$-clause gadget $(D^C,S_1^C,S_2^C)$ where $S_1^C$ contains a vertex $s_{(x,C)}$ for every $x \in X$ with $x \in C$ and $S_2^C$ contains a vertex $s_{(x,C)}$ for every $x \in X$ with $\bar{x} \in C$. We now obtain $D$ by identifying $t_{(x,C)}$ and $s_{(x,C)}$ for all $x \in X$ and $C \in \mathcal{C}$ with $x \in C$ or $\bar{x} \in C$. We show in the following that $D$ is a positive instance of $(k,\ell)$-BOGD if and only if $(X,\mathcal{C})$ is a positive instance of $(\ell+1)$-SAT.

    First suppose that $(X,\mathcal{C})$ is a positive instance of $(\ell+1)$-SAT, so there exists a satisfying assignment $\phi \colon X \rightarrow \{\true,\false\}$ for $(X,\mathcal{C})$. For every $x \in X$ with $\phi(x)=\true$ ($\phi(x)=\false$), we choose a $(k,\ell)$-factorization $(F_k^x,F_\ell^x)$ of $D^x$ with $\delta_{D^x}^-(S^x)\subseteq A(F_\ell^x)$ ($\delta_{D^x}^-(S^x)\subseteq A(F_k^x)$, respectively). For every $C \in \mathcal{C}$, let $S_0^C$ be the set containing the vertices $s_{(x,C)}$ for all variables $x$ for which either $x \in C$ and $\phi(x)=\true$ or $\bar{x}\in C$ and $\phi(x)=\false$ holds. As $\phi$ is satisfying, we have $S_0^C\neq S_2^C$. We now choose a $(k,\ell)$-factorization $(F_k^C,F_\ell^C)$ of $D^C$ with $A(F_k^C)\cap \delta_{D^C}^-(S_1^C \cup S_2^C)=\delta_{D^C}^-(S_0^C)$. Finally, we define $(F_k,F_\ell)$ by $A(F_k)=\bigcup_{x \in X \cup \mathcal{C}}A(F_k^x)$ and $A(F_\ell)=A(D)-A(F_k)$. Observe that every connected component of $F_k$ or $F_\ell$ is fully contained in $D^x$ for some $x \in X \cup \mathcal{C}$. It follows that $(F_k,F_\ell)$ is a $(k,\ell)$-factorization of $D$, so $D$ is a positive instance of BOGD.

    Now suppose that $D$ is a positive instance of BOGD, so there exists a $(k,\ell)$-factorization $(F_k,F_\ell)$ of $D$. Consider some $x \in X$. As $D^x$ is a $q(x)$-variable gadget, we obtain that $\delta_{D^x}^-(S^x)\subseteq A(F_\ell)$ or $\delta_{D^x}^-(S^x)\subseteq A(F_k)$ holds. We now define $\phi \colon X \rightarrow \{\true,\false\}$ by $\phi(x)=\true$ if $\delta_{D^x}^-(S^x)\subseteq A(F_\ell)$ and $\phi(x)=\false$, otherwise. In order to see that $\phi$ is a satisfying assignment for $(X,\mathcal{C})$, consider some $C \in \mathcal{C}$ and suppose for the sake of a contradiction that $C$ is not satisfied by $\phi$. For all $x \in X$ with $x \in C$, we obtain that the unique arc entering $s_{(x,C)}$ which is not contained in $A(D^C)$ is contained in $A(F_k)$. Hence, as $F_k$ is an out-galaxy, we obtain that $\delta_{D^C}^-(S_1^C)\subseteq A(F_\ell)$. For all $x \in X$ with $\bar{x} \in C$, we obtain that the unique arc entering $s_{(x,C)}$ which is not contained in $A(D^C)$ is contained in $A(F_\ell)$. Hence, as $F_\ell$ is an out-galaxy, we obtain that $\delta_{D^C}^-(S_2^C)\subseteq A(F_k)$. This yields $\delta_{D^C}^-(S_1 \cup S_2)\cap A(F_k)= \delta_{D^C}^-(S_2)$, a contradiction as $D^C$ is a clause gadget. Hence $\phi$ satisfies $(X,\mathcal{C})$ and so $(X,\mathcal{C})$ is a positive instance of $(\ell+1)$-SAT.
\end{proof}

\section*{Acknowledgement}
This work originates from the workshop on digraphs organized in May 2023 in Sète. We thank the organizers for this opportunity.

\bibliography{refs}

\end{document}